\documentclass{amsart}

\usepackage{amsmath,amssymb,amsfonts,amsthm}
\usepackage[ragged]{sidecap}
\usepackage{subfigure}
\usepackage[mathscr]{eucal}
\usepackage{url}
\usepackage{lscape}
\usepackage{graphicx,epsfig,amscd,graphics,xypic}
\usepackage{mathrsfs}
\usepackage{psfrag}
\usepackage[active]{srcltx}
\usepackage{amssymb, amsmath,array, amscd}
\usepackage{enumerate}
\usepackage[colorlinks]{hyperref}

\def\ssigma{ \mbox{\Large{${\sigma}_{\hspace{-0.01cm} 1 }$}\hspace{-0.05cm}}}

\theoremstyle{change}
\newtheorem{thm}{Theorem}[section]
\newtheorem{THM}{Theorem}
\newtheorem{COR}{Corollary}
\newtheorem{prop}{Proposition}[section]
\newtheorem{lemma}{Lemma}[section]
\newtheorem{claim}{Claim}[section]

\newtheorem{remark}{Remark}[section]
\newtheorem{cor}{Corollary}[section]

\newtheorem{example}{Example}[section]

\begin{document}

\sloppy

\title[Exceptional CDQL Webs]{The Classification of Exceptional CDQL Webs\\ on Compact Complex Surfaces}

\author{J.  V.  Pereira}
\address{Instituto de Matem\'atica Pura e Aplicada\\ Est.
D. Castorina, 110\\
22460-320, Rio de Janeiro, RJ, Brazil} \email{ jvp@impa.br}

\author{L.  Pirio}
\address{IRMAR \\ Campus de Beaulieu \\ 35042 Rennes Cedex, France}
\email { luc.pirio@univ-rennes1.fr  }

\maketitle

\date{\today}

\setcounter{tocdepth}{1}

\begin{abstract}
Codimension one webs are configurations of finitely many codimension one foliations in general position.
Much of the classical theory evolved around the concept of   abelian relation:
a functional relation  among the first integrals
of the foliations defining the web reminiscent of Abel's addition theorem
in classical algebraic geometry.
The abelian relations of a given web form a finite dimensional vector space  with
dimension (the rank of the web)   bounded by  Castelnuovo number $\pi(n,k)$ where $n$ is the dimension of the ambient  space
and $k$ is the number of foliations defining the web.  A fundamental  problem in  web geometry is the classification of
exceptional webs, that is, webs of maximal rank  not equivalent to the dual of a projective curve. Recently,
J.-M. Tr\'{e}preau   proved that there are no exceptional $k$-webs for  $n \ge 3$ and $k \ge 2n$. In dimension two there
are examples for arbitrary $k$  and the classification problem is wide open.

In this paper, we classify the exceptional {{C}ompletely {D}ecomposable {Q}uasi-{L}inear (CDQL)}   webs
globally defined on compact complex surfaces.
By definition, the CDQL  $(k+1)$-webs are   formed by the superposition of   $k$ linear foliations
and one non-linear foliation. For instance, we show that up to projective transformations there are exactly four countable families and thirteen sporadic
exceptional CDQL webs on $\mathbb P^2$.
\end{abstract}

\tableofcontents

\pagebreak


\section{Introduction and statement of the main results}

\subsection{Codimension one webs of maximal rank}

A {\it germ of regular $k$-web} $\mathcal W=\mathcal F_1 \boxtimes
\cdots \boxtimes \mathcal F_k$ of codimension one on $(\mathbb
C^n,0)$ is a collection of $k$ germs of smooth holomorphic
foliations $\mathcal F_i$ with tangent spaces  in general position at the origin.
By definition, the
$\mathcal F_i$'s are the {\it defining foliations} of $\mathcal W$.
If they are respectively induced by differentials 1-forms $\omega_1,\ldots,\omega_k$, then  the {\it space of abelian relations} of ${\mathcal W}$ is the
vector space
\begin{align*}
\mathcal A(\mathcal W)= \Big\lbrace {\big(\eta_i\big)}_{i=1}^k \in
{\Omega^1(\mathbb C^n,0)}^k \, \, \Big| \;
\sum_{i=1}^k\eta_i = 0\, \, \text{ and } \, \,
\,\forall i \, \, d\eta_i
=0 \,, \; \eta_i\wedge \omega_i=0 \,  \Big\rbrace \, .
\end{align*}

If $u_i:(\mathbb C^n,0) \to (\mathbb C,0)$ are local submersions
defining the foliations $\mathcal F_i$ then, after integration, the
abelian relations can be read as  functional equations of the form
$\sum_{i=1}^k g_i(u_i) =0$ for suitable germs of holomorphic functions $g_i: (\mathbb C,0) \to (\mathbb C,0)$.

The dimension of $\mathcal A(\mathcal W)$ is commonly called the
{\it rank} of $\mathcal W$ and denoted by ${rk}(\mathcal W)$.
It is a theorem of Bol (for $n=2$) and Chern (for $n\ge3$) that
\begin{equation}\label{E:bounds}
 {rk}( \mathcal W) \le \pi(n,k) = \sum_{j=1}^{\infty}
\max \big(0 , k - j(n-1) - 1\big) .
\end{equation}

A $k$-web $\mathcal W$ on $(\mathbb C^n,0)$ is
 of {\it maximal rank} if ${rk}(\mathcal W) =
\pi(n,k)$. The integer $\pi(n,k)$ is the well-known Castelnuovo's
bound for the arithmetic genus of irreducible and non-degenerated
algebraic curves of degree $k$ on $\mathbb P^n$.

One of the main topics of the theory of webs concerns the
characterization of webs of  maximal rank.  It follows from Abel's
Addition Theorem that all the webs $\mathcal W_C$ obtained from
reduced Castelnuovo curves\footnote{{That is, non-degenerate  algebraic curve $C\subset {\mathbb P}^n$ such that $p_a(C)=\pi(n,k)$ where $k=deg(C)$.}} $C$ by projective duality  are of maximal
rank (see \cite{jvpBourbaki} for instance). The webs analytically
equivalent to $\mathcal W_C$ for some non-degenerated projective curve $C$ are the so called {\it algebraizable webs}.

It can be traced back to Lie the proof that  all $4$-webs on
$(\mathbb C^2,0)$ of maximal rank are algebraizable. In \cite{Bol},
Bol proved that a maximal rank $k$-web on $(\mathbb C^3,0)$ is
algebraizable when $k\ge 6$. Recently, building up on previous work
by  Chern and Griffiths \cite{Jbr}, Tr\'{e}preau extended Bol's result and
established in \cite{Trepreau} that $k$-webs  of maximal rank on $(\mathbb C^n,0)$ are algebraizable whenever $n\ge3$ and $k\ge 2n$.

The non-algebraizable webs of maximal rank on $(\mathbb C^2,0)$ are
nowadays called {\it exceptional webs}. For almost 70 years there
was just one  example, due to Bol \cite{Bolcontra},  of exceptional planar web in the
literature. Recently a number of new examples have appeared, see
\cite{PTese,Robert,PT,MPP}. Despite these  new examples,
the classification problem for exceptional planar  webs is wide open.

\subsection{Characterization of planar webs of  maximal rank}
Although a classification  seems out of reach,
there are  methods to decide if a given  web has maximal rank.
The first result in this direction is   due to Pantazi \cite{Pantazi}. It was  published during the second
world war and  remained unknown  to the practitioners of web theory until recently, see \cite{PTese}.
Unaware of this classical result, H\'{e}naut  \cite{Henaut} worked out an alternative approach to determine if a
given web has maximal rank. Both approaches share in common the use of prolongations
of differential systems to express  the  maximality of the rank by the
vanishing  of certain differential expressions determined  by  the defining equations of the web.

It has to be noted that these results are wide generalizations of the classical criterion
of Blaschke-Dubourdieu for the maximality of the rank of $3$-webs.  If $\mathcal W=\mathcal F_1 \boxtimes \mathcal F_2
\boxtimes \mathcal F_3$ is a planar $3$-web and the foliations $\mathcal
F_i$ are defined by $1$-forms $\omega_i$ satisfying $\omega_1
+\omega_2 + \omega_3=0$ then a simple computation ensures the existence of
a unique $1$-form $\gamma$ such that $d\omega_i = \gamma
\wedge \omega_i$ for $i=1,2,3$. Although  $\gamma$ does
depend on the choice of the $1$-forms $\omega_i$ its differential $d\gamma$ is
intrinsically attached to $\mathcal W$. It is the so called {\it
curvature} $K(\mathcal W)$ of $\mathcal W$. In \cite{BD} it is proved that
a $3$-web $\mathcal W$ has  maximal rank  if and only if $K(\mathcal W)=0$.

Building on Pantazi's result, Mih\u{a}ileanu gave in \cite{Mi} a necessary condition for  a planar $k$-web
be of maximal rank : if  $\mathcal W$ has maximal rank then $K(\mathcal W)=0$. Now, the curvature
$K(\mathcal W)$ is the sum of the curvatures of all $3$-subwebs of $\mathcal W$. Recently H\'{e}naut, Ripoll and Robert
 (see \cite[p.281]{Henaut2},\cite{Ri}) have rediscovered Mih\u{a}ileanu necessary condition using
H\'{e}naut's approach.

As in the case of $3$-webs  the curvature is a holomorphic  $2$-form intrinsically attached to $\mathcal W$: it does not
depend on the choice of the defining equations of ${\mathcal W}$. Another nice feature of the curvature is that  it still makes sense, as a meromorphic $2$-form,  for global webs.
More precisely if $S$ is complex surface  then  a  {\it global $k$-web} on $S$ can be defined
as an element $\mathcal W= [\omega]$ of $\mathbb P H^0( S,
 \mathrm{Sym}^k \Omega^1_S\otimes\mathcal N)$
--- where $\mathcal N $ denotes a line-bundle and
$\mathrm{Sym}^k \Omega^1_S$ the sheaf of $k$-symmetric powers of
differential $1$-forms  on $S$ --- subjected to the following two conditions:
{\it (i)} the zero locus of $\omega$ has codimension at least two; {\it (ii)}
 $\omega(p)$ factors as the product of
pairwise linearly independent $1$-forms at some point $p \in S$.
For $k=1$ the condition {\it (ii)} is vacuous and we recover one of the usual definitions of foliations. When $k\ge 2$, the set where the condition
{\it (ii)} does not hold is the {\it discriminant} of $\mathcal W$ and
will be denoted by $\Delta(\mathcal W)$. For $k \ge 3$, the curvature  $K(\mathcal W)$
is a  global meromorphic $2$-form on $S$ with polar set contained in $\Delta(\mathcal W)$.

Elementary  arguments imply that the space of abelian relations
of $\mathcal W$, in this global setup, is a local system
over $S \setminus \Delta(\mathcal W)$, see  for instance \cite[Th\'eor\`eme 1.2.2]{PTese}. The rank of $\mathcal W$ appears now as
the rank of the local system ${\mathcal A}({\mathcal W})$.

One has to be careful when talking about defining foliations of a global web since
these will make sense only in  sufficiently small analytic  open subsets
 of $S$. When it is possible to write globally $\mathcal
W= \mathcal F_1 \boxtimes \cdots \boxtimes \mathcal F_k$ we will say
that $\mathcal W$ is  {\it completely decomposable}.

When $S$ is a
pseudo-parallelizable
surface\footnote{A complex manifold $M$ of dimension $n$ is called {\it pseudo-parallelizable} if it carries
$n$ global meromorphic $1$-forms   $\omega_1, \ldots, \omega_n$  with exterior product  $\omega_1 \wedge \cdots \wedge \omega_n$ not identically zero. },
  a  global $k$-web on $S$  can be
alternatively defined as an element $\mathcal W= [\omega]$ of the projective space
$\mathbb P_{\mathbb C(S)} ( \mathrm{Sym}^k \Omega^1_S)$
--- where $\mathbb C(S)$ is the field of meromorphic functions on $S$ and
$\mathrm{Sym}^k \Omega^1_S$ denotes  now  the $\mathbb C(S)$-vector
space of meromorphic $k$-symmetric powers of differential $1$-forms  on
$S$--- subjected to the condition that
$\omega$ factors as the product of  pairwise linearly independent
$1$-forms at some point of $S$.

\subsection{Mih\u{a}ileanu necessary condition and  $\mathcal F$-barycenters}
The present work stems from an attempt to understand geometrically Mih\u{a}ileanu's necessary
condition for the maximality of the rank. More precisely we try to understand
the conditions imposed by the vanishing of the curvature on the behavior of $\mathcal W$
over its discriminant. It has to be mentioned that the idea of analyzing  webs through
theirs discriminants is not new, see  \cite{Cerveau} and \cite{LinsNeto}.
More recently,  \cite{Henaut2} advocates the
study of webs (decomposable or not)  in  neighborhoods of theirs discriminants.

Our result in this direction  is stated in terms of $\beta_{\mathcal F}( \mathcal W)$
--- the {\it $\mathcal F$-barycenter} of a web $\mathcal W$. Suppose that $S$ is a pseudo-parallelizable surface and
 $\mathcal F \in \mathbb P_{\mathbb C(S)}(\Omega^1_S)$ is a
foliation on it. There is a naturally  defined  affine structure on
$\mathbb A^1_{\mathcal F} = \mathbb P_{\mathbb C(S)}(\Omega^1_S)
\setminus \mathcal F$.
  If $\mathcal W \in \mathbb P_{\mathbb
C(S)}(\mathrm{Sym}^k\Omega^1_S)$  is a $k$-web not containing
$\mathcal F$ as one of its defining foliations then it can be
loosely interpreted as $k$ points in $\mathbb A^1_{\mathcal F}$.
The $\mathcal F$-barycenter of $\mathcal W$ is then the foliation
$\beta_{\mathcal F}(\mathcal W)$ defined by  the
barycenter of these $k$ points in $\mathbb A^1_{\mathcal F}$. For a precise definition
and some properties of $\beta_\mathcal F(\mathcal W)$, see
Sections \ref{S:bary} and \ref{S:fbary}.

\begin{THM}\label{TT:curvatura}
Let $\mathcal F$ be a foliation and  $\mathcal W= \mathcal F_1
\boxtimes \mathcal F_2 \boxtimes \cdots \boxtimes \mathcal F_k$ be a
$k$-web, $k \ge 2$, both defined on the same domain $U \subset
\mathbb C^2$. Suppose that $C$ is an irreducible component of
$\mathrm{tang}(\mathcal F, \mathcal F_1)$ that
  is not contained in $\Delta(\mathcal W)$. The curvature
 $K(\mathcal F \boxtimes \mathcal W)$ is holomorphic over a generic point of $C$ if and only if the curve
  $C$ is $\mathcal F$-invariant or $\beta_{\mathcal F}(
\mathcal W')$-invariant, where $\mathcal W'=\mathcal F_2 \boxtimes
\cdots \boxtimes \mathcal F_k$.
\end{THM}

Theorem \ref{TT:curvatura} is the cornerstone  of our approach to the  classification of
exceptional completely decomposable quasi-linear webs (CDQL webs for short) on compact complex surfaces.

\subsection{Linear webs and CDQL webs}
Linear webs are classically defined as the ones for which all the leaves  are open subsets of lines. Here we will adopt the following global definition.
A web $\mathcal W$ on compact complex surface $S$ is { \it linear} if (a) the universal covering of $S$ is an open subset $\widetilde S$
of $\mathbb P^2$; (b) the group of deck transformations acts on $\widetilde S$ by automorphisms of $\mathbb P^2$, and; (c) the pull-back of
$\mathcal W$ to  $\widetilde S$ is linear in the classical sense\begin{footnote}{Alternatively one  could assume that
$S$ admits a $(\mathbb P^2,\mathrm{PGL}(3,\mathbb C))$-structure and that $\mathcal W$ is linear in  the  local charts of this structure.
Although more general,  this definition does not seem to encompass more examples of linear webs.  To avoid a lengthy  case by case analysis of the classification of
$(\mathbb P^2,\mathrm{PGL}(3,\mathbb C))$-structures
 \cite{Kl} we opted for  the more astringent definition above.}\end{footnote}.

A {\it CDQL $(k+1)$-web} on a compact complex surface $S$ is, by definition, the superposition of $k$ linear foliations and one non-linear foliation.

It follows from  \cite{IKO, KO}  that the only compact complex surfaces satisfying (a) and (b) are:
 the projective plane;  surfaces covered by the unit ball; Kodaira primary surfaces;
 complex tori;  Inoue surfaces; Hopf surfaces  and principal elliptic bundles over hyperbolic curves with odd first Betti number.

If $S$ is not $\mathbb P^2$ then the group of deck transformations is infinite. Because it  acts on $\widetilde S$ without fixed points, every linear foliation on $S$ is a smooth foliation.
An  inspection of Brunella's classification of smooth foliations \cite{Brsmooth} reveals that the only compact complex surfaces
admitting at least two distinct linear foliations are the projective plane, the complex tori and the Hopf surfaces. Moreover
the only Hopf surfaces admitting  four distinct linear foliations  are the primary Hopf surfaces  $H_{\alpha}$,  $|\alpha| > 1$. Here $
  H_{\alpha}$ is the quotient of $\mathbb C^2 \setminus \{ 0\}$ by the map $(x,y) \mapsto (\alpha x, \alpha y)$.

The linear foliations on  complex tori are pencils of parallel lines on theirs universal coverings. The ones
on Hopf surfaces are either pencils of parallels lines or the pencil of lines through the origin of $\mathbb C^2$. In particular
all completely decomposable linear webs on compact complex surfaces are algebraic\begin{footnote}{Beware that algebraic here means that they are locally
dual to plane curves. In the cases under scrutiny  they are dual to certain products of lines.   }\end{footnote} on theirs universal coverings.

\subsection{Classification of exceptional CDQL  webs on the projective plane}
\begin{figure}[h]
\begin{minipage}{4.0cm}
\begin{center}
\includegraphics[width=4.0cm,height=4.0cm]{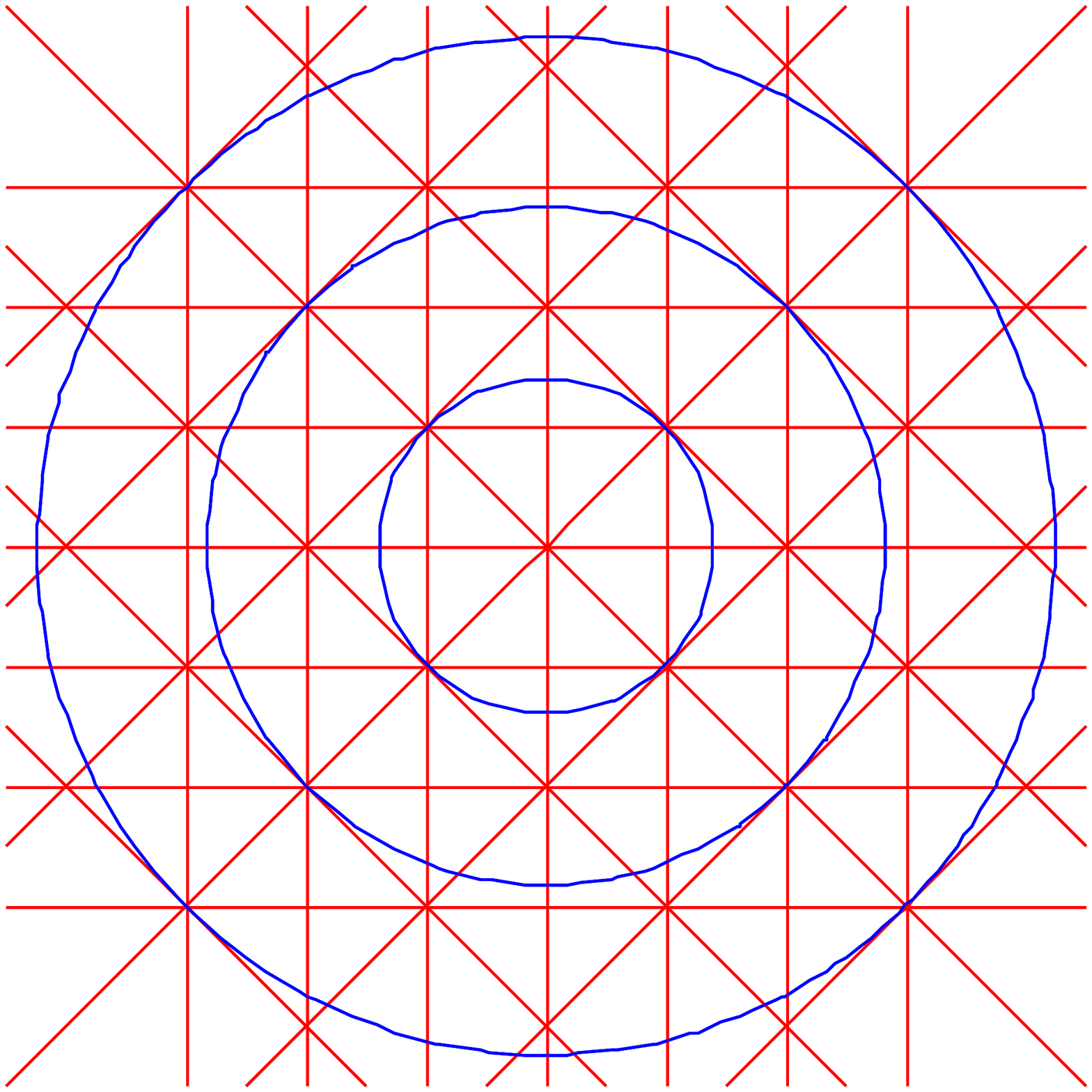}
\end{center}
\end{minipage}
\hfill
\begin{minipage}{4.0cm}
\begin{center}
\includegraphics[width=4.0cm,height=4.0cm]{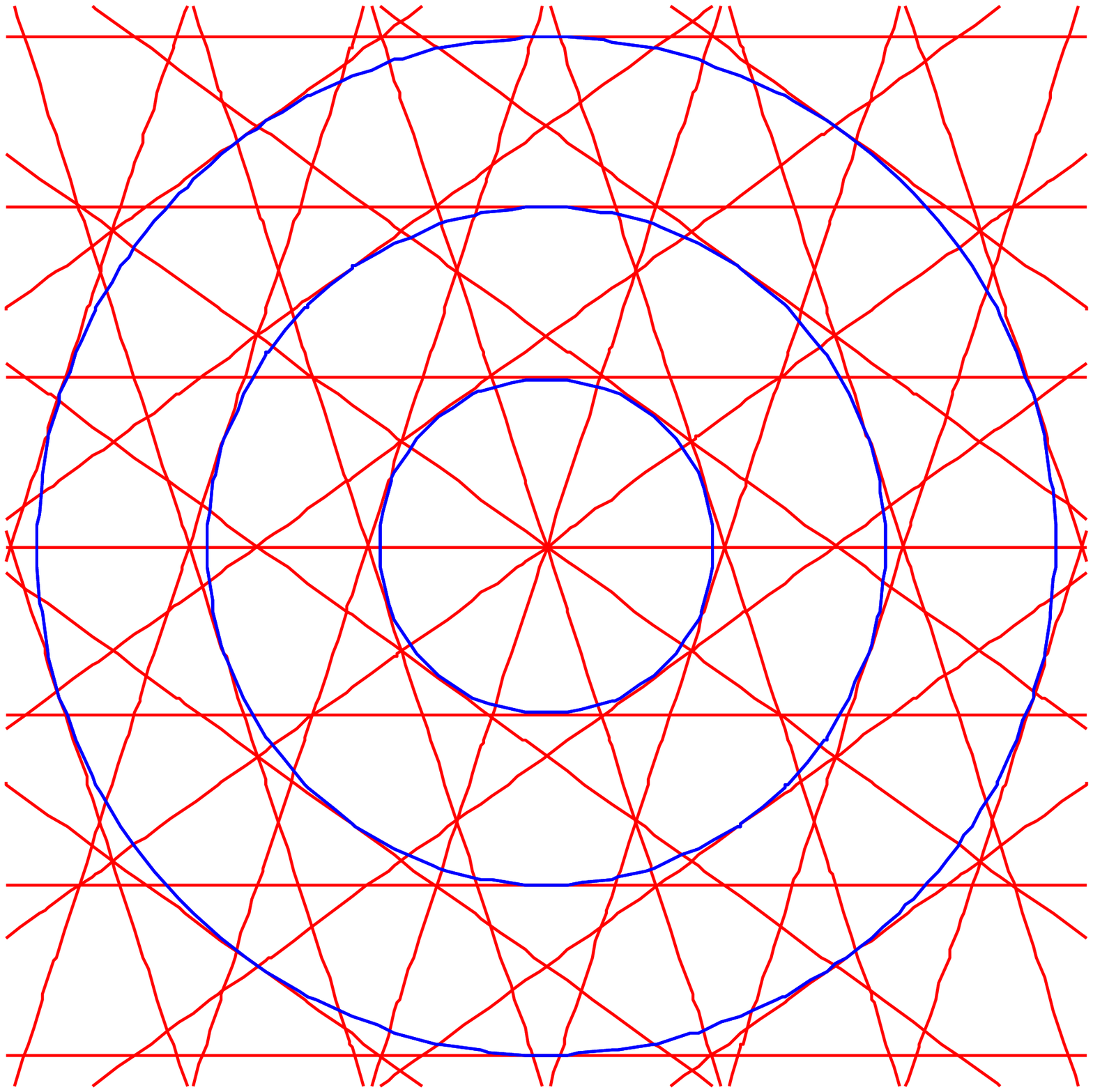}
\end{center}
\end{minipage}
\hfill
\begin{minipage}{4.0cm}
\begin{center}
\includegraphics[width=4.0cm,height=4.0cm]{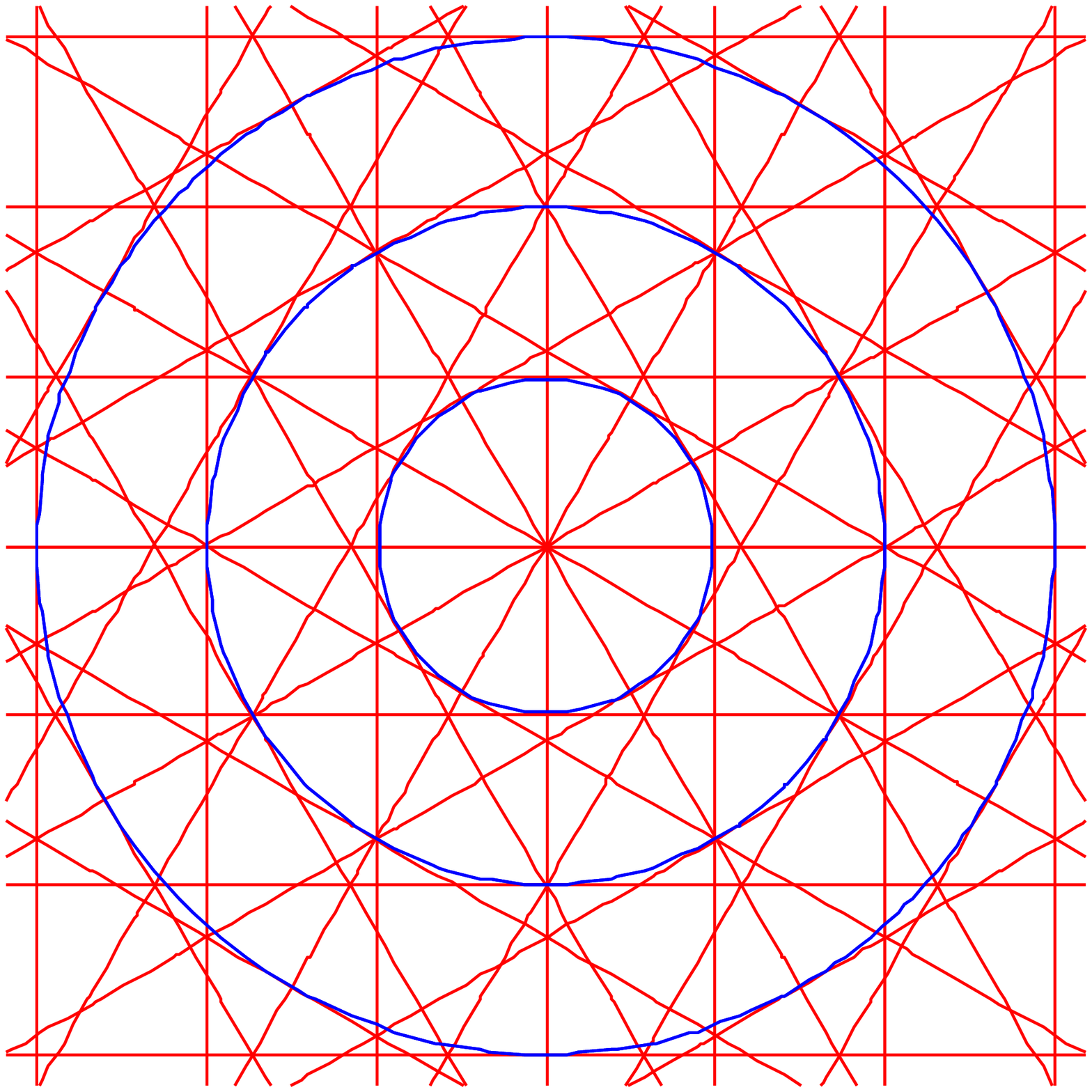}
\end{center}
\end{minipage}
\hfill
\begin{minipage}{4.0cm}
\begin{center}
\includegraphics[width=4.0cm,height=4.0cm]{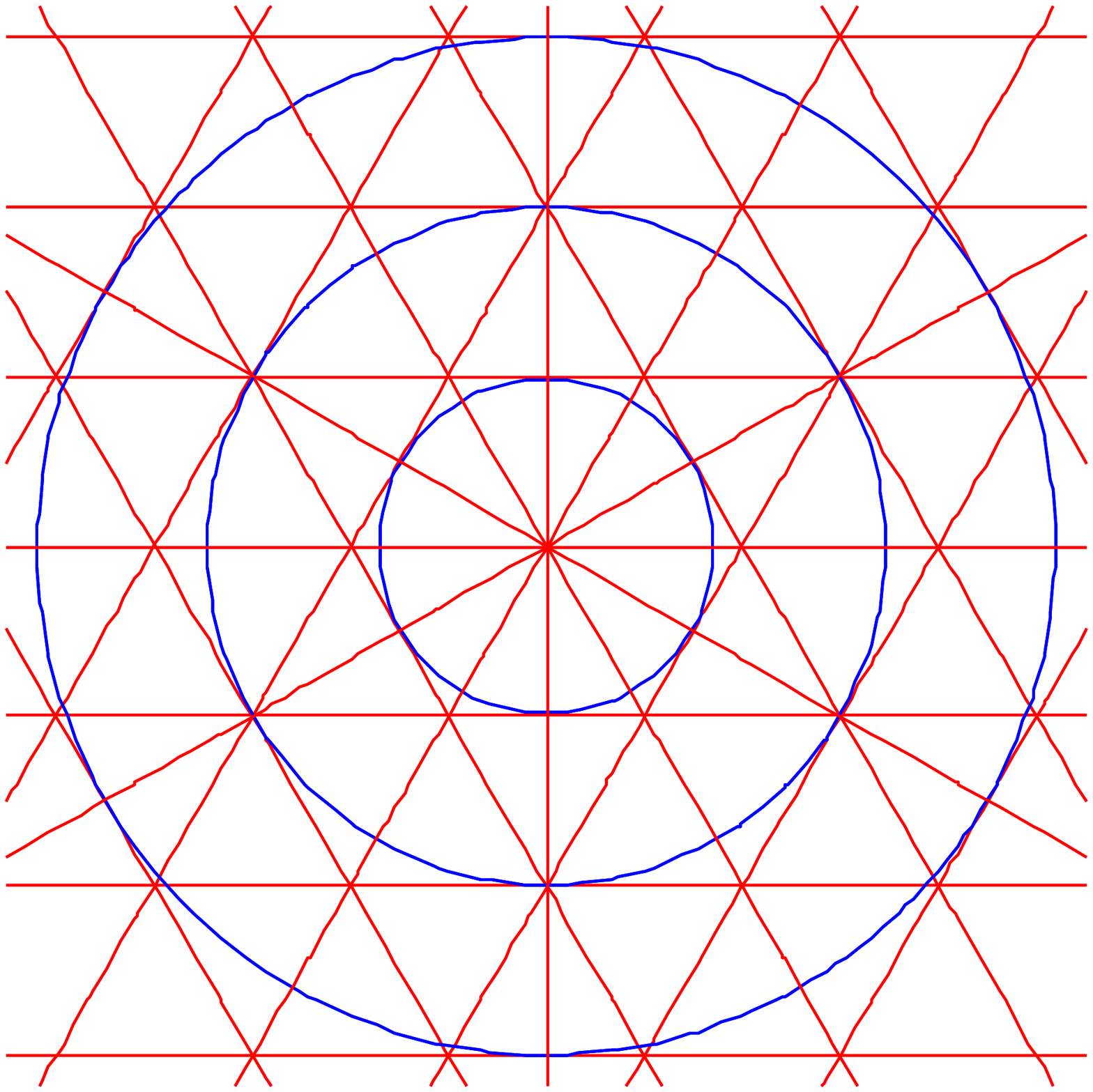}
\end{center}
\end{minipage}
\hfill
\begin{minipage}{4.0cm}
\begin{center}
\includegraphics[width=4.0cm,height=4.0cm]{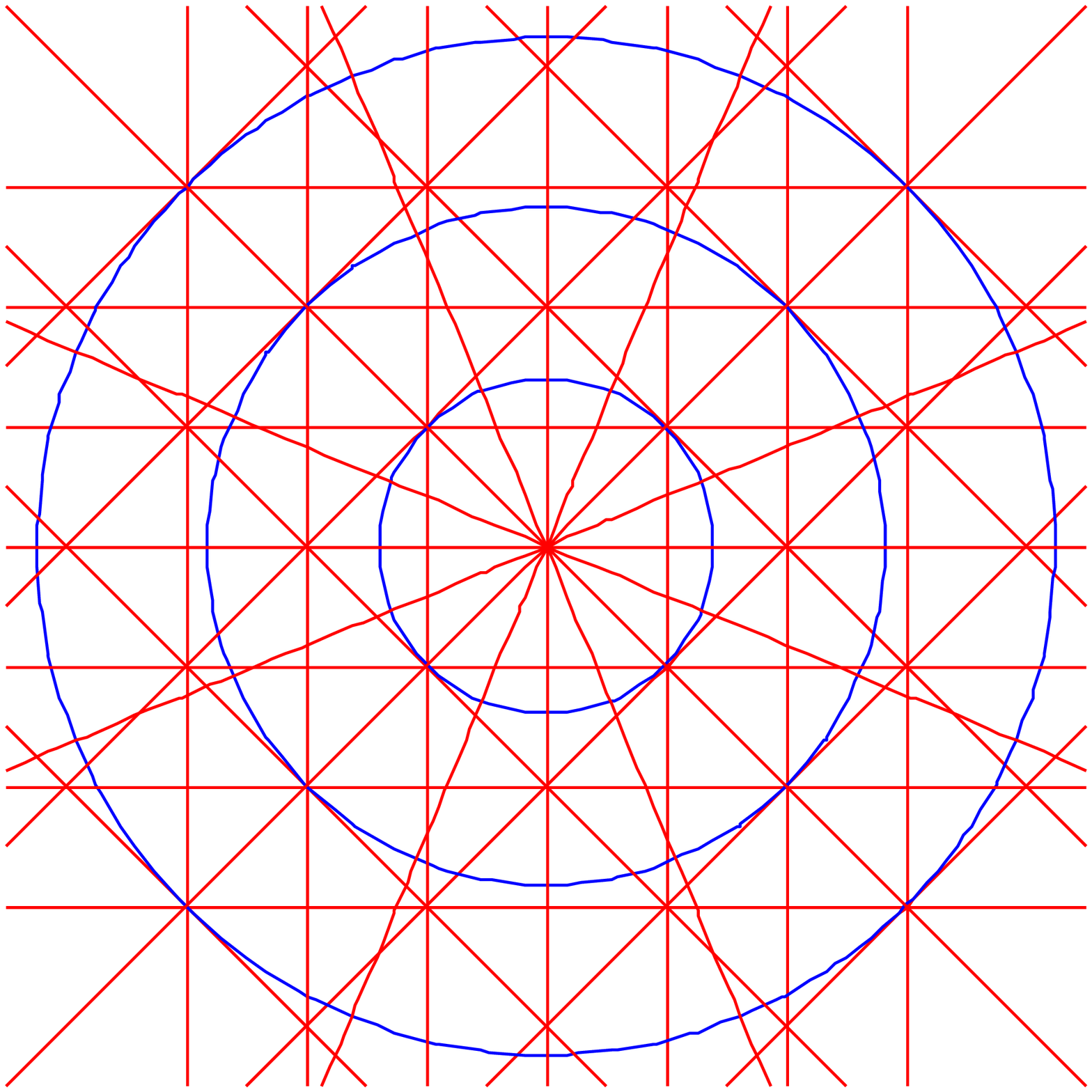}
\end{center}
\end{minipage}
\hfill
\begin{minipage}{4.0cm}
\begin{center}
\includegraphics[width=4.0cm,height=4.0cm]{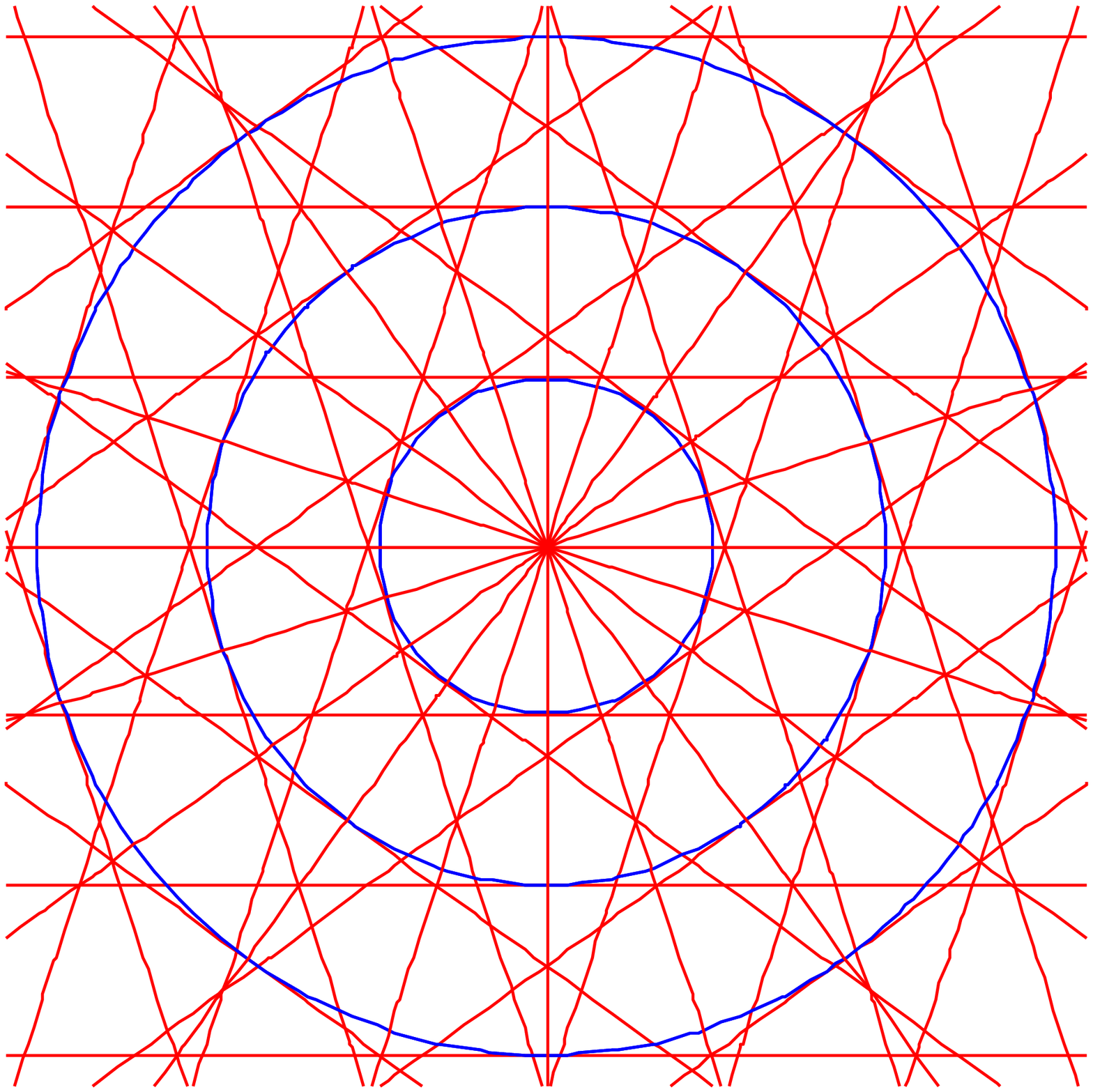}
\end{center}
\end{minipage}
\hfill
\begin{minipage}{4.0cm}
\begin{center}
\includegraphics[width=4.0cm,height=4.0cm]{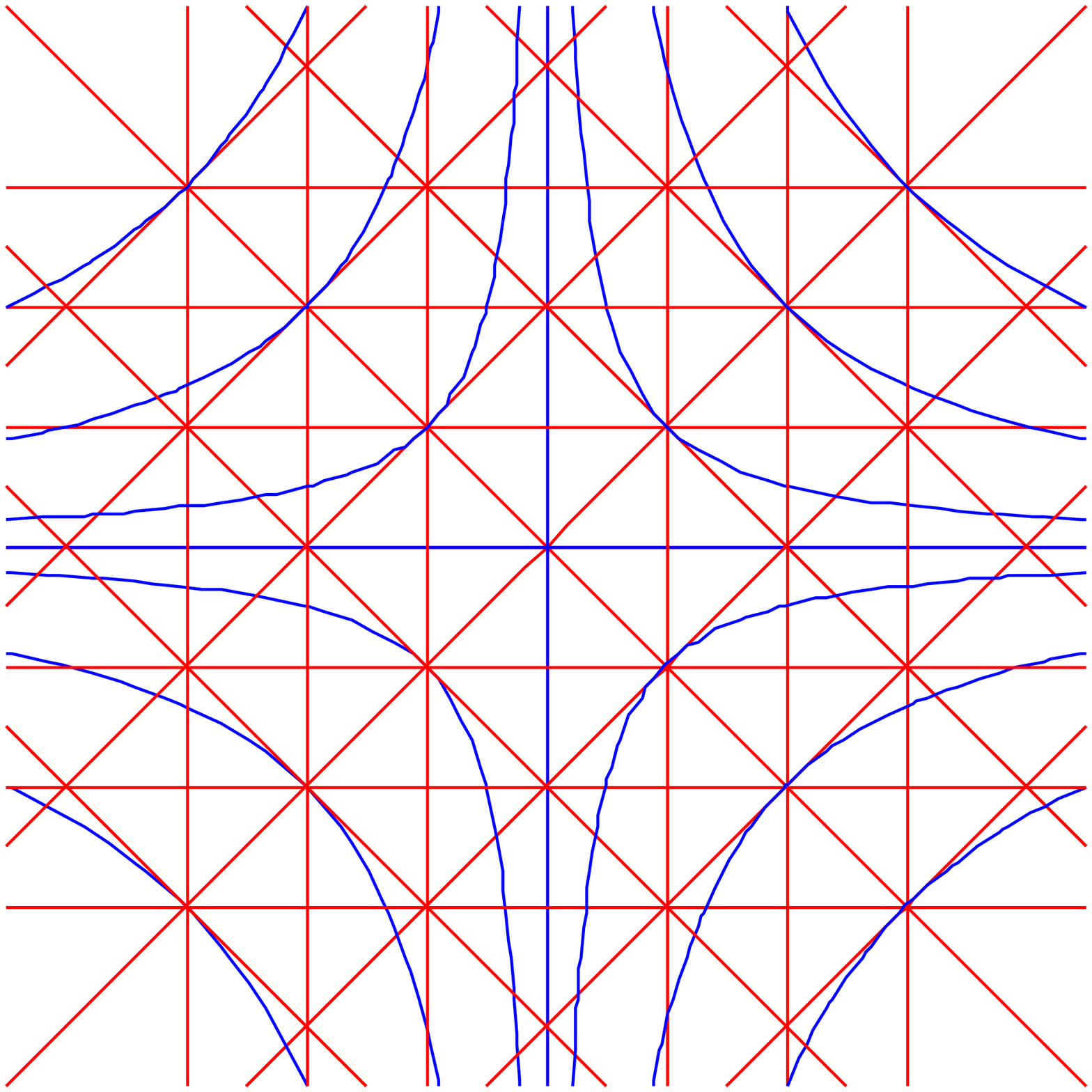}
\end{center}
\end{minipage}
\hfill
\begin{minipage}{4.0cm}
\begin{center}
\includegraphics[width=4.0cm,height=4.0cm]{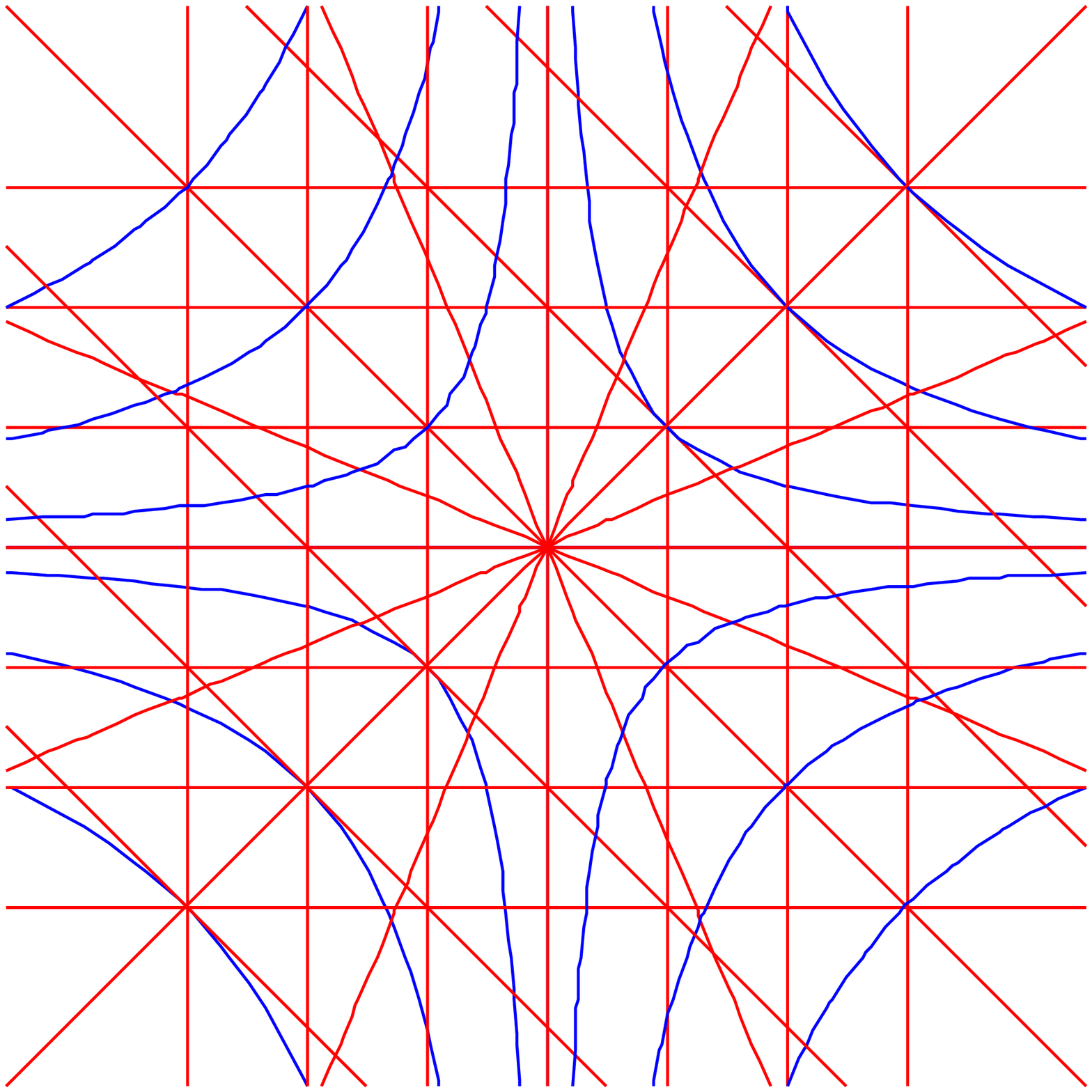}
\end{center}
\end{minipage}
\hfill
\begin{minipage}{4.0cm}
\begin{center}
\includegraphics[width=4.0cm,height=4.0cm]{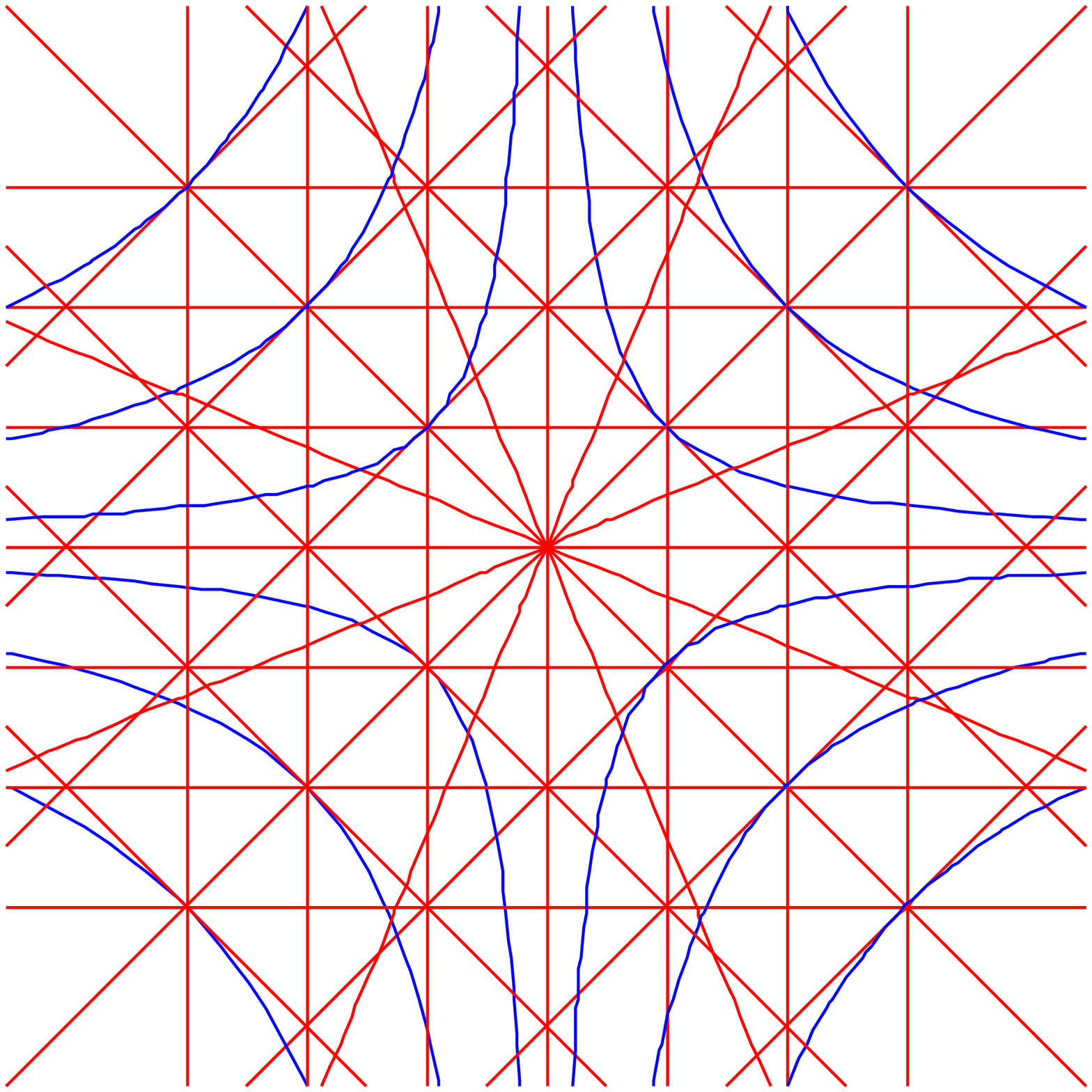}
\end{center}
\end{minipage}
\hfill
\begin{minipage}{4.0cm}
\begin{center}
\includegraphics[width=4.0cm,height=4.0cm]{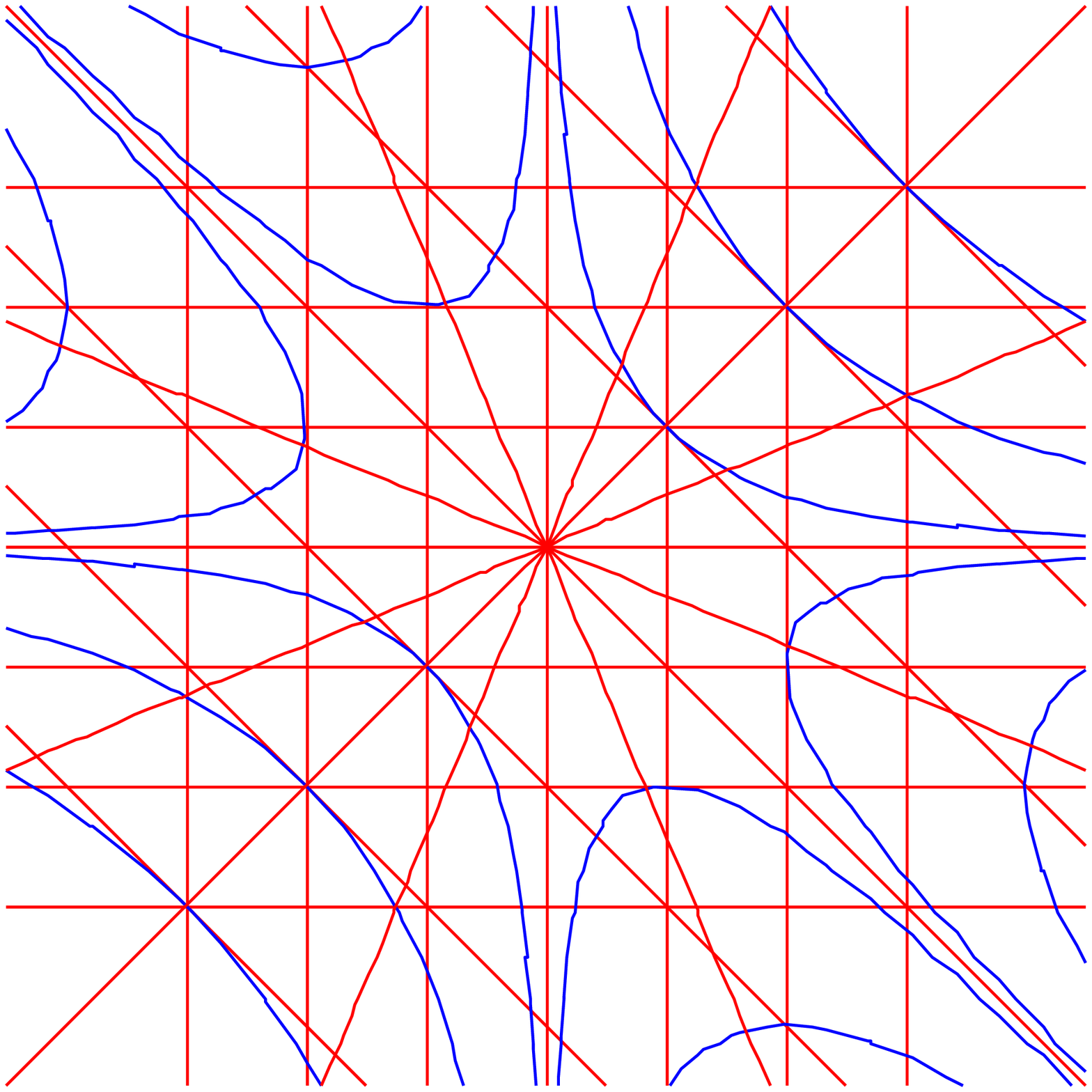}
\end{center}
\end{minipage}
\hfill
\begin{minipage}{4.0cm}
\begin{center}
\includegraphics[width=4.0cm,height=4.0cm]{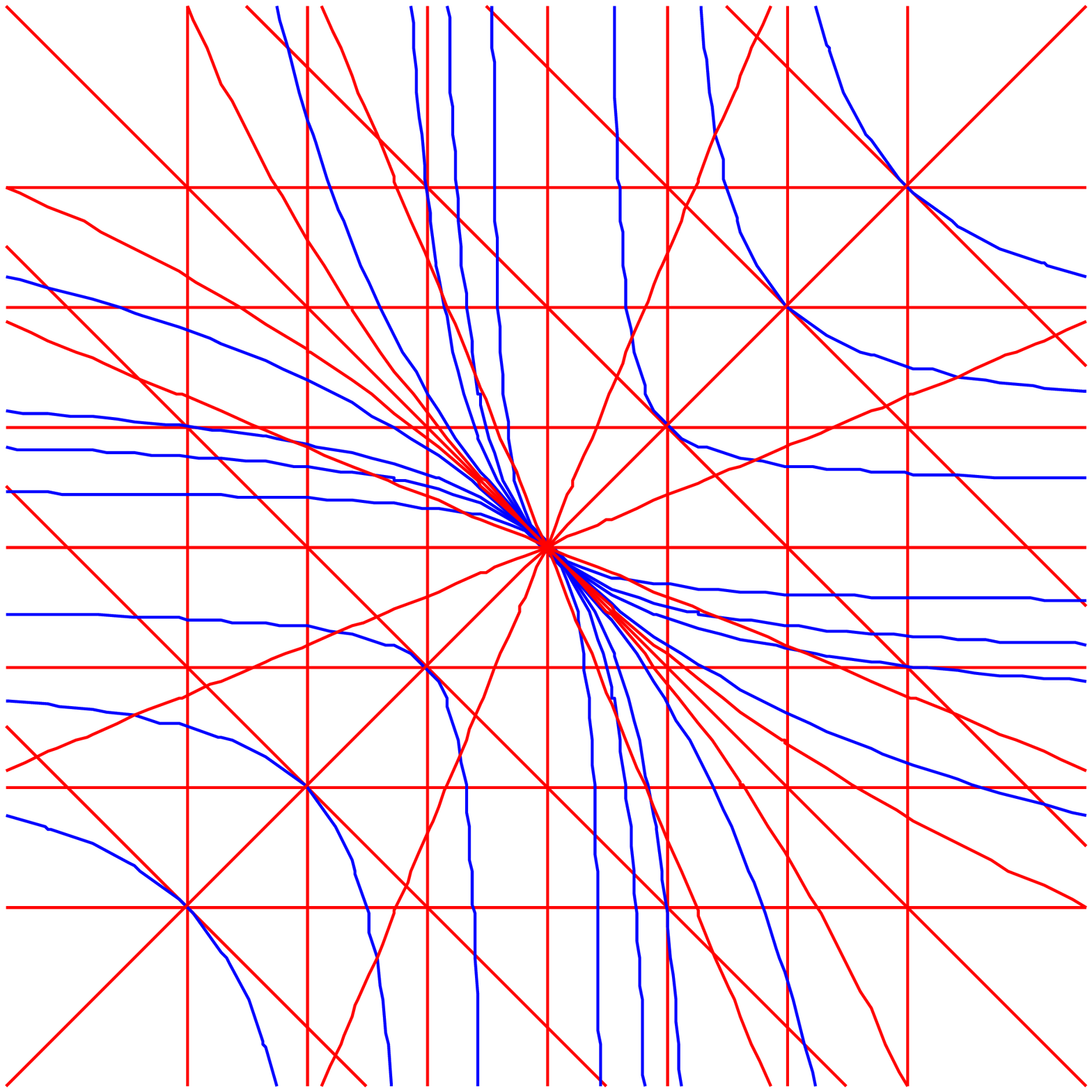}
\end{center}
\end{minipage}
\hfill
\begin{minipage}{4.0cm}
\begin{center}
\includegraphics[width=4.0cm,height=4.0cm]{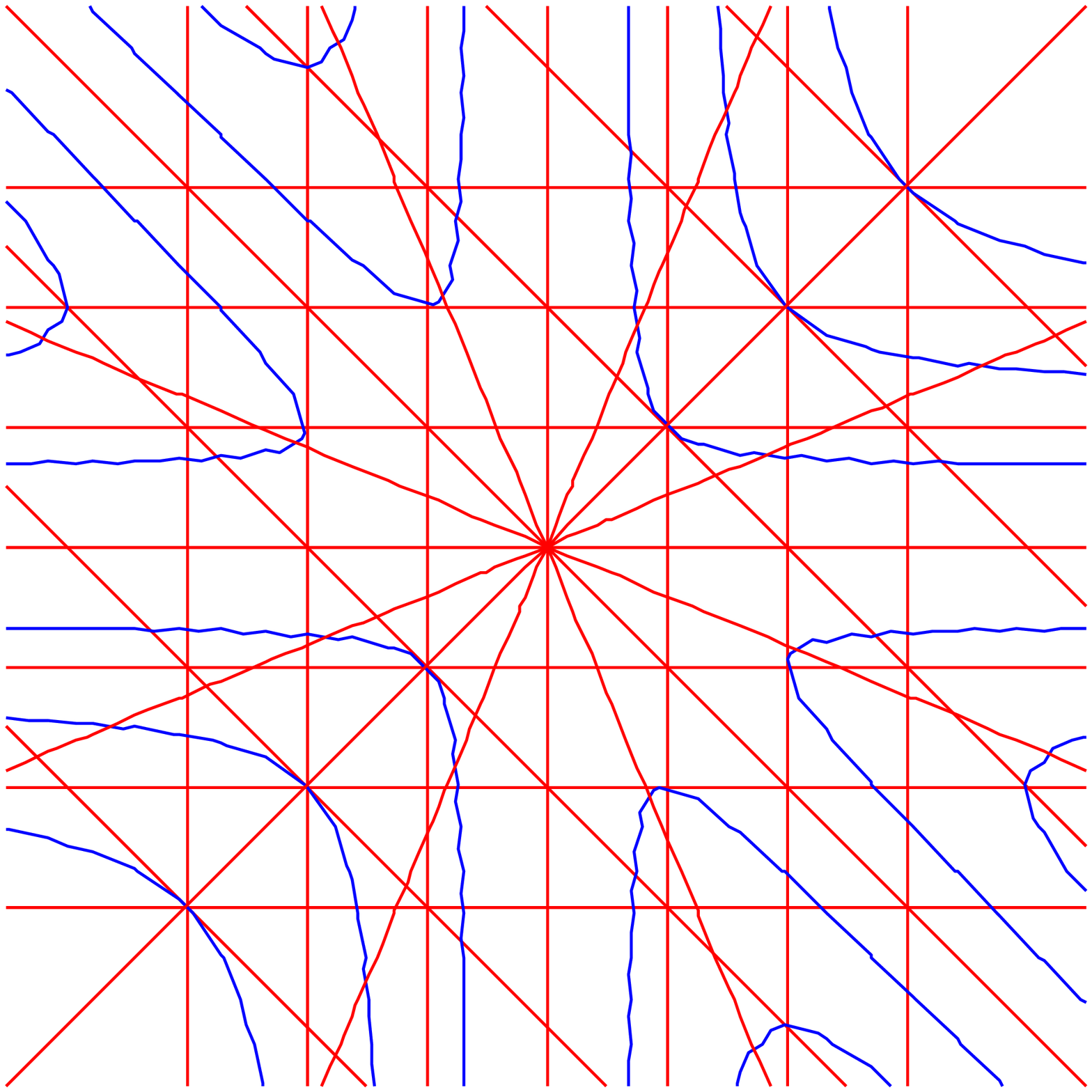}
\end{center}
\end{minipage}
 \caption[] {\label{F:cdqlfig}
{ \small A sample of   real models for  exceptional CDQL webs on $\mathbb P^2$.  In the first
and second rows, the first three members of the infinite family $\mathcal A_I^k$ and  $\mathcal A_{II}^k$ respectively.
In the third row,  from left to right, $\mathcal A^2_{III}$, $\mathcal A^1_{IV}$ and $\mathcal A^2_{IV}$. In the  fourth row:
$\mathcal A_5^a,\mathcal A_5^b$ and $\mathcal A_5^c$.   }}
\end{figure}
On $\mathbb P^2$ the CDQL webs  can be written as
$\mathcal W \boxtimes \mathcal F$ where $\mathcal W$ is a product of
pencil of lines and  $\mathcal F$ is a non-linear foliation. These webs are
determined by the  pair $(\mathcal P, \mathcal F)$ where
 $\mathcal P \subset
\mathbb P^2$ is the set of singularities of the linear  foliations defining $\mathcal W$. One key example is the already mentioned Bol's
$5$-web. It is the exceptional CDQL $5$-web on $\mathbb P^2$ with  $\mathcal F$ equal to  the
 pencil generated by two reduced  conics intersecting transversely  and $\mathcal P$ equal to  the set of four base points of this
pencil. Other examples of exceptional CDQL webs on the plane have appeared in
\cite{PTese,Robert}.

We will deduce from Theorem \ref{TT:curvatura} a   complete classification of exceptional CDQL
webs on the projective plane. In succinct terms it can be stated as follows:

\begin{THM}\label{T:2} Up to projective
automorphisms,  there are exactly four infinite families  and thirteen  sporadic  exceptional CDQL webs on  ${\mathbb P}^2$.
\end{THM}

In  suitable affine coordinates $(x,y)
\in \mathbb C^2 \subset \mathbb P^2$, the four infinite families are
\begin{align*}
\mathcal A_I^k =& \;  \big[ (dx^k - dy^k)\big] \boxtimes \big[ d(xy)  \big] && \mbox{where }\,  k\ge 4 \,; \\
\mathcal A_{II}^k = &\;  \big[(dx^k - dy^k)\, (xdy - ydx) \big] \boxtimes \big[d(xy) \big] &&   \mbox{where }\, k\ge 3 \,;\\
 \mathcal A_{III}^k = &\;  \big[(dx^k - dy^k)\, dx \, dy\big] \boxtimes \big[ d(xy)  \big] && \mbox{where }\, k\ge 2\, ;\\
\mathcal A_{IV}^k  = & \; \big[ (dx^k - dy^k)\,dx \, dy \, (xdy - ydx) \big] \boxtimes \big[ d(xy) \big]    \hspace{-1cm} && \mbox{where }\,  k\ge 1.
\end{align*}

The diagram below shows how these webs relate to each other in terms of inclusions for a fixed $k$. Moreover if $k$ divides $k'$ then
$\mathcal A_I^k, \mathcal A_{II}^k,\mathcal A_{III}^k,\mathcal A_{IV}^k$ are subwebs of  $\mathcal A_I^{k'}, \mathcal A_{II}^{k'},
\mathcal A_{III}^{k'},\mathcal A_{IV}^{k'}$ respectively.
\[
\xymatrix@R=0.1cm@C=0.8cm{  & \mathcal A_{II}^k  \ar[dr]  \\
  \mathcal A_{I}^k  \ar[rr] \ar[dr] \ar[ur] &  &   \mathcal A_{IV}^k  \, .\\
  & \mathcal A_{III}^k  \ar[ur]
  }
\]

All the webs above are invariant by the $\mathbb
C^*$-action $t\cdot(x,y) = (tx,ty)$ on $\mathbb P^2$. Among the
thirteen sporadic examples of exceptional CDQL webs on the projective plane, seven (four $5$-webs, two
$6$-webs and one $7$-web) are also invariant by the same  $\mathbb
C^*$-action. They are:
\[
\begin{array}{lclcl}
\mathcal A_5^a &=& \big[dx\, dy\, (dx+dy)\, (xdy - ydx)  \big] &\boxtimes& \big[d\big(xy(x+y)\big) \big]\, ;
\vspace{0.10cm}
 \\
\mathcal A_5^b &=& \big[
dx\, dy\, (dx+dy)\, (xdy - ydx) \big] &\boxtimes&
\left[d\big(\frac{xy}{x+y}\big)  \right] ;\vspace{0.10cm}
\\
\mathcal A_5^c &=& \big[
dx\, dy\, (dx+dy)\, (xdy - ydx)\big] &\boxtimes&
\left[d\big(\frac{x^2 + xy + y^2}{xy\,(x+y)}\big)  \right] \, ; \vspace{0.10cm}
\\
\mathcal A_5^d &=& \big[  dx\, (dx^3 + dy^3) \big] &\boxtimes&
\left[d\big(x(x^3+y^3) \big)  \right] ;
\vspace{0.10cm} \\
\mathcal A_6^a &=& \big[
 dx\, (dx^3 + dy^3)\, (xdy - ydx) \big]
&\boxtimes& \left[d\big(x(x^3+y^3) \big)  \right] ; \vspace{0.10cm}\\
\mathcal A_6^b &=& \big[ dx\, dy\, (dx^3 + dy^3) \big] &\boxtimes& \left[d(x^3+y^3)\right] ; \vspace{0.10cm}\\
\mathcal A_7 &=& \big[  dx\, dy\, (dx^3 + dy^3)\, (xdy - ydx) \big] &\boxtimes&\left[d(x^3+y^3)  \right] .
\end{array}
\]

Four of the remaining six sporadic exceptional CDQL webs {(one $k$-web for each $k\in \{5,6,7,8\}$)}  share the same non-linear foliation
$\mathcal F$: the pencil of conics through four points in general
position. For them the set $\mathcal P$ is a subset of
$\mathrm{sing}(\mathcal F)$ containing the base points of the
pencil. Up to automorphism of $\mathcal F$ there is just one choice
for each possible cardinality. They all have been previously known (see  \cite{Robert}).
\[
\begin{array}{lclcl}
\mathcal B_5 &=& \left[ dx \, dy  \, d\big(\frac{x}{1-y}
\big) \,  d\big(\frac{y}{1-x} \big) \right] &\boxtimes&
\left[ d
\big(\frac{xy }{(1-x)(1-y)}\big) \right];
\vspace{0.10cm}\\
\mathcal B_6 &=&
\mathcal B_5
 &\boxtimes& \big[ d \left(x + y  \right) \big] ; \vspace{0.10cm} \\
\mathcal B_7 &=&
\mathcal B_6
&\boxtimes&
\left[ d \big(\frac{x}{y} \big) \right];  \vspace{0.10cm}   \\
\mathcal B_8 &=&
\mathcal B_7
&\boxtimes&
\left[ d \big(\frac{1-x}{1-y} \big) \right] .
\end{array}
\]

The last two sporadic CDQL exceptional webs (one $5$-web and one $10$-web) also
share the same non-linear foliation: the Hesse pencil of cubics. Recall that
this pencil is the one generated by a smooth cubic and its
Hessian and that it is unique up to automorphisms of $\mathbb P^2$. These webs are (with   $\xi_3=\exp(2i\pi/3)$):
\[
\begin{array}{l}
\mathcal H_5\; =  \left[ (dx^3 + dy^3 ) \, d\big(\frac{x}{y}
\big)  \right] \boxtimes \left[ d \big(\frac{x^3 + y^3 + 1
}{xy}\big) \right]; \vspace{0.15cm} \\
 \mathcal H_{10} = \left[  (dx^3 + dy^3 ) \Big(\prod_{i=0}^2
d\big(\frac{y-\xi_3^i }{x } \big)\Big) \Big(\prod_{i=0}^2
d\big(\frac{x-\xi_3^i }{y } \big)\Big)\right] \boxtimes \left[
d \big(\frac{x^3 + y^3 + 1 }{xy}\big) \right] .
\end{array}
\]

The  web ${\mathcal H}_{10}$ shares a number of features with Bol's web $\mathcal B_5$. They
both have a huge group of birational automorphisms (the symmetric group ${\mathrm S}_5$ for
$\mathcal B_5$ and Hesse's group $G_{216}$ for ${\mathcal H}_{10}$),  {both are} naturally associated to nets in the sense
of Section \ref{S:nets}  and their abelian relations
can be  expressed in terms of  logarithms and dilogarithms.

\medskip

Because they have parallel $4$-subwebs whose slopes have non real cross-ratio  the webs $\mathcal A_{III}^k$, $\mathcal A_{IV}^k$ for
$k \ge 3$, $\mathcal A_5^d, \mathcal A_6^a, \mathcal A_6^b$ and $\mathcal A_7$
do not admit  real models. The web $\mathcal H_{10}$ also does not admit a real model. {To verify this fact, one} possibility  is to observe  that  the lines passing through two
of the nine  base points always contain a third and notice that this contradicts  Sylvester-Gallai Theorem \cite{Coxeter}:  for every finite set
of non collinear points in $\mathbb P^2_{\mathbb R}$ there exists a line containing exactly two points of the set.
All the other exceptional CDQL webs admit real models. Some of them are pictured in  Figure \ref{F:cdqlfig}.

\subsection{Exceptional CDQL webs on Hopf surfaces}

The classification of CDQL webs on $\mathbb P^2$ admits as a corollary  the classification of exceptional CDQL webs on Hopf surfaces.

\begin{COR}
Up to automorphism, the only exceptional CDQL webs on Hopf surfaces are quotients   of  the restrictions of
the webs $\mathcal A^*_*$ to $\mathbb C^2 \setminus \{ 0 \}$ by the group of deck transformations.
\end{COR}

The proof is automatic. One has just to remark that a foliation on a Hopf surface of type $H_{\alpha}$ when lifted to $\mathbb C^2 \setminus \{ 0 \}$
gives rise to an algebraic  foliation on $\mathbb C^2$ invariant by the $\mathbb C^*$-action $t\cdot(x,y) = (tx,ty)$.

\subsection{From global to local\ldots}
Although based on global methods, the classification of exceptional CDQL webs on $\mathbb P^2$ also yields
information about the {singularities of} local exceptional webs.

\begin{COR}\label{CC:local} Assume that $k\ge 4$.
Let $\mathcal W$ be a smooth $k$-web and  $\mathcal F$ be a
 foliation, both defined on $(\mathbb C^2,0)$. If the $(k+1)$-web $\mathcal W \boxtimes \mathcal F$ has maximal rank then one of the following situations holds:
\begin{enumerate}
\item the foliation $\mathcal F$ is of the form $\big[ H(x,y) ( \alpha d x + \beta dy) + {h.o.t.} \big]$
where $H$ is a non-zero homogeneous polynomial and $(\alpha, \beta) \in \mathbb C^2 \setminus \{ 0 \}$;\vspace{0.1cm}
\item the foliation $\mathcal F$ is of the form $\big[ H(x,y) ( y d x - x  dy) + {h.o.t.} \big]$
where $H$ is a non-zero homogeneous polynomial;\vspace{0.1cm}
\item $\mathcal W \boxtimes \mathcal F$ is exceptional and its first non-zero jet defines, up to linear automorphisms,   one of the following  webs
\[
\mathcal A_I^k, \,  \mathcal A_{III}^{k-2}, \,  \mathcal A_5 ^d\;  (\text{only when } k=4 )\;\; \mbox{ or } \;  \; \mathcal A_6^b \; (\text{only when } k=5) \, .
\]
\end{enumerate}
\end{COR}

In fact, as it will be clear from its proof, it is possible to state a slightly more general result
in the same vein. Nevertheless the result above suffices for   the classification
of exceptional CDQL webs  on complex tori.

\subsection{\ldots and back: classification of exceptional CDQL  webs on tori}
A  CDQL web on a torus is   the superposition   of a
non-linear foliation with a product of foliations induced by global holomorphic $1$-forms.
Since \'{e}tale coverings between complex tori abound and  {because} the pull-back
of exceptional CDQL webs under these are still exceptional CDQL  webs, we are naturally
lead to extend the notion of isogenies between complex tori.
Two webs $\mathcal W_1, \mathcal W_2$ on  complex tori $T_1,T_2$ are {\it isogeneous}  if there exist  a  complex torus $T$
and \'etale morphisms $\pi_i: T \to T_i$ for $i=1,2$,
such that $\pi_1^* (\mathcal W_1) = \pi_2^* (\mathcal W_2)$.

\begin{THM}
\label{T:CDQLontori}
Up to isogenies, there are exactly
three sporadic  (one for each $k \in \{ 5,6,7\}$) and one continuous family (with $k=5$) of exceptional CDQL $k$-webs  on complex tori.
\end{THM}

 The elements of the continuous family are
\[
\mathcal E_{\tau} = \big[ \, dx\,  dy \, (dx^2 - dy^2) \, \big] \boxtimes \left[ d\left(\frac{\vartheta_1(x,\tau)  \vartheta_1(y,\tau) }{\vartheta_4(x,\tau)  \vartheta_4(y,\tau) }\right)^{\!\!2}\, \right]
\]
defined, respectively, on the torus $E_{\tau}^2$ for  arbitrary $\tau \in \mathbb H{=\{ z \in \mathbb C \, | \, \Im{\rm m}(z)>0\, \}}$ where $E_{\tau} = \mathbb C /  (\mathbb Z \oplus \mathbb Z \tau)$.
The functions $\vartheta_i$ involved in the definition are the classical Jacobi theta functions, see Example \ref{Ex:theta}.

These  webs first appeared in  Buzano's work \cite{Buzano}   but their rank was not determined at that time. They were later rediscovered in \cite{PT} where it is proved that they are all exceptional and  that
$\mathcal E_{\tau}$ is isogeneous to $\mathcal E_{\tau'}$ if and only if $\tau$ and $\tau'$
belong to the same orbit under the natural action on $\mathbb H$  of the $\mathbb Z/2\mathbb Z$ extension of $\Gamma_0(2) \subset \mathrm{PSL}(2, \mathbb Z)$
generated by $\tau \mapsto -2 \tau^{-1}$. Thus the continuous  family of exceptional CDQL webs on tori is parameterized by a $\mathbb Z/2\mathbb Z$-quotient of the modular curve $X_0(2)$.

\smallskip

The sporadic  CDQL $7$-web $\mathcal E_7$ is strictly related to a particular element of the previous family. Indeed $\mathcal E_7$ is the $7$-web on  $E_{1+i}^2$
\[
\mathcal E_7 = \big [dx^2 + dy^2\big]  \boxtimes {\mathcal E}_{1+i}\,.
\]

\smallskip

The sporadic  CDQL $5$-web $\mathcal E_5$ lives naturally in $E_{\xi_3}^2$ and can be described as
\[
 \big[ dx\, dy\, (dx - dy) \, (dx + \xi_3^2\,  dy)  \big] \boxtimes \left[ d\Big( \frac{\vartheta_1(x, \xi_3)\vartheta_1(y, \xi_3 )\vartheta_1(x-  y , \xi_3)\vartheta_1(x +  \xi_3^2 \, y, \xi_3) }{\vartheta_2(x, \xi_3)\vartheta_3(y, \xi_3 )\vartheta_4(x- y ,  \xi_3)\vartheta_3(x +  \xi_3^2 \,y, \xi_3)}\Big) \right] .
\]

\smallskip

The sporadic CDQL $6$-web $\mathcal E_6$ also lives in $E_{\xi_3}^2$ and  is best described in terms of Weierstrass $\wp$-function.
\[
\mathcal E_6 = \big[\,  dx\,  dy \, (dx^3 + dy^3) \big] \boxtimes
\left[
{\wp(x,\xi_3)}^{-1}
 {dx} + {\wp(y,\xi_3)}^{-1}{dy}\right] .
\]
Although not completely evident from the above presentation,  it turns out that the foliation $[
{\wp(x,\xi_3)}^{-1}  {dx} + {\wp(y,\xi_3)}^{-1}{dy}]$ admits a rational first integral, see Proposition \ref{P:33}.

A more geometric description of these exceptional {\it elliptic webs} will be given in Section \ref{S:elliptic} together with the proof that they are
indeed exceptional.

\smallskip

The proof of Theorem \ref{T:CDQLontori} follows the same lines of the proof of Theorem \ref{T:2} but with some
 twists.  The key  extra ingredients are  Corollary \ref{CC:local} and  the following (considerably easier)  analogue for two dimensional
complex tori of \cite[Theorem 1]{PY}.

\begin{THM}\label{T:PYontori}
If $T$ is a two-dimension complex tori and $f:T \dashrightarrow \mathbb P^1$ a meromorphic  map then the number
of linear fibers of $f$, when finite, is at most six.
\end{THM}

For us the linear fibers of a rational map from a two-dimensional complex torus to a curve are the ones that are
set-theoretically equal to a union of subtori.

\subsection{Plan of the Paper} The remaining of the paper can be roughly divided in four parts.
The first goes from Section \ref{S:action} to Section \ref{S:elliptic} and is devoted to
prove that all the webs presented in the Introduction are exceptional. The highlights are Theorems \ref{T:net}
and Theorem \ref{T:ellipticnet} that show that the webs $\mathcal B_5, \mathcal H_{10}, \mathcal E_{\tau}, \mathcal E_5, \mathcal E_6$ and $\mathcal E_7$ are exceptional thanks to essentially  the  same reason. Their abelian relations are expressed in terms of logarithms, dilogarithms and their elliptic
counterparts.
Sections \ref{S:bary}, \ref{S:fbary} and \ref{S:curvatura} form the second part of the paper which is  devoted to the study of the
$\mathcal F$-barycenter of a web. Besides the proof of Theorem \ref{TT:curvatura} of the Introduction it also contains a
very precise description of  the barycenters of decomposable linear  webs centered at linear foliations  on $\mathbb P^2$.
This description  lies at the heart of our approach
 to the classification of exceptional CDQL webs on $\mathbb P^2$.
The third  part of the paper goes from Section \ref{S:constraints} to Section \ref{S:flat}
and contains the classification of exceptional CDQL webs on the projective plane.
Finally the fourth and last part is contained in the last two sections and
deals with the classification of exceptional CDQL webs on complex tori.
Beside  this classification it also contains  the proofs of Corollary \ref{CC:local}
and Theorem \ref{T:PYontori}.

\subsection{Acknowledgements}
The first author thanks Jorge Pastore  for enlightening discussions.
The second author thanks Frank Loray and the International Cooperation Agreement Brazil-France.
Both authors {are grateful to} Marco Brunella for the elegant proof of Proposition \ref{P:Brunella} and
to David Mar\'{\i}n for the explicit expression for $\beta_*$ presented in Remark \ref{R:dynamics}.

\pagebreak


\section{Abelian relations for CDQL webs invariant by $\mathbb
C^*$-actions}\label{S:action}

We start things off with the following well-known proposition.

\begin{prop}\label{P:exc}
Let $\mathcal W$ be a linear $k$-web of maximal rank and $\mathcal F$ be a non-linear foliation on $(\mathbb C^2,0)$.
The $(k+1)$-web $\mathcal W \boxtimes \mathcal F$ is exceptional if and only if it has maximal rank and $k\ge 4$.
\end{prop}
\begin{proof}
For $k\le 3$, all $(k+1)$-webs of maximal rank are algebraizable
thanks to Lie's Theorem. Suppose  that $k\ge 4$  and let
$\varphi:(\mathbb C^2,0) \to (\mathbb C^2,0)$ be a biholomorphism algebraizing $\mathcal W \boxtimes \mathcal F$.
 Since $\mathcal W$ has maximal rank,  $\varphi^* \mathcal W$ must be algebraic. According to
\cite{Henautlinear} (see also \cite[p. 247]{BB}) the
 biholomorphism    $\varphi$ must be the restriction
of a projective transformation.   It follows that $\varphi^* \mathcal F$  is
non-linear and consequently  it cannot exist an algebraization of $\mathcal W
\boxtimes \mathcal F$.
\end{proof}

As a corollary one sees that in order to prove that  a CDQL $k$-web ${\mathcal W}$ is exceptional, when  $k\ge 5$,
it suffices to verify that it has maximal rank. The most obvious way
to accomplish this task is to  exhibit a basis of the  space of its abelian
relations. In general, the explicit determination of
${\mathcal A}({\mathcal W})$ is a fairly difficult problem. To  our knowledge, the only general
method available is Abel's method for solving functional equations
(see \cite{Abel} and \cite[Chapitre 2]{PTese}). It assumes  the knowledge of first integrals for the
defining foliations of ${\mathcal W}$ and   it tends to involve rather lengthy computations.

In particular cases there are more efficient ways to determine the space of
abelian relations. For instance, if the web admits an
infinitesimal automorphism  then the results of \cite{MPP} reduce
the problem to plain linear algebra. In Section \ref{S:reviewMPP} we
recall the results of \cite{MPP} and use them in Sections \ref{S:4} and \ref{S:7} to deal with  the CDQL webs invariant by
$\mathbb C^*$-actions  described in the Introduction. We point out
that the content of Section \ref{S:reviewMPP}  plays a decisive role
in the classification of  exceptional CDQL webs of degree one
carried out in Section \ref{S:deg1}.

\subsection{Webs with infinitesimal automorphisms}\label{S:reviewMPP}

Let $\mathcal F$ be a regular foliation on $(\mathbb C^2,0)$ induced
by a $1$-form $\omega$. We say that a  vector field $X$ is an
{\it infinitesimal automorphism} of $\mathcal F$ if  $ L_X \omega \wedge
\omega = 0 \, .$ When such infinitesimal automorphism $X$ is
transverse to $\mathcal F$, that is when  $\omega(X)\neq 0$, then  the
$1$-form \[ \eta = \frac{\omega}{i_X \omega}
\]
is   closed and satisfies $L_X \eta =0$. By definition, the integral
\[
u(z) = \int^z_0 \eta
\]
is the {\it canonical first integral} of ${\mathcal F}$ (with
respect to $X$).

Assume now that $\mathcal W$ is a  regular $k$-web on $(\mathbb
C^2,0)$ induced by  $1$-forms $\omega_1, \ldots, \omega_k$ and
let  $X$ be an  infinitesimal automorphism of   all the defining
foliations of $\mathcal W$.

The Lie derivative $L_X$ induces a linear map
\begin{eqnarray}
\label{themap!}
L_X : \mathcal A(\mathcal W) & \to & \mathcal A(\mathcal W) \\
(\eta_1,\ldots ,\eta_k) & \mapsto& (L_X\eta_1 ,\ldots, L_X \eta_k)
\, .  \nonumber
\end{eqnarray}

The study of this linear map leads to the following proposition.

\begin{prop} \label{P:description}
Let $\lambda_1,\ldots,\lambda_\tau \in {\mathbb C}$ be the
eigenvalues of the map $L_X$ acting on ${\mathcal A}({\mathcal W})$
corresponding to  minimal eigenspaces with dimensions
 $\sigma_1,\ldots,\sigma_\tau $.
The abelian relations of $ {\mathcal W}$ are of the form
 $$ P_1(u_1)\,e^{\lambda_i\,u_1}\,du_1+\cdots+ P_k(u_k)\,e^{\lambda_i\,u_k}\,du_{k}=0
$$
where $P_1, \ldots, P_k$ are  polynomials of degree less   or equal
to $ \sigma_i$.
Moreover the abelian relations corresponding to
eigenvectors are precisely the ones for which the polynomials $P_i$
are constant.
\end{prop}

Proposition \ref{P:description} suggests an effective method to
determine $\mathcal A(\mathcal W)$ from the study of the linear map (\ref{themap!}). For details see \cite{MPP}.
It also follows from the study of  (\ref{themap!}) the
main result of \cite{MPP}.

\begin{thm}\label{T:1}
Let ${\mathcal W}$ be a $k$-web which admits a transverse
infinitesimal automorphism $X$. Then
\[
{rk}(\mathcal W \boxtimes \mathcal F_X) ={rk}(\mathcal
W) + (k -1)\, .
\]
In particular,  ${\mathcal W}$ is of maximal rank if and only if
  $\mathcal W \boxtimes \mathcal F_{X}$ is also of maximal rank.
\end{thm}

Below we will make use of Theorem \ref{T:1} to prove that certain webs have maximal rank without giving a complete list of their abelian relations. Nevertheless, the proof of Theorem \ref{T:1} (see \cite{MPP})
is constructive and the interested reader can easily determine  a complete list of the  abelian relations.

\subsection{Four infinite families}\label{S:4}
 \label{S:homogeneous} Recall  the definition of the webs $\mathcal A_{I}^k, \mathcal A_{II}^k, \mathcal A_{III}^k, \mathcal A_{IV}^k$:
\begin{align*}
\mathcal A_I^k =& \;  \big[ (dx^k - dy^k)\big] \boxtimes \big[ d(xy)  \big] && \mbox{where }\,  k\ge 4 \,; \\
\mathcal A_{II}^k = &\;  \big[(dx^k - dy^k)\, (xdy - ydx) \big] \boxtimes \big[d(xy) \big] &&   \mbox{where }\, k\ge 3 \,;\\
 \mathcal A_{III}^k = &\;  \big[(dx^k - dy^k)\, dx \, dy\big] \boxtimes \big[ d(xy)  \big] && \mbox{where }\, k\ge 2\, ;\\
\mathcal A_{IV}^k  = & \; \big[ (dx^k - dy^k)\,dx \, dy \, (xdy - ydx) \big] \boxtimes \big[ d(xy) \big]    \hspace{-1cm} && \mbox{where }\,  k\ge 1.
\end{align*}

The exceptionality of these webs follows from the next proposition.

\begin{prop}
For all $k\ge 1$ the webs $\mathcal A_{I}^k, \mathcal A_{II}^k, \mathcal A_{III}^k$ and $ \mathcal A_{IV}^k$ have maximal rank.
\end{prop}
\begin{proof}
Let
 $R= x \frac{\partial}{\partial x} + y \frac{\partial}{\partial y}$ be the radial vector field. Note that it is an infinitesimal automorphism of all the  webs above. Moreover
\[
\mathcal A_{II}^k = \mathcal A_{I}^k \boxtimes \mathcal F_R  \quad \text{and } \quad \mathcal A_{IV}^k = \mathcal A_{III}^k \boxtimes \mathcal F_R  \, .
\]
It follows from Theorem \ref{T:1} that $\mathcal A_{II}^k$ (resp. $\mathcal A_{IV}^k$) has maximal rank if and only if  $\mathcal A_{I}^k$ (resp. $\mathcal A_{III}^k$)
also does.

To prove that $\mathcal A_{I}^k$ has maximal rank  consider the linear
automorphism of $\mathbb C^2$,    $\varphi(x,y)=  (x, \xi_k y)$. Consider also  the
induced automorphism  of   the vector space $\mathbb C_{2k-2}[x,y]$  of homogeneous polynomials
of degree $2k-2$:
\begin{eqnarray*}
\varphi^*: \mathbb C_{2k-2}[x,y] &\longrightarrow& \mathbb C_{2k-2}[x,y] \\
p &\mapsto& p \circ \varphi \, .
\end{eqnarray*}
For $k=1$ there is nothing to prove: every $2$-web has maximal rank. Assume that $k \ge 2$.
If $\xi_k = \exp(2\pi i /k)$ then   the $(\xi_k^{k-1})$-eigenspace of
$\varphi^*$ has dimension one and is generated by $(xy)^{k-1}$.

If $V\subset \mathbb C_{2k-2}[x,y]$ denotes  the vector subspace generated
by the homogeneous polynomials $(x-\xi_k^i y)^{2k-2}$ with $i$ ranging from $0$ to $k-1$, then
${\varphi^*}$ preserves $V$ and the characteristic polynomial of
${\varphi^*}_{|V}$ is equal to $t^k -1$. It follows that there exists ${p \in
V\setminus \{0\}}  $ such that $\varphi^* p = (\xi_k^{k-1})\, p$. Since the eigenspace of  $\varphi^*$ associated to the eigenvalue
$\xi_k^{k-1}$
 has dimension one,
$p$ must be a complex multiple of $(xy)^{k-1}$. Therefore, there exist complex constants $\mu_1,\ldots,\mu_k$ such  that
\[
  (xy)^{k-1}  = \sum_{i=1}^k \mu_i \,(x-\xi_k^i  y)^{2k-2} \; .
\]
This  identity can interpreted as an abelian relation of $ \mathcal A_{I}^k$.
If we apply the second-order differential  operator
$\frac{\partial^2}{\partial x\partial y}$
to it   we obtain another abelian relation
\[
  (k-1)^2  (xy)^{k-2}  = \sum_{i=1}^k \mu_i (2k-2)(2k-1)\xi_k^i(x-\xi_k^i  y)^{2(k-1)-2} \, .
\]
When $k\ge 3$,  this abelian relation is clearly linearly independent from the previous one.
Iteration of this procedure shows  that
\[
\dim \frac{\mathcal A(\mathcal A_I^k )}{\mathcal A\big([ dx^k - dy^k]\big)} \ge k-1 \, .
\]
Since $[dx^k - dy^k]$ is an algebraic $k$-web its  rank is $(k-1)(k-2)/2$. Thus $\dim \mathcal A(\mathcal A_I^k ) =
k\,(k-1)/2$ and $\mathcal A_I^k $ is
indeed of maximal rank. Theorem \ref{T:1} implies that the  $(k+2)$-web   $\mathcal A_{II}^k $ is also of maximal rank.

The proof that $\mathcal A_{III}^k$ and $\mathcal A_{IV}^k$ are of maximal rank is  analogous. As before, it suffices to show that the
$(k+3)$-web $\mathcal A_{III}^k$ has maximal rank.

Consider now the   induced
automorphism $\varphi^*$ on the space  $\mathbb C_{2k}[x,y]$ of homogeneous polynomials of degree $2k$. The $1$-eigenspace of $\varphi^*$ has dimension three and is generated by $x^{2k}, y^{2k}$ and $(xy)^k$.
If $V\subset \mathbb C_{2k}[x,y]$ denotes now  the vector subspace generated
by the polynomials $(x-\xi_k^i y)^{2k}$ with $i = 0, \ldots, k-1$, then  the characteristic polynomial of
${\varphi}^*_{|V}$ is also equal to $t^k -1$. Thus
there exists an abelian relation of  $\mathcal A_{III}^k$ of the form
\[
  (xy)^{k}  =   \sum_{i=1}^k \mu_i (x-\xi_k^i  y)^{2k} +
 \mu_{k+1} x^{2k} + \mu_{k+2} y^{2k}
\, .
\]
Applying the operator $\frac{\partial^2}{\partial x\partial y}$
and iterating as above one deduces that
\[
\dim \frac{\mathcal A(\mathcal A_{III}^k )}{\mathcal A\big([ dx\, dy\, (dx^k - dy^k)]\big)} \ge k \, .
\]
Taking into account the logarithmic  abelian relation
\[
 \log (xy) = \log x + \log y
\]
we conclude that $\dim \frac{\mathcal A(\mathcal A_{III}^k )}{\mathcal A([ dx\,dy\,(dx^k - dy^k)])} \ge k+1.$ Since $[dxdy(dx^k - dy^k)]$ has rank $k\,(k+1)/2$,
it follows that the $(k+3)$-web $\mathcal A_{III}^k$ also has maximal rank.
\end{proof}

\subsection{The seven sporadic exceptional CDQL webs invariant by $\mathbb C^*$-actions}
 \label{S:7} Recall from the Introduction the other seven webs invariant by the $\mathbb C^*$-action \mbox{$t\cdot(x,y) \mapsto (tx,ty)$}:
\[
\begin{array}{lclcl}
\mathcal A_5^a &=& \big[dx\, dy\, (dx+dy)\, (xdy - ydx)  \big] &\boxtimes& \big[d\big(xy(x+y)\big) \big]\, ;
\vspace{0.1cm}
 \\
\mathcal A_5^b &=& \big[
dx dy (dx+dy) (xdy - ydx) \big] &\boxtimes&
\left[d\big(\frac{xy}{x+y}\big)  \right] ;\vspace{0.1cm}
\\
\mathcal A_5^c &=& \big[
dx dy (dx+dy) (xdy - ydx)\big] &\boxtimes&
\left[d\big(\frac{{x^2 + xy + y^2}}{xy(x+y)}\big)  \right]  ; \vspace{0.1cm}
\\
\mathcal A_5^d &=& \big[  dx dy (dx+dy)(dx - \xi_3 dy) \big] &\boxtimes&
\left[d\big(xy(x+y)(x-\xi_3y) \big)  \right] ;
\vspace{0.1cm} \\
\mathcal A_6^a &=& \big[
 dx dy (dx+dy)(dx - \xi_3 dy) (xdy - ydx) \big]
&\boxtimes& \left[d\big(xy(x+y)(x-\xi_3y) \big)  \right] ; \vspace{0.1cm}\\
\mathcal A_6^b &=& \big[ dx dy (dx^3 + dy^3) \big] &\boxtimes& \left[d(x^3+y^3)\right] ; \vspace{0.15cm}\\
\mathcal A_7 &=& \big[  dx dy (dx^3 + dy^3) (xdy - ydx) \big] &\boxtimes&\left[d(x^3+y^3)  \right] .
\end{array}
\]

Of course, they all   share the same infinitesimal automorphism: the radial
vector field $R$. Because
\[
\mathcal A_{6}^a = \mathcal A_{5}^d \boxtimes \mathcal F_R  \qquad \text{and } \qquad \mathcal A_{7} = \mathcal A_{6}^b \boxtimes \mathcal F_R  \, ,
\]
Theorem \ref{T:1} implies that the maximality of the rank of $\mathcal A_{6}^a$ (resp. $ \mathcal A_{7} $) is equivalent
to the maximality of the rank of   $\mathcal A_{5}^d$ (resp. $ \mathcal A_{6}^b $). Thus, to prove that all the seven webs above are exceptional, it suffices to
prove that $\mathcal A_5^a, \mathcal A_5^b,\mathcal A_5^c,\mathcal A_5^d,\mathcal A_6^b$  have maximal rank.
For this sake we list  below a basis for a subspace of the space of abelian relations of these webs that
is transversal to the space of abelian relations of the maximal linear subweb contained in  {each of} them.

\subsubsection{Abelian Relations for $\mathcal A_5^a$.} If $g_0= xy(x+y)$, $g_1=x$, $g_2=y$, $g_3=x+y$ and $g_4=\frac{x}{y}$
then the sought  abelian relations for $\mathcal A_5^a$ are
\[
\begin{array}{lclclclcl}
\ln g_0 &=& \ln g_1 &+& \ln g_2 &+& \ln g_3\,  \vspace{0.1cm} \\
\ln^2 g_0 &=& 3\ln^2 g_1 &+& 3 \ln^2 g_2 &+& 3\ln^2 g_3 &-& \varphi(g_4)\, \vspace{0.1cm}  \\
3 \,g_0 &=& -g_1^3  &-& g_2^3 &+& g_3^3
\end{array}
\]
where $\varphi(t) = \ln^2 t + \ln^2( t + 1) + \ln^2(t^{-1} + 1)$.

\subsubsection{Abelian Relations for $\mathcal A_5^b$.}
 If $g_0={xy}/{(x+y)}$ and $g_1, g_2, g_3, g_4$ are  as above then
\[
\begin{array}{lclclclcl}
\ln g_0 &=& \ln g_1 &+& \ln g_2 &-& \ln g_3 \, \vspace{0.1cm}\\
\ln^2 g_0 &=& \ln^2 g_1 &+& \ln^2 g_2 &-& 3\ln^2 g_3 &-& \varphi(g_4)\, \vspace{0.1cm}\\
 g_0^{-1} &=& g_1^{-1}  &+&  g_2^{-1}
\end{array}
\]
where $\varphi(t) = \ln^2 t - \ln^2(t+1) - \ln^2(t^{-1} + 1)$.

\subsubsection{Abelian Relations for $\mathcal A_5^c$.}
 If {$g_0={(x^2 + xy + y^2 )}/{\big(xy(x+y)\big)}$} and $g_1, g_2, g_3, g_4$ are  as above then
\[
\begin{array}{lclclclcl}
\ln g_0 &=&  &&  &+& \ln g_3 &+& \ln( g_4 + g_4^{-1} + 1)\;  \vspace{0.1cm}\\
 g_0 &=&  g_1^{-1} &+& g_2^{-1} &-& g_3^{-1}\;  \vspace{0.1cm}\\
 g_0^2 &=& g_1^{-2} &+&  g_2^{-2} &-&  g_3^{-2} .
\end{array}
\]

\subsubsection{Abelian Relations for $\mathcal A_5^d$.}
Notice that $\mathcal A_5^d$ is equivalent to $$\big[  dx\, dy \, ( dx + dy) \, (dx - \xi_3 dy)  \big] \boxtimes
\left[d\big( xy(x+y)(x- \xi_3 y) \big)  \right]$$ under a linear change of coordinates.
 If $g_0=xy(x+y)(x- \xi_3 y)$, $g_4=x - \xi_3\,y$  and $g_1, g_2, g_3$ are  as above then
\[
\begin{array}{ccccccccc}
  \ln g_0  & = & \ln g_1  & + & \ln g_2 & + &  \ln g_3 & +  & \ln g_4 \;   \vspace{0.1cm} \\
  12\, g_0 & = & \left( - 2 - \xi_3 \right)\, g_1^4 & + &  (1 + 2\, \xi_3)\,  g_2^4   & + & \left( 1 - \xi_3  \right)\, g_3 ^4 & + &  (1+2\,\xi_3)\, g_4 ^4 \;  \vspace{0.1cm} \\
  28\, g_0^2 & = & \left( 1 + \xi_3 \right)\,{ g_1}^8 & - &    {g_2}^8   & - & \xi_3\, {g_3} ^8 & - &    {g_4} ^8 \;   .
\end{array}
\]

\subsubsection{Abelian Relations for $\mathcal A_6^b$.}
If $g_0=x^3 + y^3 $, $g_4= x + \xi_3\,y$, $g_5 = x+ \xi_3^2\, y$   and $g_1, g_2, g_3$ are  as above then
\[
\begin{array}{ccccccccccc}
 g_0  & =&{ g_1}^3 &+ & {g_2}^3  &  &  &  &   && \vspace{0.1cm}\\
  \ln g_0  & = & &  &  &  &  \ln g_3 & +  & \ln g_4 &+& \ln g_5 \vspace{0.1cm}\\
   30\, g_0^2 &=& 27\,{g_1}^6 &+& 27\,{g_2}^6 &+& {g_3}^6 &+& {g_4}^6 &+& {g_5}^6 \vspace{0.1cm}\\
   84\, {g_0}^3 &=& 81\, {g_1}^9 &+& 81\,{g_2}^9 &+& {g_3}^9 &+& {g_4}^9 &+& {g_5}^9 \,.
\end{array}
\]

\section{Abelian relations for planar webs associated to nets}\label{S:net}
 The determination of $\mathcal A(\mathcal B_5)$ is due to  Bol, see \cite{Bol}. The determination of $\mathcal A(\mathcal B_6), \mathcal
A(\mathcal B_7)$ and $\mathcal A(\mathcal B_8)$ is treated in
\cite{Robert} (see also \cite{PTese,Piriopoly} for the determination of
$\mathcal A(\mathcal B_6)$ and $\mathcal A(\mathcal B_7)$  through
Abel's method).  In this section  we will prove that  $\mathcal H_5$ and $\mathcal H_{10}$ --- the two remaining exceptional CDQL webs on ${\mathbb P}^2$
presented in the Introduction --- have maximal rank.
We adopt  {here} an approach similar to the one used by Robert in
\cite{Robert} and that can be traced back to \cite{HainMcPherson}. We look for the abelian relations among
$k$-uples of Chen's iterated integrals of logarithmic $1$-forms with
poles on certain hyperplane arrangements.  It turns out that this  particular
class of webs carry logarithmic and dilogarithmic abelian relations thanks to purely combinatorial reasons.

\subsection{Webs associated to nets}
\label{S:nets}
Let $r\ge 3$ be  an integer.
Recall from \cite{Ynet} that a {\it $r$-net} in
$\mathbb P^2$ is a pair $(\mathcal L,\mathcal P)$ where $\mathcal L$
is a finite set of lines partitioned into $r$ disjoint subsets $\mathcal L=
{\sqcup}_{i=1}^r \mathcal L_i$ and $\mathcal P$ is a finite set of
points subjected to the two conditions:
\begin{enumerate}
\item for every $i\neq j$ and every $\ell \in \mathcal L_i$, $\ell'
\in \mathcal L_j$, we have that $\ell \cap \ell' \in \mathcal P$;
\item for every $p \in \mathcal P$ and every $i=1,2, \ldots, r$,
there exists a unique $\ell \in \mathcal L_i$ passing through $p$.
\end{enumerate}
The definition implies that $\mathcal P$ has cardinality $m^2$ and that
 the cardinalities  of the sets
$\mathcal L_i$ do not depend on $i$ and  are all  equal to $m = {{\mathrm Card}(\ell \cap \mathcal P)}$ for any $\ell \in \mathcal L$. We say that $\mathcal L$
is a {\it $(r,m)$-net}.

For every pair $(\alpha,\beta)  \in \{ 1,\ldots, r-1\}^2$, we have a
function
\[
n_{\alpha}^{\beta} :\mathcal L_{\alpha} \times \mathcal L_r  \to
\mathcal L_{\beta}
\]
that assigns to $(\ell,\ell') \in \mathcal L_{\alpha} \times
\mathcal L_{r} $ the line in $\mathcal L_{\beta}$ passing through
$\ell \cap \ell'$. Notice that for a fixed $\ell \in \mathcal L_r$
the functions $n_{\alpha}^{\beta}( \cdot,\ell):\mathcal L_{\alpha}
\to \mathcal L_{\beta}$ are bijective.

It follows  from the definition of a $r$-net ({\it cf.} \cite{Ynet}) that
there exists a rational function $F:\mathbb P^2 \dashrightarrow
\mathbb P^1$ of degree $m$ with $r$  values $c_1, \ldots, c_r \in
\mathbb P^1$ for which $F^{-1}(c_i)$ can be identified with
$\mathcal L_i$. Although there is some ambiguity in the definition
of $F$ (we can compose it with an automorphism of $\mathbb P^1$)
the induced foliation is uniquely determined and will be denoted by
$\mathcal F(\mathcal L)$. Similarly, if $(\mathcal L, \mathcal P)$
is a $(r,m)$-net then we will denote   by $\mathcal W(\mathcal L)$ the
CDQL $(m^2+1)$-web $\mathcal W(\mathcal P) \boxtimes \mathcal
F(\mathcal L)$, where $\mathcal W(\mathcal P)$ is the completely decomposable
linear $m^2$-web formed by the superposition of the pencils of lines through the points of $\mathcal P$ .

Among the thirteen sporadic examples of exceptional CDQL webs presented in the Introduction, two are webs
associated to nets. The first one is Bol's web $\mathcal B_5$  which
is associated to a $(3,2)$-net with  $\mathcal P$ equal to  four
points in general position and  $\mathcal L$  equal to the set of lines
joining any two of them. In this case $\mathcal F(\mathcal L)$ is the pencil
of conics through the four points. The other example is the CDQL
$10$-web $\mathcal H_{10}$. It is associated to a $(4,3)$-net with $\mathcal P $ equal
to  the set of base points of the Hesse pencil,  $\mathcal L$ equal
to  the set of lines through any two of them and $\mathcal
F(\mathcal L)$ equal to the Hesse pencil.

The result
below implies that both $\mathcal B_5$ and $\mathcal H_{10}$ are exceptional.

\begin{thm}\label{T:net}
If $\mathcal L$ is a $(r,m)$-net then
\[
{rk}\big(\mathcal W(\mathcal L)\big) \ge \frac{(m^2 -1)(m^2-2)}{2} +
(r-1)^2 -1 \, .
\]
In particular if $\mathcal L$ is a $(3,2)$-net or a $(4,3)$-net then
$\mathcal W(\mathcal L)$ has maximal rank.
\end{thm}
\begin{proof}
Since the $m^2$-subweb $\mathcal W(\mathcal P)$ is linear, it has
maximal rank. To prove the theorem it suffices to show
that
\[
\dim \frac{\mathcal W(\mathcal L)}{\mathcal W(\mathcal P)} \ge
(r-1)^2 -1 \,.
\]

Set $\mathcal L_i = \{ \ell^{(i)}_1, \ldots, \ell^{(i)}_m \}$ and
let $L^{(i)}_j$ be a linear homogenous polynomial in $\mathbb
C[x,y,z]$ defining $\ell^{(i)}_j$. Let  $p_{ij}= \ell^{(1)}_i \cap
\ell^{(r)}_j$ and  $\mathcal L_{ij}$ be the  subset  of $\mathcal L$
formed by the lines through  $p_{ij}$. Notice that $\mathcal P =
\cup_{i,j} \{p_{ij} \} $.

For a suitable choice of the linear forms $L^{(i)}_j$ the  rational function $F : \mathbb P^2 \dashrightarrow \mathbb P^1$
associated to the~arrangement~satisfies
\begin{equation}\label{E:**}
F-c_{\alpha}  = \displaystyle{\frac{ \prod_{i=1}^m L^{(\alpha)}_i} {
\prod_{i=1}^m L^{({r})}_i}} \;
\end{equation}
for every $\alpha < r$.

Let $V = H^0 (\mathbb P^2, \Omega^1_{\mathbb P^2} (
\log \mathcal L))$ (resp. $V_{ij} = H^0 (\mathbb P^2,
\Omega^1_{\mathbb P^2} ( \log \mathcal L_{ij}))$)  be the vector space of logarithmic $1$-forms with
poles in $\mathcal L$  (resp. $\mathcal L_{ij}$). Every element in
$V_{ij}$ vanishes when restricted to the leaves of  the  foliation $\mathcal
L_{p_{ij}}$ induced by the pencil of lines through
$p_{ij}$. If
$\ell_{n_{\alpha}^{\beta}(i,j)}$ denotes the line
$n_{\alpha}^{\beta}(\ell^{(\alpha)}_i,\ell^{(\beta)}_j)$ then  the logarithmic $1$-forms
\[
\frac{dL^{(\alpha)}_{n_1^{\alpha}(i,j)}}{L^{(\alpha)}_{n_{1}^{\alpha}(i,j)}}-
\frac{dL^{(r)}_{j}}{L^{(r)}_{j}}
\]
with $\alpha$ ranging from $1$ to $r-1$, can be taken as a basis of
$V_{ij}$.

 It can be promptly verified
that the union of the subspaces $V_{ij} \subset V$ spans $V$. Indeed
$V$ is generated by elements of the form $\omega= \frac{dL}{L} -
\frac{dL'}{L'}$ where $L, L'$ are linear forms cutting out $\ell,
\ell' \in \mathcal L$. If $\ell \cap \ell'=p_{ij} \in \mathcal P$
then $\omega \in V_{ij}$. Otherwise $\ell$ and $\ell'$ belongs to the
same set $\mathcal L_i$. Then we choose $\ell'' \in \mathcal L_j$, $j
\neq i$, and write
\[
\omega = \left(\frac{dL}{L} - \frac{dL''}{L''}\right)  - \left(
\frac{dL'}{L'} - \frac{dL''}{L''}\right) \, .
\]

It follows from (\ref{E:**}) that the $1$-forms $\frac{dF}{F-c_{\alpha}}$, $\alpha =1, \ldots, r-1$,
 belong to $V$. Therefore there exists $\omega_{ij}^{(\alpha)} \in V_{ij}$ such
 that
 \[
\frac{dF}{F-c_{\alpha}} + \sum_{i,j} \omega_{ij}^{(\alpha)} = 0 \, .
 \]
These equations can be interpreted as elements of $\mathcal
A(\mathcal W(\mathcal L))$. Since the $1$-forms
$\frac{dF}{F-c_{\alpha}}$ are linearly independent, the classes
of these equations span a $(r-1)$-dimensional subspace  $A_0 \subset \frac{\mathcal
A(\mathcal W(\mathcal L))}{\mathcal A(\mathcal W(\mathcal P))}$.

\smallskip

If  $\alpha \neq \beta$ then    $
  \bigcup_{j=1}^{m} \ell_{n_{\alpha}^{\beta}(i,j)} = \mathcal
  L_{\beta}$ for every fixed $i \in \{1, \ldots, m \}$.
Using this fact to work out  the expansion of
$\displaystyle{\frac{dF}{F-c_{\alpha}} \otimes
\frac{dF}{F-c_{\beta}}} $ yields
\begin{eqnarray*}
\displaystyle{\frac{dF}{F-c_{\alpha}} \otimes
\frac{dF}{F-c_{\beta}}} &=  & \underbrace{\sum_{i,j} {\left( \frac{d
  L^{(\alpha)}_i }{ L^{(\alpha)}_i } - \frac{
dL_j^{(r)}}{L_j^{(r)}} \right) \otimes \left( \frac{d
L^{(\beta)}_{n_{\alpha}^{\beta}(i,j)} }{
L^{(\beta)}_{n_{\alpha}^{\beta}(i,j)} } - \frac{
dL_j^{(r)}}{L_j^{(r)}} \right)}}_{\in \bigoplus V_{ij} \otimes
V_{ij}}\\ &+& \sum_{i\neq j} \frac{dL^{(r)}_i}{L^{(r)}_i} \otimes
\frac{dL^{(r)}_j}{L^{(r)}_j} - (r-2)\sum_{i }
\left(\frac{dL^{(r)}_i}{L^{(r)}_i}\right)^{\otimes 2}
 \, .
\end{eqnarray*}

Thus for  ordered pairs  $(\alpha,\beta)$ and $(\gamma,\delta)$ with
distinct entries in $\{1, \ldots, r-1\}$, one have an identity
 \[
 \left({\frac{dF}{F-c_{\alpha}} \otimes
\frac{dF}{F-c_{\beta}}} - \displaystyle{\frac{dF}{F-c_{\gamma}}
\otimes \frac{dF}{F-c_{\delta}}}\right) + \sum_{i,j}
\omega_{ij}^{(\alpha\beta\gamma\delta)} = 0
\]
for suitable $ \omega_{ij}^{(\alpha\beta\gamma\delta)} \in V_{ij}
\otimes V_{ij}$. It follows from Chen's theory of iterated integrals
(see  \cite[Theorem 4.1.1]{Chen} and \cite[Th\'{e}or\`{e}me 2.1]{Robert})
that after integration these identities can be interpreted  as
elements of $\mathcal A(\mathcal L)$.
\begin{footnote}{Notice that  on  the universal covering of $\mathbb P^2 \setminus |\mathcal L|$, one have \\
$
 {\int} \left({\frac{dF}{F-c_{\alpha}} \otimes
\frac{dF}{F-c_{\beta}}} - \displaystyle{\frac{dF}{F-c_{\gamma}}
\otimes \frac{dF}{F-c_{\delta}}}\right)
= \log( F - c_{\alpha})  \frac{dF}{F-c_{\beta}} - \log({F-c_{\gamma}})\frac{dF}{F-c_{\delta}}
$\,.}
 \end{footnote}

Moreover their classes span a
subspace $A_1 \subset \frac{\mathcal A(\mathcal W(\mathcal
L))}{\mathcal A(\mathcal W(\mathcal P))}$ of dimension $(r-1)(r-2)
-1$. Since $A_0 \cap A_1 = 0$, the theorem follows.
\end{proof}

It has to be noted that Theorem \ref{T:net} has a rather limited
scope. Indeed, the Hesse net is the only $r$-net in $\mathbb P^2$ known with $r\ge 4$ and
recently J. Stipins has proved that there {is} no $r$-net
 in $\mathbb P^2$ {if} $r\ge 5$ (for a proof that there is no $r$-net {in} $\mathbb P^2$ {when} $r\ge 6$, see \cite{PY}).
 Nevertheless  Theorem \ref{T:net}
might give some clues on how to  approach the problem about the
abelian relations of  webs associated to hyperplane arrangements
proposed in \cite{PY}. We refer to this paper and the
references therein  for further examples of nets.

The maximality of the rank of $\mathcal H_5$ follows from similar reasons. If $\mathcal L$ is  the Hesse arrangement of lines then
an argument similar to the one used in the proof of Theorem \ref{T:net} shows that
$V={H}^0(\mathbb P^2, \Omega^1(\log \mathcal L))$ can be generated by logarithmic $1$-forms inducing the defining foliations of the
maximal linear subweb of $\mathcal H_5$. Since the Hesse pencil has four linear fibers it follows that
\[
\dim \frac{\mathcal A( \mathcal H_5)}{ \mathcal A\big([(xdy-ydx)\, (dx^3 + dy ^3)]\big) } \ge 3 \, .
\]
Consequently $\mathcal H_5$ has maximal rank.

\subsection{Explicit abelian relations for $\mathcal H_5$}
Alternatively, one can also  establish directly that the rank of $\mathcal H_5$ is maximal.
Indeed the functions  $g_0=(x^3 + y^3 +1)/(xy)$, $g_1=\xi_3\, x+y$, $g_2=x+y$,  $g_3=x+\xi_3\, y$ and $g_4=x/y+y/x$ are  first integrals of $\mathcal H_5$
 and they verify the abelian relations:
\[
\begin{array}{rccccccccc}
  \ln \big( \frac{g_0-3}{g_0-3\,\xi_3}\big)  =&
\hspace{-0.3cm}
\ln \big(\frac{g_1+(\xi_3)2}{g_1+1}  \big)     &\hspace{-0.2cm}  + &\hspace{-0.3cm}
\ln \big(\frac{g_2+ 1}{g_2+\xi_3}\big)
 &\hspace{-0.3cm} + & \hspace{-0.3cm} \ln \big(\frac{g_3+(\xi_3)2}{g_3+1} \big)   &&      \vspace{0.1cm}\\
 \ln \big( \frac{g_0-3\, \xi_3}{g_0-3\,(\xi_3)2}\big)   =&
\hspace{-0.3cm}
\ln \big(\frac{g_1+1}{g_1+\xi_3}  \big)      & \hspace{-0.2cm}+ &
\hspace{-0.2cm}
\ln \big(\frac{g_2+ \xi_3}{g_2+(\xi_3)2}\big)
 &\hspace{-0.3cm} + & \hspace{-0.3cm} \ln \big(\frac{g_3+1}{g_3+\xi_3} \big)   &&     \vspace{0.1cm}\\
\ln \big( \xi_3\,( g_0-3)\big)   = &
\hspace{-0.2cm}
-\ln \big(g_1+(\xi_3)2\big)
     & \hspace{-0.2cm}+ &\hspace{-0.3cm} \ln
\big(\frac{g_22}{1+g_2}  \big)
 &\hspace{-0.3cm} - & \hspace{-0.2cm} \ln (g_3+(\xi_3)2 )  & \hspace{-0.25cm}  - \,
\ln \big({g_4}+\frac{1}{2} \big) \, .
 \vspace{0.1cm}
\end{array}
\]

These   abelian relations  span a 3-dimensional vector space  ${\mathcal A}_1$ such that
$$\mathcal A( \mathcal H_5)=\mathcal A\big([(xdy-ydx)\, (dx^3 + dy ^3)]\big)\oplus {\mathcal A}_1.$$


\section{Abelian relations for the elliptic CDQL webs}\label{S:elliptic}
In this section we will prove that the elliptic CDQL webs presented in the introduction are exceptional.
The  abelian relations of them  that do not come from the maximal linear subweb are all captured
by Theorem \ref{T:ellipticnet} below.

The analogy with Theorem \ref{T:net} is evident and
 probably not very surprising for the specialists  in polylogarithms since, according to the terminology of  Beilinson and Levine {\cite{BL}},
the integrals $\int dz$  and $\int d \log \vartheta$ ($\vartheta$ being a theta function) can be considered  as  elliptic analogs of the classical logarithm and dilogarithm.

\subsection{Rational maps on complex tori with many linear fibers}\label{S:Anet}

Let $T$ be a two-dimensional complex torus  and $F:T \dashrightarrow {\mathbb P}^1$
be a {meromorphic} map. We will say that a fiber $F^{-1}(\lambda)$ is
{\it linear} if it is supported on a union of subtori.

Notice that each subtorus $E$ of $T$ determines a unique linear foliation with $E$ and its translates
being the leaves. We will say that a  linear web $\mathcal W$ {on $T$} supports a fiber $F^{-1}(\lambda)$
if  it contains all the linear foliations determined by the irreducible components
of  $F^{-1}(\lambda)$.

\begin{thm}\label{T:ellipticnet}
Let  ${\mathcal F}$ be the foliation induced by a meromorphic    map $F:T \dashrightarrow \mathbb P^1$.  If $\mathcal W$ is a linear $k$-web with $k\geq 3$ that  supports $m$ distinct linear fibers of $F$,  then
\[
\dim \frac{\mathcal A(\mathcal W\boxtimes \mathcal F)}{\mathcal
A(\mathcal W)} \ge m-1 \, .
\]
\end{thm}

\smallskip

Before proving  Theorem \ref{T:ellipticnet} let us briefly review some basic facts about theta functions. For
details see for instance \cite[Chapitre IV]{Debarre}. If $V$ is a complex vector
space and $\Gamma \subset V$ is a lattice then a theta function
associated to $\Gamma$ is any entire function $\vartheta$ on $V$  such that
for each $\gamma \in \Gamma$ there exists a linear form $a_{\gamma}$
and  a constant $b_{\gamma}$ such that
\[
\quad
\vartheta (z+ \gamma) = \exp\big(2\, i\,\pi\, (a_{\gamma}(z) + b_{\gamma})\big)
\vartheta(z)    \quad \text{for every } z\in V .
\]
Any effective divisor on the complex torus $T= V/\Gamma$ is the zero
divisor of some theta function. Moreover if the divisors of
two theta functions, say $\vartheta$ and $\tilde \vartheta$, coincide then their quotient is a
trivial theta function, that is
\[
\frac{\vartheta}{\widetilde{\vartheta}} ( z ) = \exp\big(P(z)\big)
\]
where $P:V \to \mathbb C$ is a polynomial of degree at most two.

\begin{example}\rm\label{Ex:theta}
If $(\mu , \nu) \in \{ 0,1\} ^2$ and  $\tau \in \mathbb H$ then the entire functions  on $\mathbb C$
\[
\vartheta_{\mu,\nu}(x,\tau) =  \sum_{n=-\infty}^{+\infty}  (-1)^{n\nu} \exp \left(  i\, \pi  \big( n + \frac{\mu}{2} \big)^2 \tau + 2\, i\, \pi  \big( n + \frac{\mu}{2} \big)x  \right)\, .
\]
satisfy the following  relations
\begin{eqnarray}
\label{E:thetamultformula}
\vartheta_{\mu,\nu}(x+ 1,\tau) &=& (-1)^{\mu}\vartheta_{\mu,\nu}(x,\tau)\, \\
 \vartheta_{\mu,\nu}(x+ \tau,\tau) &=& (-1)^{\nu}\exp\big(-i\pi(2z+\tau)\big)\, \vartheta_{\mu,\nu}(x,\tau)\,.
\nonumber
\end{eqnarray}
It is then clear that they are examples of theta functions with respect to the lattice $\mathbb Z \oplus \mathbb Z \tau \subset \mathbb C$.
The theta functions $\vartheta_i$ that appeared in the Introduction are nothing more than
\[
\vartheta_1 = -i\, \vartheta_{1,1} ,\quad \vartheta_2=\vartheta_{1,0}, \quad   \vartheta_3=\vartheta_{0,0} \quad \text{and} \quad   \vartheta_4=\vartheta_{{0},1}.
\]

If $E_{\tau}$ denotes the elliptic curve $\mathbb C / (\mathbb Z \oplus \mathbb Z \tau) $ then the zero divisors of the functions $\vartheta_i = \vartheta_i(\cdot, \tau)$ are
\[
  \big(\vartheta_1\big)_0 = 0, \quad   \big(\vartheta_2\big)_0 = \frac{1}{2} , \quad \big(\vartheta_3\big)_0 = \frac{1+ \tau}{2} \quad \text{and} \quad    \big(\vartheta_4\big)_0 = \frac{\tau}{2} \, .
\]
\end{example}

\smallskip

\begin{proof}[Proof of Theorem \ref{T:ellipticnet}]
With notation  as above, suppose  that  $T = V / \Gamma$. If
$F^{-1}(\lambda)$ is a linear fiber then  one can write
\[
F^{-1}(\lambda) = D_1^{\lambda} + \cdots + D_{r(\lambda)}^{\lambda}
\]
where each divisor $D_i^{\lambda}$ (for $i=1,\ldots,r(\lambda)$) is supported on a union    of translates of
a subtori $E_i^{\lambda}$. Therefore there exist     complex vector spaces  $V_i^{\lambda}$ of dimension one, linear maps
$p_i^{\lambda} :V \to V_i^{\lambda}$ and lattices
$\Gamma_i^{\lambda} \subset V_i^{\lambda}$ such that
\begin{enumerate}
\item $p_i^{\lambda}(\Gamma) \subset \Gamma_i^{\lambda}$;\vspace{0.1cm}
\item$D_i^{\lambda}$ is the pull-back
by the map $[p_i]: T \to
V_i^{\lambda}/\Gamma_i^{\lambda}$
of a divisor on
$V_i^{\lambda}/\Gamma_i^{\lambda}$.
\end{enumerate}
Notice that $p_i^{\lambda}$ can be interpreted as a linear form on
$V$ and its differential $dp_i^{\lambda}$ as a $1$-form defining the
linear foliation determined by $E_i^{\lambda}$.

Composing $F$ with an automorphism of $\mathbb P^1$ we can assume
that the linear fibers are $F^{-1}(\lambda_1), F^{-1}({\lambda_2}),
\ldots, F^{-1}(\lambda_{m-1})$ and $F^{-1}(\infty)$. Thus, for $j$ ranging
 from $1$ to $m-1$,  we can write
\begin{equation}\label{E:theta}
  \displaystyle{  F - \lambda_j = \exp\!\big(P_j(z)\big)\,\frac{ \prod_i [p_i^{\lambda_j}]^* \vartheta_i^{\lambda_j}
  }{ \prod_i [p_i^{\infty}]^* \vartheta_i^{\infty}
  }}
\end{equation}
where the $P_j$'s are  polynomials of degree  at most two
and $\vartheta_i^{\lambda_j}$ are theta functions on
$V_i^{\lambda_j}$ associated to the lattices $\Gamma_i^{\lambda_j}$.
Taking the logarithmic derivative of (\ref{E:theta}),~we~obtain
\begin{equation}\label{E:logtheta}
  \frac{dF}{  F - \lambda_j} = d P_j(z) + \sum_i  [p_i^{\lambda_j}]^* d \log \vartheta_i^{\lambda_j}
   -  \sum_i [p_i^{\infty}]^* d \log \vartheta_i^{\infty} \, .
\end{equation}

Since $\mathcal W$ is a $k$-web with $k\ge 3$  there exist three
pairwise linearly independent linear forms $p_1,p_2,p_3$  among the
$p_i^{\lambda_j}$ such that $dP_j$ can be written as a linear
combination of $dp_1,dp_2, p_1dp_1, p_2dp_2, p_3dp_3$. It follows
that (\ref{E:logtheta}) is an abelian relation for $\mathcal W
\boxtimes \mathcal F$. Since  the logarithmic $1$-forms
$\frac{dF}{  F - \lambda_1},  \ldots, \frac{dF}{  F - \lambda_{m-1}}$ are
linearly independent over $\mathbb C$,  the abelian relations described in (\ref{E:logtheta})
are also linearly independent and generate a subspace of $\mathcal A(\mathcal W \boxtimes \mathcal F)$ of dimension
$m-1$ intersecting $\mathcal A(\mathcal W)$ trivially. The theorem follows.
\end{proof}

In the next three subsections we will derive the exceptionality of the CDQL webs $\mathcal E_5^{\tau},\mathcal E_5,\mathcal E_6$
and $\mathcal E_7$ from Theorem \ref{T:ellipticnet}. Along the way a clear geometric picture of these webs will  emerge.

\subsection{The harmonic $5$-webs $\mathcal E_5^{\tau}$ and the superharmonic $7$-web $\mathcal E_7$}\label{S:ell57}
For $\tau \in {\mathbb H}$,
let $E_{\tau}$ be the elliptic curve $\mathbb C / (\mathbb Z + \tau \mathbb Z)$ and $T_{\tau}$ be the complex torus
$E_{\tau}^2$. For every $\tau \in  \mathbb H$, the $5$-web $\mathcal E_5^{\tau} =  [ dx \,dy \,(dx^2 - dy^2) \,dF_{\tau}]$ is
naturally defined on $T_{\tau}$ where
\begin{equation}
\label{E:Ftau}
  F_{\tau}(x,y) = \left(\frac{\vartheta_1(x,\tau)  \vartheta_1(y,\tau) }{\vartheta_4(x,\tau)  \vartheta_4(y,\tau) }\right)^{\!\!2}  .
\end{equation}

The $7$-web $\mathcal E_7 = [ dx \,dy \,(dx^2 - dy^2)\, (dx^2 + dy^2) \,dF_{1+i}]$ in its turn is naturally defined on $T_{1+i}$.

For every $\alpha, \beta \in \mathrm{End}(E_{\tau})$, denote by $E_{\alpha,\beta}$ the elliptic curve described
by the image of the morphism
\begin{equation}
\label{E:Ealphabeta}
\begin{array}{llcl}
\varphi_{\alpha,\beta} :& E_{\tau} &\longrightarrow& T_{\tau} \\
&x &\longmapsto& (\alpha\cdot  x, \beta\cdot  x) .
\end{array}
\end{equation}
For example $E_{1,0}$ is the {\it horizontal} elliptic curve through $0 \in T_{\tau}$, $E_{0,1}$ is the {\it vertical} one,
$E_{1,1}$ is the {\it diagonal} and $E_{1,-1}$ is the {\it anti-diagonal}. The translation of $E_{\alpha,\beta}$ by an element $(a,b) \in T_{\tau}$
will be denoted by $L_{(a,b)} E_{\alpha,\beta}$.

Let $D_1 = E_{1,0} + E_{0,1}$  and $D_2 = L_{(0,\tau/2)}E_{1,0}  +  L_{ (\tau/2,0)}E_{0,1}$ be divisors in $T_{\tau}$.
Notice that the rational function $F_{\tau}$  is such that $\mathrm{div}(F_{\tau}) = 2D_1 - 2D_2$.
The indeterminacy set of  $F_{\tau}$ is
\[
\mathrm{Indet}(F_{\tau}) =\big\{ (\tau/2,0) , (0,\tau/2) \big\}.
\]
Blowing-up the two indeterminacy points of $F_{\tau}$ we obtain a surface $\widetilde{T_{\tau}}$ containing two pairwise disjoint
divisors $\widetilde{D_1}$, $\widetilde{D_2}$: the strict transforms of $D_1$ and $D_2$ respectively. Let $D_3 = L_{(\tau/2,0)} E_{1,1}   +  L_{(0,\tau/2)} E_{1,-1}$.   The pairwise intersection of the
supports of the divisors $D_1$, $D_2$ and $D_3$ are all equal, that is
\[
|D_1| \cap |D_3| = |D_2| \cap |D_3| = |D_1| \cap |D_2| = \mathrm{Indet}(F_{\tau}) \, .
\]
Therefore $\widetilde{D_3}$, the strict transform of  $D_3$, is a divisor {in $\widetilde{T_{\tau}}$} with support disjoint
from the supports of $\widetilde{D_1}$ and $\widetilde{D_2}$. The lifting of $F_{\tau}$ to $\widetilde{T_{\tau}}$ must map the support of $\widetilde{D_3}$ to $ \widetilde{F_{\tau}}\big(\widetilde{T_{\tau}} - (|\widetilde{D_1}| \cup |\widetilde{D_2}|) \big) =  \mathbb P^1 \setminus \{ 0, \infty\}=\mathbb C^* $. The maximal
principle implies that the image must be a point. Since $\widetilde{D_3}\cdot\widetilde{D_3} =0$, $D_3$ must be a connected component of a fiber of $F_{\tau}$.
Since $D_3$ is numerically equivalent to $2D_1$ and $2D_2$ it turns out that $\widetilde{D_3}$ is indeed a fiber of $F_{\tau}$.
Moreover, because $\widetilde{D_3}$ is connected and reduced, the generic fiber of $F_{\tau}$ is irreducible. In particular, the linear equivalence class of the
divisor $\frac{1}{2} \mathrm{div}(F_{\tau}) =D_1 - D_2$ is a non-trivial $2$-torsion point in $\mathrm{Pic}_0(T_{\tau})$.

So far we have proved that $F_{\tau}$ has at least three linear fibers. For generic $\tau$ it can be verified that $3$ is the exact number of linear fibers of $F_{\tau}$. But if $\tau = 1+i$~then
\[
D_4 = L_{((1+i)/2,0)} E_{1,i} +  L_{ (0, (1+i)/2)}E_{1,-i}
\]
is such that  $|D_1|\cap| D_4 |=| D_2|\cap| D_4 |=| D_3| \cap| D_4| = |D_1 |\cap |D_2|=\mathrm{Indet}(F_{1+i})$. The arguments  above imply that $F_{1+i}$ has
at least $4$ linear fibers.

Theorem \ref{T:ellipticnet} can be applied to the  $5$-webs
$\mathcal E_5^{\tau}= [ dx\, dy\, (dx^2 - dy^2)\, dF_{\tau}]$ (resp. to the $7$-web $\mathcal E_7 = [ dx\, dy\, (dx^2 - dy^2)\,(dx^2 + dy^2)\, dF_{1+ i}]$)
to ensure that ${rk}(\mathcal E_5^{\tau}) \ge 3 + 2 = 5$ (resp. ${rk}(\mathcal E_7) \ge 10  + 3 = 13$).
To prove that $\mathcal E_5^{\tau}$ and $\mathcal E_7$ are exceptional it remains to find two extra abelian relations for the latter web and one for the~former.

The {\it missing} abelian relations  are also captured by Theorem \ref{T:ellipticnet}.
The point is that the torsion of $D_1 - D_2$ is {\it hiding}  two extra linear fibers of $F_{1+i}$ and
one extra linear fiber of $F_{\tau}$. More precisely, since $D_1 - D_2$ is a non-trivial $2$-torsion element of $\mathrm{Pic}_0(T_{\tau})$, there exists a complex torus  $X_{\tau}$, an \'{e}tale covering $\rho:X_{\tau} \rightarrow T_{\tau}$ and a rational function $G_{\tau} : X_{\tau} \dashrightarrow \mathbb P^1$, with irreducible generic fiber, fitting in the  commutative diagram:
\[
\xymatrix{ X_{\tau} \ar@{-->}[d]_{G_{\tau}} \ar[r]^{\rho} & T_{\tau}\ar@{-->}[d]^{F_{\tau}} \\
\mathbb P^1 \ar[r]^{z  \mapsto z^2}  & \mathbb P^1 \; .  }
\]
The \'{e}tale covering above can  be assumed to lift to the identity over the universal covering of $X_{\tau}$ and $T_{\tau}$. In other words,
if $T_{\tau} = \mathbb C^2 / \Gamma_{\tau}$  for some lattice $\Gamma_{\tau} \subset \mathbb C^2$ then $X_{\tau}$ is induced by a sublattice of $\Gamma_{\tau}$. In
particular
$${\rho^* \big[dxdy(dx^k- dy^k)\big]=[dxdy(dx^k- dy^k)\big] \quad \text{for every }  k \ge 1. }$$

Clearly $G_{\tau}$ has at least four linear fibers supported by the linear web {$ [dxdy(dx^2 - dy^2)]$} and $G_{1+i}$ has at least six linear fibers
supported by the linear web {$[dxdy(dx^4 - dy^4)]$}. Theorem \ref{T:ellipticnet} implies that $\rho^* \mathcal E_5^{\tau}$
and $\rho^* \mathcal E_7$ are webs of maximal rank. Since the rank is locally determined the same holds for $\mathcal E_5^{\tau}$
and $\mathcal E_7$.

\subsection{The equianharmonic $5$-web $\mathcal E_5$}

Let $F: T_{\xi_3}\dashrightarrow \mathbb P^1$ be the
rational function
\begin{equation}
\label{E:defFE5}
 F= \frac{\vartheta_1(x, \xi_3)\,\vartheta_1(y, \xi_3 )\,\vartheta_1(x-  y , \xi_3)\,\vartheta_1(x +  \xi_3^2 \,y, \xi_3) }{\vartheta_2(x, \xi_3)\,\vartheta_3(y, \xi_3 )\,\vartheta_4(x- y , \xi_3)\,\vartheta_3(x +  \xi_3^2\, y, \xi_3)} \; .
\end{equation}

\begin{prop}
\label{P:E5}
The  function $F$
has four  linear
fibers on $T_{\xi_3}$. Moreover each of these fibers is supported on the linear $4$-web $\mathcal W = [ dx \, dy \, (dx-dy)\,  ( dx + \xi_3^2\, dy) ]$.
\end{prop}
\begin{proof}
Consider the divisor  $D_1 = E_{1,0} + E_{0,1} + E_{1,1}  + E_{1,-\xi_3}$. Notice that $D_1$ can be given by the vanishing of
\[
f_1(x,y)= \vartheta_1(x, \xi_3)\, \vartheta_1(y, \xi_3 )\, \vartheta_1(x-  y , \xi_3)\, \vartheta_1(x +  \xi_3^2 y, \xi_3).
\]
The complex torus  $T_{\xi_3}$ has sixteen $2$-torsion points and the support of $D_1$ contains
thirteen of them.  The $2$-torsion points that are not contained in  $|D_1|$ are
\[
p_2= \Big( \frac{1}{2}, \frac{1+ \xi_3}{2} \Big) , \quad p_3= \Big( \frac{\xi_3}{2}, \frac{1}{2} \Big) \quad \text{and} \quad
p_4 = \Big( \frac{1+ \xi_3}{2}, \frac{ \xi_3}{2} \Big)  \, .
\]

If we set  $D_i = L_{p_i} D_1$ (the translation of $D_1$ by $p_i$) for $i=2,3,4$,  then
the support of $D_i \cap D_j$  (with $j\neq i$) does not depend on $(i,j)$
and is the set of 12 non-trivial  2-torsion points
of $T_{\xi_3}$ contained in $D_1$. Notice that $D_2$ can be given by the vanishing of
\[
f_2(x,y)= \vartheta_2(x,\xi_3)\, \vartheta_3(y, \xi_3 )\, \vartheta_4(x-  y , \xi_3)\, \vartheta_3(x +  \xi_3^2 y, \xi_3).
\]
The quotient $F=f_1(x,y)/f_2(x,y)$ is the rational function we are interested in.

Blowing up the $12$ indeterminacy points of $F$  one sees that the strict transforms of
the divisors $D_i$ are connected and pairwise disjoint divisors of self-intersection
zero. This is sufficient to prove that each of the divisors $D_i$ is a linear fiber of $F$ and that $F$ has generic fiber irreducible as in the analysis
of the webs $\mathcal E_5^{\tau}$ and $\mathcal E_7$.
Clearly each one
of these fibers is supported on the linear web  $\mathcal W$.
\end{proof}

The proposition above combined with  Theorem  \ref{T:ellipticnet} implies at once that the web $$\mathcal E_5= [ dx \, dy \, (dx-dy)\,  ( dx + \xi_3^2\, dy) ] \boxtimes [dF]$$ is exceptional.

\subsection{The equianharmonic $6$-web $\mathcal E_6$} It remains to analyze the $6$-web
\[
\mathcal E_6 = \big[ dx \, dy\,  (dx + dy) \, (dx + \xi_3\, dy) \, ( dx +\xi_3^2\, dy) \big] \boxtimes \big[ dx/ {\wp(x)} +dy/  {\wp(y)}\big]
\]
on $T_{\xi_3} = E_{\xi_3}^2$.  We will proceed exactly as in the previous cases.

\begin{prop}\label{P:33}
The foliation  $\mathcal F=[ dx/ {\wp(x)} +dy/  {\wp(y)}]$ on $T_{\xi_3}$
admits a rational first integral $F:T_{\xi_3} \dashrightarrow \mathbb P^1$ with generic fiber irreducible and
with three linear fibers, one reduced and two of  multiplicity three. Moreover these three linear fibers are supported on
the linear web $\mathcal W= [ dx \, dy\,  (dx + dy) \, (dx + \xi_3\, dy) \, ( dx +\xi_3^2\, dy) ]$.
\end{prop}
\begin{proof}
Recall that if $\Gamma \subset \mathbb C$ is a lattice then the
Weierstrass $\wp$-function  associated to $\Gamma$ is defined as
\begin{equation}
\label{E:wp}
\wp(z , \Gamma) =\frac{1}{z^2} +  \sum_{\gamma \in \Gamma
\setminus \{ 0 \} } \left( \frac{1}{(z- \gamma)^2} -
\frac{1}{\gamma^2} \right)  .
\end{equation}
It is an  entire meromorphic function with poles of order two on
$\Gamma$ and for a fixed $\Gamma$, the function
$\wp(\cdot,\Gamma)$ is $\Gamma$-periodic, that is
$\wp(\cdot,\Gamma)$ {descends to} a  meromorphic function on the elliptic curve
$E(\Gamma)= \mathbb C / \Gamma$ with a unique pole of order two
at zero.

Recall also that $\wp$ is homogeneous of degree $-2$, that is,
 for any $\mu \in \mathbb C^*$
\begin{equation}\label{E:homogeneo}
\wp(\mu z, \mu \Gamma) = \mu^{-2}{\wp(z , \Gamma)}\, .
\end{equation}

Set  $\Gamma = \mathbb Z
\oplus \mathbb Z \xi_3$ in what follows. Because $\xi_3 \Gamma = \Gamma$,  multiplication
by $\xi_3$ induces an automorphism of $E= E(\Gamma)$, of order $3$  with two   fixed points  besides the origin:
\[
 p_+ = \frac{2 + \xi_3 }{3 } + \Gamma \quad \text{ and } \quad p_-=  \frac{1 + 2\,\xi_3  }{3 } + \Gamma  .
\]

The relation (\ref{E:homogeneo})  implies that
\[
 \wp\big(p_{\pm},  \Gamma\big) = \tau^{-2} {\wp\big(p_{\pm},
 \Gamma\big)} \, .
\]
It follows that  $p_{+}$ and $p_{-}$ are two zeroes of $\wp(\cdot, \Gamma)$.
Since $\wp(\cdot, \Gamma)$ has a unique pole of order two there are
no other zeroes. The points $0, p_+, p_-$ form a subgroup $T$ of the $3$-torsion group
$E(3)$  of $E$.

The $1$-form $\omega=  dx/ {\wp(x)} +dy/  {\wp(y)}     $ is a logarithmic $1$-form with polar set at
$E_{\pm} = \{p_{\pm}\} \times E$ and $E^{\pm} = E \times
\{p_{\pm}\}$. The residues of $\omega$ along  $E_-$ and $E^-$ are equal and so are
those along $E_+$ and $E^+$.  Moreover   the  residue of $\omega $ along  $E_-$ is
 the opposite of its    residue along $E_+$.

The singular set of the foliation $\mathcal F=[\omega]$ consists of five
points: $p_{00}=(0,0)$ and
$$
p_{--}=(p_-,p_-)\, , \quad
p_{-+}=(p_-,p_+)\, , \quad
p_{+-}=(p_+,p_-)\, ,  \quad p_{++}=(p_+,p_+)\; .
$$

The inspection of the first non-zero jet of the closed $1$-form $\omega$ at the singularity $p_{00}$ reveals that ${\mathcal F}$
admits a local first integral analytically equivalent to $x^3 + y^3$ at this point.
The two
singularities $p_{-+},p_{+-}$ are radial ones with local meromorphic
first integrals analytically equivalent to $x/y$. Finally the last two
singularities $p_{--}$ and $p_{++}$ have local holomorphic first integrals
analytically equivalent to $xy$.

If $\alpha \Gamma  = \Gamma$ then
(\ref{E:homogeneo}) implies that
\[
\varphi_{1,\alpha}^* \omega = \frac{ \alpha \, dx}{ \wp( \alpha x,
\Gamma)} + \frac{ dx}{ \wp(x, \Gamma)} = \frac{\alpha^3+
1 }{ \wp( x, \Gamma)} \, dx \; .
\]
A simple consequence is that the  three separatrices of $\mathcal F$ through $p_{00}$ are the elliptic
curves  $E_{(1,-1)}, E_{(1, -\xi_3)}$ and  $E_{(1,-\xi_3^2)}$.

 Let  $\pi:\widetilde{ T_{\xi_3} }  \to T_{\xi_3}$ be the blow-up of $T_{\xi_3}$ at the radial singularities of
$\mathcal F$ and denotes by $\widetilde{ \mathcal F}$ the transformed
foliation.  If
\[
D_1= E_+ + E^+  \, ,  \quad D_2= E_- + E^- \, ,  \quad  \quad D_3 =  E_{(1,-1)} + E_{(1, -\xi_3)} + E_{(1,-\xi_3^2)}
\]
and  $\widetilde{D_1}, \widetilde{D_2},\widetilde{D_3}$  {designate their} respective strict transforms then
\[
\widetilde{D_1}^2 = \widetilde{D_2}^2 = \widetilde{D_3}^2 = 0 \, .
\]

The polar set of $\pi^* \omega$ has  two connected
components, one supported on $|\widetilde{D_1}|$ and the other on $|\widetilde{D_2}|$.
The  divisor $\widetilde{D_3}$  has connected  support  and is disjoint from the polar set of $\pi^*
\omega$. It follows from the  Hodge index Theorem that this divisor is numerically equivalent
to multiples of {$\widetilde{D_1}$}  and {$\widetilde{D_2}$}.
Of course this can be  verified by a direct computation. Indeed,
\[
3\, (E_- + E^-) \equiv 3\, (E_+ + E^+) \equiv
 E_{(1,-1)} + E_{(1, -\xi_3)} + E_{(1,-\xi_3^2)}
\]
where $\equiv$ denotes numerical equivalence.

We can apply
 \cite[Theorem 2.1]{Totaro} (see also \cite[Theorem 2]{jvpjag}) to conclude that the divisors $\widetilde{D_1},\widetilde{D_2}$ and $\widetilde{D_3}$
are fibers of a fibration on $\widetilde{T}$. Consequently   $\mathcal F$ admits a first integral
$F: T_{ \xi_3} \dashrightarrow \mathbb P^1$ satisfying $F^{-1} (0) = 3 D_1$, $F^{-1} (\infty) = 3 D_2$ and $F^{-1} (1) =  D_3$.
As before, $F$ has  generic fiber irreducible because $\widetilde{D_3}$ is connected and  reduced.
\end{proof}

To conclude we proceed as in Section \ref{S:ell57}. The  proof of Propostion \ref{P:33}
shows that the linear equivalence class of $D_1 - D_2$ is a non-trivial $3$-torsion point
of $\mathrm{Pic}_0(T_{\xi_3})$. Therefore there exists a complex torus $X$, an  \'{e}tale covering $\rho: X \to T_{\xi_3}$ and
a rational function $G: X\dashrightarrow \mathbb P^1$ with generic irreducible fiber fitting into the commutative diagram
\begin{equation}
\label{E:G^3=F}
\xymatrix{ X \ar@{-->}[d]_{G} \ar[r]^{\rho} & T_{\xi_3}\ar@{-->}[d]^{F}  \\
\mathbb P^1 \ar[r]^{z  \mapsto z^3}  & \mathbb P^1 \; .  }
\end{equation}
Notice that $G$ has $5$ linear fibers (three of them over $F^{-1}(1)$) and, as in Section \ref{S:ell57}, Theorem \ref{T:ellipticnet}
implies that $\mathcal E_6$ has maximal rank and therefore is exceptional.

\subsection{Explicit abelian relations for elliptic exceptional CDQL webs}
The results of {the four preceding subsections give   geometrical descriptions} of the abelian relations of the  webs $\mathcal E_5^\tau$ (for $\tau \in \mathbb H$), $\mathcal E_5$, $\mathcal E_6$ and $\mathcal E_\tau$. Closed explicit forms for the abelian relations of {these} elliptic exceptional CDQL webs
can {be deduced} from theirs proofs.

\subsubsection{Explicit abelian relations for $\mathcal E_5^{\tau}$.}
\label{AE5explicit}
We recall the description of ${\mathcal A}({\mathcal E}_5^\tau)$ that have been obtained in  \cite{PT}. We fix $\tau \in {\mathbb H}$ and  set $G=F_\tau^{1/2}$ (see
(\ref{E:Ftau})), $g_1=x$, $g_2=y$, $g_3=x+y$ and $g_4=x-y$. Then  the following multiplicative abelian relations hold:
\begin{align}
\label{RAE5}
G  \;= & \;\; \frac{\vartheta_1(g_1,\tau)\, \vartheta_1(g_2,\tau)}{\vartheta_4(g_1,\tau)\,\vartheta_4(g_2,\tau)}
  \vspace{0.1cm} \nonumber  \\
1-G\;  = & \;\;
\frac{\vartheta_3\big(\frac{g_3}{2},\frac{\tau}{2}\big)\,\vartheta_4\big(\frac{g_4}{2},\frac{\tau}{2}\big)}{\vartheta_4(g_1,\tau)\,\vartheta_4(g_2,\tau)      }
 \vspace{0.1cm}\\
1+G \;  = &\;\;
\frac{\vartheta_4\big(\frac{g_3}{2},\frac{\tau}{2}\big)\,\vartheta_3\big(\frac{g_4}{2},\frac{\tau}{2}\big)}{\vartheta_4(g_1,\tau)\,\vartheta_4(g_2,\tau)      }
\; .\nonumber
\end{align}\vspace{0.15cm}

\subsubsection{Explicit abelian relations for $\mathcal E_7$.}
We fix $\tau=1+i$ in this section and we note  $H=F_{\tau}^{1/2}=F_{1+i}^{1/2}$. Let $g_1,\ldots,g_4$ designate  the same functions than above and set $g_5=i\,x+y$, $g_6=x+i\,y$. The   relations
(\ref{RAE5}) of the subweb ${\mathcal E}_5^\tau$ are of course three abelian relations for ${\mathcal E}_7$.
To obtain the last two, we just substitute $ix$ to $x$ in (\ref{RAE5}) and use the transformation formulas for thetas functions admiting complex multiplication (see Section \S8 of  \cite[Chap.V]{Chandrasekharan} for instance)  to get:
\begin{align*}
\label{RAE5}
1-i\,H\;  = & \;\;
\frac{\vartheta_3\big(\frac{g_5}{2},\frac{\tau}{2}\big)\,\vartheta_4\big(i\,\frac{g_6}{2},\frac{\tau}{2}\big)}{\vartheta_4(i\,g_1,\tau)\,\vartheta_4(g_2,\tau)      }
 \vspace{0.1cm}\\
1+i\, H \;  = &\;\;
\frac{\vartheta_4\big(\frac{g_5}{2},\frac{\tau}{2}\big)\,\vartheta_3\big(i\,  \frac{g_4}{2},\frac{\tau}{2}\big)}{\vartheta_4(i\, g_1,\tau)\,\vartheta_4(g_2,\tau)      }
\; .
\end{align*}

\subsubsection{Explicit abelian relations for $\mathcal E_5$.}

To simplify the formulae, we shall abreviate $\xi_3$ by $\xi$, will write $\vartheta_i(z)=\vartheta_i(z,\xi)$ ($i=1,\ldots,4$) and will set $q=e^{i\pi\xi_3}$ in this subsection.  We will also use the notations introduced in the proof of Proposition~\ref{P:E5}.

Let $F$ be the rational function (\ref{E:defFE5}), that is  $F=f_1/f_2$ with
\begin{align*}
 f_1(x,y)=&\, \vartheta_1(x)\,\vartheta_1(y )\,\vartheta_1(x-  y )\,\vartheta_1(x +  \xi^2 y) \\
\mbox{ and } \quad f_2(x,y)= &\,  \vartheta_2(x)\,\vartheta_3(y)\,\vartheta_4(x- y)\,\vartheta_3(x +  \xi^2 y)\; .
\end{align*}
Since $f_1\big(x+\frac{\xi}{2}\,,\,y+\frac{1}{2}\big)
=\, i\,q^{-1/2}e^{i\pi(y-2x)} \vartheta_4(x)
\vartheta_2(y)
\vartheta_3(x-y)
\vartheta_2(x+\xi^2 y)$
 (see \cite[p. 63-64]{Chandrasekharan}),   the linear divisor  $D_3=L_{p_3}(D_1)$ on $T_{\xi}$ is cut out by $$f_3(x,y)= \vartheta_4(x)
\vartheta_2(y)
\vartheta_3(x-y)
\vartheta_2(x+\xi^2 y)\; . $$
One verifies that  $f_3\equiv a_3\,f_1+b_3\,f_2$ where   $a_3=i\frac{  \vartheta_2(0) \vartheta_4(0)}{ \vartheta_3(0)}$ and
$b_3=\frac{  \vartheta_2(0)}{ \vartheta_3(0)} $.
Consequently,  $D_3$ is the linear fiber $F^{-1}(c_3)$ where $c_3=-b_3/a_3=i/\vartheta_4(0)$. According to (the proof of) Theorem \ref{T:ellipticnet},
there is an associated logarithmic abelian relation. Explicitly, it  is (in multiplicative form)
$$ a_3\,F+b_3=\frac{\vartheta_4(x)
\vartheta_2(y)
\vartheta_3(x-y)
\vartheta_2(x+\xi^2 y)   }{\vartheta_2(x)\,\vartheta_3(y )\,\vartheta_4(x- y )\,\vartheta_3(x +  \xi^2\, y)   }\; . $$

In the same way, one proves that  the linear divisor  $D_4=L_{p_4}(D_1)$
is~cut~out~by
$$f_4(x,y)= \vartheta_3(x)
\vartheta_4(y)
\vartheta_2(x-y)
\vartheta_4(x+\xi^2 y)\; . $$
One verifies that  $f_4\equiv a_4\,f_1+b_4\,f_2$ where
$a_4 =i\frac{  \vartheta_2(0) }{ \vartheta_3(0)}$
 and $b_4 =\frac{  \vartheta_4(0)}{ \vartheta_3(0)} $.
So  $D_4=F^{-1}(c_4)$
where $c_4=i\vartheta_4(0)/\vartheta_2(0)$. The associated logarithmic abelian relation is
$$ a_4\,F+b_4=\frac{\vartheta_3(x)
\vartheta_4(y)
\vartheta_2(x-y)
\vartheta_4(x+\xi^2 y)   }{\vartheta_2(x)\,\vartheta_3(y )\,\vartheta_4(x- y )\,\vartheta_3(x +  \xi^2\, y)   }\;. $$

\subsubsection{Explicit abelian relations for $\mathcal E_6$.}
We  shall also abbreviate  $\xi_3$ by $\xi$ in this subsection
and use the notations introduced in the proof of Proposition \ref{P:33}.
Let $\wp(z)$ be the
Weierstrass $\wp$-function (\ref{E:wp}) associated to the lattice
$\Gamma=\mathbb Z\oplus \mathbb Z \,\xi $.
It satisfies the   differential equation
\begin{equation}
\label{E:eqdiffwp}
 \wp'(z)^2=4\,\wp(z)^3-  
\left( \frac{\Gamma({1}/{3})^{3}}{2\,\pi}  \right)^6
\end{equation}
 and  $(\wp)=(p_+)+(p_-)-2(0)$ as divisors  on the elliptic curve $E=\mathbb C/\Gamma$.

We want to make explicit the abelian relations of
\begin{equation*}
   \mathcal E_6=\big[ dx \, dy\,  (dx + dy) \, (dx + \xi\, dy) \, ( dx +\xi^2\, dy) \big] \boxtimes \big[ dx/ {\wp(x)} +dy/  {\wp(y)}\big]
\end{equation*}
defined on $E^2$.
Let $f$ be the elliptic function  defined by
\begin{equation*}
f(x)=\frac{ \wp'(x)-\wp'(p_+)}{\wp'(x) -\wp'(p_-)  }\,.
\end{equation*}
Using (\ref{E:eqdiffwp}), one verifies by a straight-forward computation that $F=f(x)f(y)$ is a first integral  for the  foliation  $[ dx/ {\wp(x)} +dy/  {\wp(y)}]$. We claim that this rational function corresponds exactly to the first integral  deduced in the proof of Proposition \ref{P:33} (also denoted by $F$ there).
One verifies that $(f)=3\,(p_+)-3\,(p_-)$.

Recall (from \cite[Chap. IV]{Chandrasekharan} for instance) the
definition of the Weierstrass sigma function  associated to a lattice $\Lambda \subset \mathbb C$:
$$
\sigma(z,\Lambda)= z\hspace{-0.3cm} \prod_{ \lambda\in \Lambda\setminus \{0\}} \hspace{-0.1cm} \Big( 1- \frac{z}{\lambda}  \Big)\, e^{ \frac{z}{\lambda}+ \frac{z^2}{2\,\lambda^2} } \, .
$$

\begin{lemma}
Let  $\ssigma$ be the Weierstrass sigma function associated to the lattice
$${\Gamma}_1=(2+\xi)\,\Gamma =
 (2+\xi) \mathbb Z\oplus (1+2\,\xi) \mathbb Z\,.$$ If  $E_1=\mathbb C/{\Gamma}_1$ and
\begin{equation*}
g(x)=-\, \frac{  {\ssigma}\big(x-p_+\big) \ssigma\big(x-\xi\,p_+\big) \ssigma\big(x-\xi^2\,p_+\big)   }{ \ssigma\big(x-p_-\big) \ssigma\big(x-\xi\,p_-\big) \ssigma\big(x-\xi^2\,p_-\big)   }
\,. \end{equation*}
then the  product $G=g(x)g(y)$ is a function
that makes commutative the diagram (\ref{E:G^3=F}). More precisely, let $X = E_1^2$  and set $\rho=(\mu,\mu): X\rightarrow E^2$ where $\mu: E_1\rightarrow E$ denotes the isogeny  of degree three induced by the natural inclusion $\Gamma_1 \subset \Gamma$.  Then
\begin{enumerate}
 \item the functions $g$ and $G$ are  rational functions  on $E_1$ and $X$ respectively;\vspace{0.15cm}
\item  they satisfy  $g^3=f\circ \mu$   and   $G^3=F \circ \rho$ on $X$.
\end{enumerate}

\end{lemma}
\begin{proof}
Item {\it(1)} follows at once from  formulae (\ref{E:thetamultformula}). To establish item {\it(2)} one proceeds as usual by comparing the zeroes and the poles of $g^3$ and $ f\circ \mu$ on $E_1$.
\end{proof}

Using the function $G$  one can give closed explicit formulae for the  non-elementary  abelian relations of
\begin{equation*}
\big[ dx \, dy\,  (dx + dy) \, (dx + \xi\, dy) \, ( dx +\xi^2\, dy) \big] \boxtimes \big[ dG\big].
\end{equation*}
The  simplest  is certainly (in multiplicative form)
\begin{equation}
 G=g(x)\,g(y)\, .
\end{equation}
If we set  $g_3=x+y$, $g_4=x+\xi\,y$ and  $g_5=x+\xi^2\,y$ then the other three are
\begin{align}
\label{E:RAE61}
1-G=&  \;  \epsilon_0 \,  \frac{ \ssigma\big(g_3\big)\ssigma\big(g_4 \big)
\ssigma\big(g_5 \big)
}{\prod_{\ell=0}^2 \ssigma\big(x-\xi^\ell p_-\big) \ssigma\big(y-\xi^\ell p_-\big)  }        \\
\label{E:RAE62}
1-\xi\,G=&  \epsilon_1 \frac{ \ssigma\big(g_3+ \xi^2\big)\ssigma\big(g_4+
 \xi \big)
\ssigma\big(g_5+ 1  \big)
}{\prod_{\ell=0}^2 \ssigma\big(x-\xi^\ell p_-\big) \ssigma\big(y-\xi^\ell p_-\big)  } \vspace{0.3cm} \\
\label{E:RAE63}
1-\xi^2\,G=& \epsilon_2 \frac{ \ssigma\big(g_3 -\xi^2\big)\ssigma\big(g_4 -\xi\big)
\ssigma\big(g_5-1 \big)
}{\prod_{\ell=0}^2 \ssigma\big(x-\xi^\ell p_-\big) \ssigma\big(y-\xi^\ell p_-\big)  }\; .
\end{align}
where $ \epsilon_0$, $ \epsilon_1$ and $ \epsilon_2$ are complex constants. Notice that (\ref{E:RAE62}) and (\ref{E:RAE63}) can be obtained from
(\ref{E:RAE61}) by using the relations
$ g(x+1)=g(x+\xi)=\xi \, g(x) $.

\begin{remark}\rm
Since $1-F=(1-G)(1-\xi\,G)(1-\xi^2\,G)$, multiplying (\ref{E:RAE61}), (\ref{E:RAE62}) and (\ref{E:RAE63}) one get a multiplicative abelian relation of $\mathcal E_6$  involving $1-F$. After several simplifications (left to the reader), we find the relation
\begin{align}
\label{E:RA(1-F)1}
1-F=&  \; -\wp'(p_+)\sigma(p_-)^6 \,  \frac{ \sigma(g_3)\sigma(g_4 )
\sigma(g_5)
}{\prod_{\ell=0}^2 \sigma(x-\xi^\ell p_-) \sigma(y-\xi^\ell p_-)  }
\end{align}
where $\sigma$ designates the Weierstrass sigma function associated to the lattice $\Lambda$.
\vspace{0.1cm}

Since $\wp'(x)-\wp'(p_\pm)= \frac{ 2}{\sigma(p_\pm)^3 }  \;  \frac{
\sigma(x-p_\pm) \sigma(x-\xi\,p_\pm) \sigma(x-\xi^2\,p_\pm)
}{   \sigma(x)^3  }$ on $E$, one have also
\begin{align}
\label{E:RA(1-F)2}
1-F
= \; &
\frac{\wp'(p_+) \sigma(p_-)^6}{2}\,
\frac{  \sigma(x)^3 \sigma(y)^3 \,\big( \wp'(x)+\wp'(y)\big)      }
{\prod_{\ell=0}^2 \sigma(x-\xi^\ell p_-) \sigma(y-\xi^\ell p_-)  } .
\end{align}

Comparing (\ref{E:RA(1-F)1}) and (\ref{E:RA(1-F)2}) yields the relation
\begin{equation*}
 -\frac{1}{2}\Big(\wp'(x)+\wp'(y)\Big) =
\frac{ \sigma\big(x+y\big)\sigma\big(x+\xi\,y \big)
\sigma\big(x+\xi^2\,y\big)
}{\sigma(x)^3 \sigma(y)^3 } \, .
\end{equation*}
This is the recently discovered addition formula (6.6) of \cite{EMO}.
\end{remark}


\section{The barycenter transform}\label{S:bary}
\subsection{The $[v]$-barycenter of a configuration}\label{S:vbary} Let $V$ be a two-dimensional vector space over $\mathbb
C$ equipped with a non-zero alternating two-form  $\sigma \in
\bigwedge^2 V^*$. For a fixed  $k\geq 1$ and $v \in V$ {distinct from 0}, consider the map
\begin{equation}\label{E:adef}
\begin{array}{lclcl}
\alpha_v &:& \qquad V^k &\longrightarrow& V \\
 & & (v_1,\ldots, v_k) &\longmapsto& {
\sum_{i=1}^k    \big(\prod_{j\neq i} \sigma(v,v_j)\big)\,           v_i} \, .
\end{array}
\end{equation}

These maps have the following properties:
\begin{enumerate}
\item  $\alpha_v' =
\lambda^{k-1} \alpha_v \, \, $ if $ \, \, \sigma'= \lambda \,\sigma$ {with $\lambda \in \mathbb C^*$};
\item  $\alpha_{\lambda v} = \lambda^{k-1}
\alpha_v$  for every $\lambda \in \mathbb C^*$;
\item   $\alpha_v$ is symmetric;
\item $\alpha_v(v_1,\ldots,
v_k) = 0$ if and only if {there exist   $i$ and $ j$ distinct such that $v_i$,
$v_j$ and $v$ are multiples of each other or if  one of the $v_i$'s is zero}.
\end{enumerate}

\smallskip

 The projectivization  of $\alpha_v$ is a rational map
$\beta_{[v]} : \mathbb P(V)^k
\dashrightarrow \mathbb P(V)$
that admits a nice geometric interpretation: if $[v_i]\neq [v]$
for every $ i \in \{ 1,\ldots, k\}$ then
$\beta_{[v]}([v_1],\ldots,[v_k])$ is nothing but  the barycenter of
$[v_1],\ldots, [v_k]$ seen as points of the affine line $\mathbb C
\cong \mathbb P(V) \setminus \{ [v] \}$. Unlike $\alpha_v$,
$\beta_{[v]}$ does not depend on the choice of $\sigma$. The point
$\beta_{[v]}([v_1],\ldots,[v_k])$ will be referred as the {\it
$[v]$-barycenter of $[v_1],\ldots, [v_k]$}.

\smallskip

The naturalness of  $\beta_v$ is testified  by its
$\mathrm{PSL}(V)$-equivariance, that is,  for every $g \in \mathrm{PSL}(V)$,
$\beta_{gv}( g v_1, \ldots, gv_k ) = g \beta_v(v_1, \ldots, v_k) $.

\subsection{Symmetric versions} Since $\beta_{[v]}$ is a symmetric function it factors
through the natural map $\mathbb P(V) ^k \to \mathbb
P(\mathrm{Sym}^k V)$. Still denoting by $\beta_{[v]}$ the resulting
rational map from $\mathbb P(\mathrm{Sym}^k V)$ to $\mathbb P(V)$, it has been observed in \cite{Frankel} (see also \cite{DMM}) that
$\beta_{[v]}$ admits the affine expression
\begin{equation}\label{E:baryafim}
 \beta_x \big( p(t) \big) = x - k \frac{p(x)}{p'(x)}
\end{equation}
where $x \in \mathbb C$ and  the roots of the degree $k$ polynomial $p
\in\ \mathbb C[t]$
correspond  to  $k$ points in
an affine chart $\mathbb C \subset \mathbb P(V)$.

There are also symmetrized versions of the above maps. Namely we can
define
\[
\begin{array}{lclcl}
\alpha &:& \mathrm{Sym}^k V  &\longrightarrow& \mathrm{Sym}^k V  \\
 & & v_1\cdot v_2\cdots v_k &\longmapsto& \displaystyle{
\prod_{i=1}^k \alpha_{v_i} ( v_1,\ldots, \widehat{v_i}, \ldots, v_k)} \, .
\end{array}
\]
Its   projectivization
\[
\beta: \mathbb P (\mathrm{Sym}^k V) \dashrightarrow \mathbb P
(\mathrm{Sym}^k V)
\]
is a $\mathrm{PSL}(V)$-equivariant rational map.

An affine expression for $\beta$ is also presented in
\cite{Frankel}. If all the  $k$ points belong to the same affine
chart $\mathbb C \subset \mathbb P(V)$ then
\begin{equation}\label{E:symbaryafim}
\beta\left(p(t) \right) = \mathrm{Resultant}_z \big( p(z), (t-z)
p''(z) +2(k-1)p'(z) \big) \,
\end{equation}
where $p \in \mathbb C[t]$ is a degree $k$ polynomial whose  roots
correspond to $k$ points in $\mathbb C \subset \mathbb P(V)$.

\begin{remark}\label{R:dynamics}\rm For
$k=2$, the rational map $\beta:\mathbb
P(\mathrm{Sym}^kV) \dashrightarrow \mathbb P(\mathrm{Sym}^kV)$ is nothing more than the identity map. For $k=3$ it is
 still rather simple: it is  a birational involution
of $\mathbb P^3$ with indeterminacy locus equal to a cubic
rational normal curve. For $k=4$ it is already more interesting from
the dynamical point of view. Recall that for four unordered points
of $\mathbb P^1$ there is a unique invariant, the so
called $j$-invariant. It can be interpreted as  a rational map $j:
\mathbb P ( \mathrm{Sym}^4V) \dashrightarrow \mathbb P^1$ whose
generic fiber contains an orbit of the natural
$\mathrm{PSL}(V)$-action of $\mathbb P( \mathrm{Sym}^4 V)$ as an
open and dense subset. Therefore there exists a rational map
$\beta_*: \mathbb P^1 \to \mathbb P^1$ that fits in the commutative
diagram below:
\[
 \xymatrix{
 \mathbb P ( \mathrm{Sym}^4V) \ar@{-->}[d]_{j} \ar@{-->}[r]^{\beta}  & \mathbb P ( \mathrm{Sym}^4V)\ar@{-->}[d]^{j}
  \\
\mathbb P^1  \ar@{->}[r]^{\beta_*} & \mathbb P^1\,.}
\]

We learned from David Mar\'{\i}n that there exists a choice of coordinates in $\mathbb P^1$ where
\[
\beta_*(z) = \frac{z^2(z+540)^3}{(5z-216)^4} \, .
\]
It can be immediately verified that $\beta_*$ is a post-critically
finite map. We do not know if a similar property holds for the  map
$\beta$ when $k\ge 5$. For a  more comprehensive discussion about
the dynamic of $\beta$ see \cite{Frankel}.
\end{remark}
\subsection{Multiplicities of the first iterate}

Notice that  $\mathbb P (\mathrm{Sym}^k V)$ can be naturally
identified with the set of degree $k$ effective divisors on $\mathbb
P(V)$. In particular, it makes sense to talk about the support of an
element in $\mathbb P(\mathrm{Sym}^k V)$.

\begin{lemma}\label{L:mult}
If  $q_1, \ldots, q_k \in \mathbb P(V)$ are pairwise distinct points
then every point in the support of $\beta(q_1 \cdots  q_k)$ appears
with multiplicity at most  $k-2$.
\end{lemma}
\begin{proof}
Let $q$ be in $S$, where $S$ stands for the  support of
$\beta(q_1,\ldots, q_k)$.  Choose an affine coordinate system in
$\mathbb P(V)$ where $q = 0 $ and $q_i=t_i$ for $i=1,\ldots,k$. If we set $p(t) = \prod
(t-t_i)$ and $p_i(t) = p(t) / (t-t_i)$ then
\[
p_i (t_i ) = p'(t_i) \quad \text{and} \quad p_i'(t_i) = \frac{1}{2} p''(t_i) \, .
\]
Thus, according to (\ref{E:baryafim}), the points in the support of $S$ are of the form
\[
t_i + 2\,(1-k) \,\frac{p'(t_i)}{p''(t_i)} \, .
\]

If $q $ appears in $\beta(q_1 \cdots  q_k)$  with multiplicity at least  $k-1$ then
\[
t_i + 2\,(1-k) \, \frac{p'(t_i)}{p''(t_i)} =0 \implies t_i\,p''(t_i) =
2\,(k-1)\, p'(t_i)
\]
is verified
for at least $(k-1)$ distinct values of  $i$. Since both $p'(t)$ and
$tp''(t)$ are polynomials of degree $k-1$  they
must differ by a nonzero constant. Thus there  exist $\lambda \in
\mathbb C^* $ such that $y(t)=p'(t)$  satisfies the differential  equation
\[
y'(t) = \, \frac{\lambda }{t}\,  y(t)\; .
\]
Hence $p'(t) = C\cdot \exp\left(\frac{\lambda}{2} t^2 \right)$ is not a polynomial. This
contradiction proves the lemma.
\end{proof}


\section{The $\mathcal{F}$-barycenter of a web}
\label{S:fbary}

If $V$ is a two-dimensional vector space over an arbitrary field $F$ of characteristic zero
then it is still possible to define the  $[v]$-barycenter of an element $\mathbb P(\mathrm{Sym}^k V)$.
This can be inferred directly from  equation ({\ref{E:symbaryafim}})    in Section  \ref{S:vbary}.

More explicitly, one {can  specialize}  (\ref{E:baryafim}) to the $\mathcal{F}$-barycenter of a $k$-web ${\mathcal W}$ when
there are at our disposal global rational coordinates $x,y$ on ${S}$.
Assume that ${\mathcal F}=[dx+a\,dy]$ with $a\in {\mathbb C}(S)$. If ${\mathcal W}$ is defined by an implicit differential equation  $F(x,y,dy/dx)=0 $ where $F(x,y,p)$ is  a
polynomial of degree $k$ in $p$ with coefficients in ${\mathbb C}(S)$, then
\begin{equation*}
\label{E:betaF(W)}
\beta_{\mathcal F}\big( {\mathcal W}\big)=\left[
 dx + \bigg(
a-k \, \frac{F(a)}{ \frac{\partial F}{\partial p}(a)}
\bigg)\, dy
\right] .
\end{equation*}

Note also that  the
$\mathrm{PSL}(V)$-equivariance of the barycenter transform yields
\[
 \beta _{\varphi^* \mathcal F} ( \varphi^* \mathcal W) = \varphi^*\big(
 \beta_{\mathcal F} ( \mathcal W ) \big)
\]
for any $\varphi \in \mathrm{Diff}(S)$. Therefore the ${\mathcal F}$-barycenter of $\mathcal W $ is a foliation that is geometrically attached to the pair $({\mathcal F},{\mathcal W})$ and, as such, can be defined on an arbitrary surface by patching together over local coordinate charts the  construction presented above.

\begin{remark}\rm
In \cite{Nakai}, Nakai defines  the {\it
dual $3$-line configuration} of a configuration $L = L_1 \cup L_2
\cup L_3$ of three concurrent lines in the plane: it is ``the unique invariant $3$-line configuration distinct
from $L$ invariant by the group generated by three involutions
respecting  the line $L_i$ and  $L$." The {\it dual $3$-web}
$\mathcal W^*$ of a $3$-web $\mathcal W$ is then defined as the one
obtained by integrating the dual $3$-line configuration of the
tangent $3$-line fields of $\mathcal W$.

It turns out that $\mathcal W^*$ is nothing more than the barycenter
transform of  $\mathcal W)$ in our terminology. Since
$\beta$ is an involution
in the
 case of $3$-points, it follows that $(\mathcal
 W^*)^* = \mathcal W$, a fact already noted by Nakai. Moreover he
 also observed that $K(\mathcal W) = K(\beta(\mathcal
 W))$, see \cite[Theorem 4.1]{Nakai}. In particular, a $3$-web $\mathcal W$ is flat if and only if
 $\beta(\mathcal W)$ is flat \cite[Corollary 4.1]{Nakai}.

We have verified, with the help of a computer algebra system, that the identity  $K(\mathcal W)= K(\beta(\mathcal W))$
also holds for $4$-webs as soon as the four defining foliations of $\beta(\mathcal W)$ are distinct.   For $5$-webs the
situation is different: the barycenter transforms of  most algebraic $5$-webs do not have zero curvature. These  blind constatations
are  crying for  geometric interpretations.
\end{remark}

\subsection{Barycenters of completely decomposable linear
webs }\label{S:riccati}

Let $p_0, \ldots, p_k $  be
$(k+1)$ pairwise distinct points in  $\mathbb P^2$. For any $i=0,\ldots,k$, let  $\mathcal L_i$ denotes the foliation
of $\mathbb P^2$ tangent to the pencil of lines through $p_i$. In
what follows, we give  a description of the foliation $\beta_{\mathcal F}(
\mathcal W)$ when $\mathcal F=\mathcal L_0$ and $\mathcal W =
\mathcal L_1 \boxtimes \cdots \boxtimes \mathcal L_k$.

In the simplest case the points $p_0, \ldots, p_k$ are aligned. If
one chooses an affine coordinate system where all the  $p_i$'s  belong
to the line at infinity then the foliations $\mathcal L_i$ are
induced by constant $1$-forms and so is the $\mathcal F$-barycenter
of $\mathcal W$. The corresponding foliation  $\beta_{\mathcal
F}(\mathcal W)$ is tangent to the pencil of lines through the
point $\beta_{p_0}(p_1, \ldots, p_k)$ on the line at infinity.

If we think of $\mathbb C^2$ as the universal covering of
a two-dimensional complex torus $T$ then, if  $p_0, \ldots, p_k$ are at the line at
infinity, the foliations $\mathcal L_i$ are pull-backs of linear foliations
on $T$ under the covering map.  In this geometric picture the line at infinity is identified with $\mathbb P  H^0(T,\Omega^1_T)$ and the
 the linear foliations $\mathcal L_i$ are defined
by  points $p_i$ in $\mathbb P  H^0(T,\Omega^1_T)$. The $\mathcal F$-barycenter of $\mathcal W$  is the linear foliation on $T$
determined by $\beta_{p_0}(p_1, \ldots, p_k) \in \mathbb P  H^0(T,\Omega^1_T)$.

In the next simplest case $p_1, \ldots, p_k$ are on the same
line  while $p_0$ is not. In an affine coordinate system where
$p_1, \ldots, p_k$ are on the line at infinity and $p_0$ is at the
origin, the $\mathcal F$-barycenter will be induced by the $1$-form
$
\sum_{i=1}^k d \log L_i
$,
where $L_i$ is a linear polynomial vanishing on the line $\overline{
p_0 p_i}$. In particular, {the product}
\begin{equation}\label{E:Li}
\prod_{i=1}^k L_i
\end{equation}
 is a first integral of the foliation $\beta_{\mathcal F}(\mathcal W)$.

In order to  describe the  $\mathcal F$-barycenter of $\mathcal W$ without further restrictions
on the points $p_0, \ldots, p_k$, let  $\Pi: (S,E) \to \mathbb ({\mathbb P}^2,p_0)$ be the blow-up of $p_0$; $\pi: S \to \mathbb P^1$
be the fibration on $S$ induced by the lines through $p_0$;  $\mathcal G$ be the foliation
$\Pi^* \beta_{\mathcal F}(\mathcal W)$; and  $\ell_i$ be the strict transform of the
line  $\overline{p_0p_i}$ under $\Pi$  for $i=1,\ldots,k$.

\begin{lemma}\label{L:ric}
If  the points $\{p_0, \ldots, p_k\}$ are not aligned then
the foliation $\mathcal G$ is a Riccati foliation with respect to $\pi$, that is, $\mathcal G$ has no tangencies
with the generic fiber of $\pi$. Moreover $\mathcal G$  has the  following properties:
\begin{enumerate}
\item the exceptional divisor $E$ of $\Pi$ is $\mathcal G$-invariant;
\item the only fibers of $\pi$ that are $\mathcal G$-invariant  are the lines $\ell_i$, for
$i=1,\ldots, k$;
\item the singular set of $\mathcal G$ is contained in the lines $\ell_i$, for $i=1, \ldots, k$;
\item over each line $\ell_i$ there are two singularities of $\mathcal G$. One is a complex saddle
at  the intersection of $\ell_i$ with $E$, the other is a complex node at the $p_0$-barycenter
of  $\{p_1, \ldots, p_k\} \cap \ell_i$. Moreover,  if $r_i$ is the  cardinality of $\{p_1, \ldots, p_k\} \cap \ell_i$ then
the quotient of eigenvalues of the saddle (resp. node) over $\ell_i$ is $-r_i/k$ (resp.  $r_i/k$);
\item the monodromy of $\mathcal G$ around $\ell_i$ is finite  of order $k/\gcd(k,r_i)$;
\item the only separatrices of $\beta_{\mathcal F}(\mathcal W)$ through $p_0$ are
 the lines $\overline{p_0p_i}$, $i=1,\ldots,k$.
\end{enumerate}
\end{lemma}
\begin{proof}
Let $(x,y) \in \mathbb C^2 \subset \mathbb P^2$  be an affine coordinate system where  $\mathcal F=\mathcal L_0 =[dx]$ (that is  $p_0=[0:1:0]$) and  $\mathcal
L_i = \big[ (x-x_i)\, dy - (y-y_i)\, dx \big]$ (that is $p_i=[x_i:y_i:1]$) for $i=1,\ldots,k$.
It is convenient to assume also  that $y_i \neq 0$ for $i=1,\ldots, k$.

 By definition  $\beta_{\mathcal F}(\mathcal W)$ is
\begin{equation}\label{E:ric}
 \beta_{\mathcal F}(\mathcal W) = \left[   k\,  dy - \Big(\sum_{i=1}^k\frac{y-y_i}{x-x_i}\Big)
 dx\right].
\end{equation}
Since $\Pi:(S,E) \to \mathbb ({\mathbb P}^2,p_0)$ is the blow-up of a point at infinity, {the coordinates $(x,y)$  still define affine coordinates on an affine chart of $S$}. Notice that the fibration $\pi : S \to \mathbb P^1$ is nothing more than  $\pi(x,y)=x$ in these new coordinates.

If we set  $z=1/y$ then $(x,z)\in \mathbb C^2$ is another affine chart of $S$. The intersection of the exceptional divisor
$E=\pi^{-1}(p_0) $ with this chart is equal to $\{ z=0\}$. Notice that  in the new coordinates $(x,z)$ we have
\begin{equation}\label{E:ricbis}
\mathcal G = \left[ k\,dz+z\Big( \sum_{i=1}^k \frac{1-zy_i}{x-x_i}\Big)\,dx \right]  .
\end{equation}
It is  clear from the equations (\ref{E:ric}) and (\ref{E:ricbis}) that:
 $\mathcal G$ has no tangencies with the  generic fiber of $\pi$, that is, $\mathcal G$ is a Ricatti foliation; {\it (1)} the exceptional divisor
$E$ is $\mathcal G$-invariant; {\it (2)} the only $\mathcal G$-invariant fibers of $\pi$ are the lines $\ell_i$; and  {\it (3)} the singularities of $\mathcal G$
are contained in  the lines $\ell_i$.

To prove items {\it (4)} and {\it (5)}  suppose, without loss of generality, that $\ell_1 \cap \{ p_1, \ldots, p_k \} =  \{ p_1, \ldots, p_{r_1} \}$ and that $x_1=0$. In particular
 $x_i \neq 0$ for $i> r_1$. Therefore, in the open set $U= \{ (x,z) \in \mathbb C^2  \, \vert \, \,  | x  | \ll 1 \}$   we can write
 \[
 \mathcal G  = \left[ k x u(x)    \,dz+  \left( z\left( r - \sum_{i=1}^{r_1} zy_i \right)+ zx v(x,z) \right)  \,dx \right]   ,
 \]
where $u$ is an unity in $\mathcal O_U$ that does not depends on $z$ and $v \in \mathcal O_U$ is a regular function. It follows that the singularities of $\mathcal G$
 over $\ell_1= \overline{ \{ x_1=0\}}$ are
 $(0,0)$ and
 \[
  \left(0\,,\frac{r}{\sum_{i=1}^{r_1} y_i} \right)  .
\]
Notice that this last point is the $p_0$-barycenter of $\{ p_1, \ldots, p_{r_1} \}$ on $\ell_1$.

The local  expression for $\mathcal G$ over $U$ also shows that $\mathcal G$ is induced by a vector field $X$ with linear part at $(0,0)$ equal to
\[
   k x \frac{\partial}{\partial x} -    r_1 z \frac{\partial}{\partial z} \, .
\]
Clearly the quotient of eigenvalues in the direction of $\ell_1$ is $-r_1/k$. Since  the points $\{p_0, \ldots, p_k\}$ are not
aligned  $r_1 < k$ and, consequently,  $-r_1/k \in \mathbb Q \setminus \mathbb Z$.     Since $\ell_1$ has zero self-intersection it follows from Camacho-Sad index Theorem that the quotient of eigenvalues (in the direction of the fiber of $\pi$) of the other singularity of $\mathcal G$
on $\ell_1$ is $r_1/k$. Since this number is not an integer
  it follows (see \cite[page 52]{Br}) that this singularity is a complex node.  Moreover the monodromy around $\ell_1$ is analytically conjugated
  to  $z \mapsto \exp(2\pi\, i \,r_1 / k  ) z $. This concludes the proof of {\it (4)} and {\it (5)}.

Finally, to settle {\it (6)} notice that the  singular points of $\mathcal G$ contained in $E$ are  complex saddles. A classical result by Briot and Bouquet
says that these singularities admit exactly two separatrices. In our setup one separatrix corresponds to $E$ and the other corresponds to one of the lines $\ell_i$.
Thus  {{\it (6)} follows} and so does the lemma.\vspace{-0.5cm}
\end{proof}
\begin{center}
\begin{figure}[h]
\psfrag{p0}[][][1]{${p_0} $}
\psfrag{p1}[][][1]{$p_1 $}
\psfrag{p2}[][][1]{$p_2 $}
\psfrag{p3}[][][1]{$p_3 $}
\psfrag{p4}[][][1]{$p_4 $}
\includegraphics[scale=0.7]{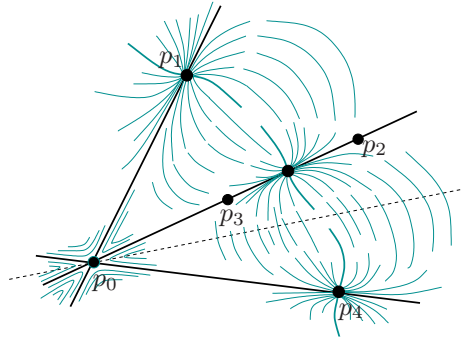}
\caption{The $\mathcal L_{p_0}$-barycenter of the linear web $\mathcal L_{p_1} \boxtimes \cdots \boxtimes \mathcal L_{p_4}$ }
\end{figure}
\end{center}

It is interesting to notice that the generic leaf of $\beta_{\mathcal F}(\mathcal W)$ is
transcendental in general. Indeed, the cases when there are more algebraic
leaves than the obvious ones (the lines $\overline{p_0p_i}$)  are conveniently characterized by the following proposition.

\begin{prop}\label{P:riccati}
The foliation $\beta_{\mathcal F}(\mathcal W)$ has an algebraic leaf
distinct from the lines $\overline{p_0 p_i}$ if and only if all
its singularities distinct from $p_0$ are aligned. Moreover in this
case all its leaves are algebraic.
\end{prop}
\begin{proof}
Since the  Riccati foliation $\mathcal G$ leaves the exceptional divisor $E$ invariant, it  has affine monodromy.
It follows from Lemma \ref{L:ric} item {\it (5)} that its monodromy group is generated by  elements of finite order.

Suppose that  $\mathcal G$ has
an algebraic leaf $L$ distinct from $E$ and the lines $\ell_i$. The existence of such leaf implies that
the monodromy group $G\subset \mathrm{Aut}(\mathbb P^1) $  of $\mathcal G$  must have
a periodic point corresponding to the intersections of $L$ with a generic fiber of $\pi$. Since $G$  already has a fixed point
(thanks to the $\mathcal G$-invariance
of $E$) it follows from Lemma \ref{L:ric} item {\it (5)} that $G$ is conjugated to a  finite subgroup
of $\mathbb C^* \subset \mathrm{Aff(\mathbb C)} \subset \mathrm{Aut}(\mathbb P^1)$. This is sufficient to
show that $\mathcal G$ admits a holomorphic first integral defined on the complement of the $\mathcal G$-invariant fibers of $\pi$.  Lemma \ref{L:ric} item {\it (4)}
implies that $\mathcal G$ is conjugated to  $[r_i ydx -  k xdy]$ in a
neighborhood of  $\ell_i$ and, consequently,  the restriction of $\mathcal G$ to this neighborhood   has  a local meromorphic first integral. Putting all together
it follows that $\mathcal G$ has a global  rational first integral.

Notice that  $\mathcal G$ admits two distinguished leaves that correspond to the two fixed points of the monodromy.
One of these is the exceptional divisor $E$ and the other is an algebraic curve $C$ invariant by $\mathcal G$ such that $\pi_{|C}: C \to \mathbb P^1$ is a
one to~one~covering.

For every $i=1,\ldots, k$, the distinguished leaf  $C$ must intersect the line $\ell_i$ at a singularity
of $\mathcal G$ away from $E$ (by Lemma \ref{L:ric} item {\it (6)}). In a neighborhood of these singularities
$\mathcal G$ has a meromorphic first integral of the form
${y^k}{x^{-r_i}}$ where $r_i$ is the cardinality of $\{p_1, \ldots,
p_k\} \cap \ell_i$ and the local coordinates $(x,y)$  are
 such  that $[dx]$  defines the reference fibration. The restriction of the projection $(x,y) \mapsto x$
to any local leaf not contained in $\{ x \, y =0 \}$ is a
$\frac{k}{\gcd(k,r_i)}$ to $1$ covering of $\mathbb D^*$. Therefore, in these local coordinates around $\ell_i$,
the distinguished leaf $C$ must be contained in $\{ y=0\}$. Notice that the Camacho-Sad index of the leaf $\{y=0\}$ is
$\frac{r_i}{k}$. Summing over the lines $\ell_i$  we obtain from the Camacho-Sad index Theorem that
$ C^2 = 1 \, .$ Since $C$ does not intersects $E$ (Lemma \ref{L:ric} item (6)) it follows that $\Pi(C)$ has self-intersection one.
Thus  $\Pi(C)$ is
a line containing all the singularities of $\beta_{\mathcal F}(\mathcal W)$ different from $p_0$. The proposition follows.
\end{proof}

\begin{cor}\label{C:riccati}
If the foliation $\beta_{\mathcal F}(\mathcal W)$ has an irreducible
algebraic leaf $C$ distinct from the lines $\overline{p_0 p_i}$ then
$C$ is a line or
\begin{equation}
\label{F:C:riccati}
 \deg C = \frac{\sum_{i=1}^m r_i}{\gcd(r_1,\ldots, r_m)}
\end{equation}
where $\{ \ell_1, \ldots, \ell_m\} = \cup_{i=1}^k \overline{p_0
p_i}$ and $r_i$ is the cardinality of $\ell_i \cap \{ p_1, \ldots,
p_k\}$ for $i=1,\ldots,m$.  In particular the degree of $C$ is bounded from below by
$m$.
\end{cor}
\begin{proof}
It follows from  Proposition \ref{P:riccati} that the singularities of $\beta_{\mathcal F}(\mathcal W)$ distinct from $p_0$ are
all contained in an invariant line $\ell$.  We can assume that $\ell$ is the line at infinity
in an affine chart $(x,y) \in \mathbb C^2 \subset \mathbb P^2$. We can also assume that $p_0=(0,0)$
 and, as a by product,  that the $m$ lines
${ \overline{p_0p_i}}$ are cut out by homogeneous linear polynomials $L_1, \ldots, L_m$.
It can easily verified that the polynomial $P = L_1 ^{r_1} \cdots L_m^{r_m}$
is a first integral for $\beta_{\mathcal F}(\mathcal W)$. Of course, if $s_i = r_i / \gcd(r_1,...,r_m)$ then
$( L_1 ^{s_1} \cdots L_m^{s_m} )^{ \gcd(r_1,...,r_m)}  = P $ and therefore $Q= L_1 ^{s_1} \cdots L_m^{s_m}$ is also
 a polynomial first integral for  $\beta_{\mathcal F}(\mathcal W)$. To conclude one has just to observe that
 $Q - c$ is irreducible when $c\neq 0$. Indeed the curve $\overline{\{ Q = c \}} $ is smooth on $\mathbb C^2$ and {has exactly one irreducible branch at each of its
points of  intersection with the line at infinity}.

\end{proof}

\section{Curvature}\label{S:curvatura}
To settle the notation  we recall the
definition of curvature for a completely decomposable $(k+1)$-web $
{\mathcal W}=
\mathcal F_0 \boxtimes \mathcal F_1\boxtimes \cdots \boxtimes
\mathcal F_k$. We start by considering $1$-forms $\omega_i$ with
isolated singularities such that $\mathcal F_i = [\omega_i]$.
For every  triple $(r,s,t)$ with $0\le r<s<t\le k$ we define
\[
\eta_{rst} = \eta(\mathcal F_r \boxtimes \mathcal F_s \boxtimes
\mathcal F_t)
\]
as the unique {meromorphic} $1$-form such that
\[
\displaystyle{\left\lbrace \begin{array}{lcl} d
(\delta_{st}\, \omega_r) &=&
\eta_{rst} \wedge \delta_{st}\, \omega_r \\
d (\delta_{tr}\, \omega_s) &=&
\eta_{rst} \wedge \delta_{tr}\, \omega_s  \\
d (\delta_{rs}\,\omega_t) &=&
\eta_{rst} \wedge \delta_{rs}\, \omega_t  \\
\end{array} \right. }
\]
where $\delta_{ij}=\sigma(\omega_i,\omega_j)$  and $\sigma$ is the
alternating two-form characterized by \[\omega_i \wedge \omega_j =
\sigma(\omega_i,\omega_j)\,  dx \wedge dy .\]

Notice that the $1$-forms $\omega_i$ are not uniquely defined but
any two differ by an invertible function. Therefore, although
dependent on the choice of the $\omega_i$'s, the $1$-forms
$\eta_{rst}$ are well-defined modulo the addition of a closed
holomorphic $1$-form. The {\it curvature} of the web $\mathcal W= \mathcal F_0
\boxtimes \mathcal F_1\boxtimes \cdots \boxtimes \mathcal F_k$ is
thus defined by the formula
\[
K({\mathcal W})=
K(\mathcal F_0 \boxtimes \mathcal F_1\boxtimes \cdots \boxtimes
\mathcal F_k) = d \, \eta({\mathcal W})
\]
where
\[
\eta({\mathcal W})=
\eta(\mathcal F_0 \boxtimes \mathcal F_1\boxtimes \cdots \boxtimes
\mathcal F_k) = \sum_{0\le r<s<t\le k}  \eta_{rst} \, .
\]
Clearly,  $K({\mathcal W})$ is a {meromorphic} $2$-form intrinsically attached to ${\mathcal W}$. More precisely
 for any invertible holomorphic map $\varphi$, one has
\[
K(\varphi^*{\mathcal W})=\varphi^* \big(
K({\mathcal W})
  \big)\, .
\]

We will say that  a $k$-web ${\mathcal W}$ is {\it flat} if its curvature $K({\mathcal W})$ vanishes identically. {This
extends} to every $k\geq 3$ a classical terminology used for a long time for $3$-webs.

\subsection{On the regularity of the curvature}

Our main motivation to introduce the $\mathcal F$-barycenter of a
web  $\mathcal W$ steams  from an attempt to characterize the
absence of poles of $K(\mathcal W)$ over a generic point of an
irreducible component of $\Delta(\mathcal W)$.

In order to state concisely our result in this direction  we introduce the following notation.
If $\mathcal F$ is one of the defining foliations of a $(k+1)$-web $\mathcal W$, then  we define the $k$-web   $\mathcal W - \mathcal F$ by the relation
\[
\mathcal W = ( \mathcal W - \mathcal F) \boxtimes \mathcal F \, .
\]
We will also profit from the usual definition of the tangency between two foliations: if $\mathcal F_1$ and $\mathcal F_2$ are two distinct holomorphic
foliations then $\mathrm{tang}(\mathcal F_1, \mathcal F_2)$ is the divisor locally defined by the vanishing of
\[
   \omega_1 \wedge \omega_2 = 0
\]
where $\omega_i$ are  holomorphic $1$-forms with isolated zeros locally defining $\mathcal F_i $ for $i=1,2$.

\begin{thm}\label{T:curvatura}
Let $\mathcal F$ be a foliation and  $\mathcal W= \mathcal F_1
\boxtimes \mathcal F_2 \boxtimes \cdots \boxtimes \mathcal F_k$ be a completely  decomposable
$k$-web, $k \ge 2$, both defined on the same domain $U \subset
\mathbb C^2$. Suppose that $C$ is an irreducible component of
${\mathrm tang}(\mathcal F, \mathcal F_1)$  that
 is not contained in $\Delta(\mathcal W)$. The curvature
 $K(\mathcal F \boxtimes \mathcal W)$ is holomorphic over a generic point of $C$ if and only if
 the curve  $C$ is invariant by $\mathcal F_1$ or by  $\beta_{\mathcal F_1}(
\mathcal W - \mathcal F_1)$.
\end{thm}

We will need the following lemma.

\begin{lemma}\label{L:curvatura}
If $C$ is an irreducible {component} of ${\mathrm tang}(\mathcal
F,\mathcal F_1)$ that is not contained in $\Delta(\mathcal W)$ then
$\eta(\mathcal F \boxtimes \mathcal W)$ is holomorphic\footnote{{Recall that  $\eta(\mathcal F \boxtimes \mathcal W)$ is defined up to the addition of a closed holomorphic 1-form.
Thus  the holomorphy of $\eta(\mathcal F \boxtimes \mathcal W)$ is well-defined.}}  over the
generic point of  $C$ if and only if $C$ is $\beta_{\mathcal
F_1}(\mathcal W - \mathcal F_1)$-invariant.
\end{lemma}
\begin{proof}
From the hypothesis we can choose a local coordinate system over a
generic point of $C$ such that
\begin{eqnarray*}
\mathcal F &=& \{ \omega_0 = dx + b\, dy =0 \}\, ,  \\
\mathcal F_1 &=& \{ \omega_1= dx =0 \} \\
\mbox{and } \quad \mathcal F_i &=& \{ \omega_i = a_i\,  dx + dy =0 \} \quad
\mbox{for  } \; i=2,\ldots,k\, .
\end{eqnarray*}

A straight-forward computation shows that {for every $i$}:
\[
\eta_{01i} = \frac{ \frac{\partial b}{\partial x} - a_i \frac{
\partial b}{\partial y} -b\left( a_i \frac{\partial b}{\partial x} + \frac{\partial a_i}{\partial y}\right) } {b\,(a_i\,b
-1)} \, \, dx -
 \frac{a_i\, b\,  \frac{\partial b}{\partial y}  + \frac{\partial
a_i}{\partial y} }{ a_i\, b - 1} \, \, dy \, .
\]
Over a generic point of $C$ we have that $C$ coincides with the zero
locus of $b$. Thus  $C$ is not contained in the polar set of
$\sum_{i=2}^k \eta_{01i}$ if and only if  the expression
\[
\sum_{i=2}^k \frac{ \frac{\partial b}{\partial x} - a_i \frac{
\partial b}{\partial y}  } {(a_i\,b
-1)}
\]
is divisible by $b$. But  \[ b \, \, \text{ divides } \, \,
\sum_{i=2}^k \frac{ \frac{\partial b}{\partial x} - a_i \frac{
\partial b}{\partial y}  } {(a_ib
-1)} \, \, \iff \, \, b \, \, \text{ divides } \sum_{i=2}^k
{\Big(} \frac{\partial b}{\partial x} - a_i \frac{
\partial b}{\partial y}{\Big)}\, .
\]
The right hand side above is  {equivalent} to
\[
 b \, \, \text{ divides } \, \, {\left( \Big( \sum_{i=2}^k a_i \Big)\,dx+(k-1)\,dy\right)  \wedge db}  .
\]

From the very definition of the barycenter (see equation (\ref{E:adef})) it follows that
\[
\beta_{\mathcal F_1}(\mathcal W - \mathcal F_1) = \left[ \sum_{i=2}^k \Big(\prod_{
\stackrel{
j=2 }{  j\neq i}}^k \delta_{1i}\Big)\,
\omega_i \right] = \left[ \Big(\sum_{i=2}^k a_i \Big)dx + (k-{1})\, dy \right]  .
\]
Notice that the $1$-form  $(\sum_{i=2}^k a_i )dx + (k-{1}) dy$      has no singularities.
Thus $\sum_{i=2}^k \eta_{01i}$  is holomorphic
on $C$ if and only if  $C$ is invariant by $\beta_{\mathcal F_1}(\mathcal W - \mathcal F_1
)$.

Since $C$ is not contained in $\Delta(\mathcal F_r \boxtimes
\mathcal F_s \boxtimes \mathcal F_r)$ for every set $\{r,s,t\}$ that
does not contain $\{0,1\}$ it follows  that $\eta_{rst}$ is
holomorphic on $C$. The Lemma follows.
\end{proof}

\begin{proof}[{\bf Proof of Theorem  \ref{T:curvatura}}]
In  the notation  of the proof of Lemma \ref{L:curvatura}
\begin{equation}\label{E:soma}
d\omega_0 = \frac{1}{k-1}\left(\sum_{i=2}^k (\eta_{01i} - d \log
\delta_{1i}  )\right) \wedge \omega_0 \; .
\end{equation}

The definition of $\eta(\mathcal W)$ laid down in the beginning of this section implies that
\[
  \sum_{i=2}^k \eta_{01i} = \eta(\mathcal F \boxtimes \mathcal W) - \eta(\mathcal W) \, .
\]
Because $C$ is not contained in $\Delta(\mathcal W)$ both $\eta(\mathcal W)$ and $\sum_{i=2}^k  \eta_{01i}$ are
 holomorphic at the generic point of $C$.

Suppose first that $K(\mathcal F \boxtimes \mathcal W)$ is
holomorphic over the generic point of $C$. If $C$ is $\mathcal F$-invariant then there is nothing to prove.
Thus assume that $C$ is not $\mathcal F$-invariant. If $p$ is a generic point of $C$ and
$\alpha$ is a holomorphic primitive of $d\eta(\mathcal F \boxtimes \mathcal W)$ on a neighborhood of
$p$ then
\[
  \eta(\mathcal F \boxtimes \mathcal W) - \alpha = \frac{d f(b)}{b^n} + dg \,
\]
where $f$ and $g$ are holomorphic functions on a neighborhood of $p$
and $n$ is a positive integer. Therefore (\ref{E:soma}) implies
\[
 d\omega_0 = \frac{1}{k-1}\left( \frac{d f(b)}{b^k} + \alpha'  \right) \wedge
 \omega_0\, ,
\]
for some   holomorphic $1$-form $\alpha'$. Since $d\omega_0$ is
holomorphic and, by assumption,  $\{b=0\}$ is not $\mathcal F$-invariant the only
possibility is that $f\equiv 0$. Therefore  $\eta(\mathcal F \boxtimes \mathcal W)$ is holomorphic along
$C$. It follows from Lemma \ref{L:curvatura} that $C$ is
$\beta_{\mathcal F_1}(\mathcal W - \mathcal F_1)$-invariant.

Suppose now that $C$ is left invariant by  $\mathcal F$ or
$\beta_{\mathcal F_1}(\mathcal W - \mathcal F_1)$. In the latter case the result
follows from Lemma \ref{L:curvatura}. In the former case we can
assume, for a fixed $i \in \{ 2, \ldots, k\}$,  that $C= \{ x=0\}$,   $\omega_0=dx +
x^n u dy$, $\omega_1 = dx$ and $\omega_i = dy$  where $u$ does not vanish identically on $C$. A straight-forward computation shows that
\[
d \eta_{01i} = {\frac {u  {\frac {\partial ^{2} u}{\partial x
\partial y}}  -  {\frac {\partial u }{\partial
y}}   {\frac {\partial u}{\partial x}}
  }{ u ^{2}}} \, .
\]
Thus the $2$-forms
$d\eta_{01i}$ are holomorphic for every $i= 2, \ldots, k$. Because
\[
K(\mathcal F   \boxtimes \mathcal W) = \sum_{i=2}^k d \eta_{01i} + d \eta(\mathcal W) \, ,
\]
and the righthand side is a sum of holomorphic $2$-forms, the curvature $K(\mathcal F \boxtimes \mathcal W)$
 is also holomorphic and the theorem follows.
\end{proof}

\subsection{Specialization to CDQL webs on complex tori}

 Theorem \ref{T:curvatura} completely characterizes in geometric terms
the flat  CDQL webs on two-dimensional complex tori.

\begin{thm}\label{T:curvaturatoro}
Let $\mathcal W= \mathcal L_1 \boxtimes \cdots \boxtimes \mathcal L_k$ be  a linear $k$-web and $\mathcal F$ be a non-linear foliation on a complex torus $T$.
If $k\ge 2$ then $K(\mathcal W \boxtimes \mathcal F)=0$ if and only if  { any} irreducible component  of $\mathrm{tang}(\mathcal F,\mathcal L_i)$  is invariant by $\mathcal F$ or {by}  $\beta_{\mathcal L_i} (\mathcal W - \mathcal L_i)$ for every $i= 1, \ldots, k$ .
\end{thm}
\begin{proof}
Notice that the discriminant of $\mathcal W$ is empty. Therefore the hypotheses of Theorem \ref{T:curvatura} are all satisfied.

If  every irreducible component  of $\mathrm{tang}(\mathcal F,\mathcal L_i)$  is invariant by $\mathcal F$ or $\beta_{\mathcal L_i} (\mathcal W - \mathcal L_i)$ for every $i= 1, \ldots, k$ then Theorem \ref{T:curvatura} implies that
$K(\mathcal W)$ is a holomorphic $2$-form. Since every foliation on $ T$ is induced {by} a global meromorphic $1$-form,  one can proceed as in the beginning
 of this Section to
define a global meromorphic $1$-form $\eta$ on $T$ such that $K(\mathcal W) = {d}\eta$. The result follows from the next proposition.
\end{proof}

\begin{prop}
\label{P:Brunella}
Let $\omega$ be a meromorphic $1$-form on  a
compact  K\"{a}hler manifold
$M$. If $d \omega$ is holomorphic then $\omega$ is closed.
\end{prop}
\begin{proof} We learned the following proof  from Marco Brunella.
 Notice that although $\omega$ is not closed a priori, the holomorphicity of $\Omega=d\omega$ ensures
that  its residues along the irreducible components $Z_i$ of its polar set   are well-defined complex numbers.
If $S$ is {a} real subvariety of $M$ of real
{dimension} $2$, then Stoke's Theorem implies that
\[
  \int_S \Omega = \sum_{i=1}^m {\rm res}_{Z_i}(\omega) \cdot \big( S \cdot Z_i\big)  \,
\]
where $S \cdot Z_i$ stands for the topological intersection number of $S$ with $Z_i$. It follows that the class of $\Omega$, seen as a current, lies in ${H}^{1,1}(M,\mathbb C)$.

On the other hand,  $\Omega$ being a closed holomorphic 2-form, its class lies also in ${H}^{{2,0}}(M,\mathbb C)$. But ${H}^{1,1}(M, \mathbb C) \cap {H}^{{2,0}}(M,\mathbb C) = 0$
since $M$  is K\"{a}hler. This implies that  $\Omega$ is zero and consequently $\omega$ is   closed.
\end{proof}

Theorem \ref{T:curvaturatoro} admits the following   consequence.

\begin{cor}\label{C:w-f}
Let $\mathcal W$ be a linear $k$-web and $\mathcal F$ be a foliation both defined on the same complex torus $T$.
Suppose that $\mathcal W$ decomposes as $\mathcal W_1 \boxtimes \mathcal W_2$ in such a way that
$\mathcal W_1$ and $\mathcal W_2$ are not foliations. Suppose also that for every defining foliation
$\mathcal L$ of $\mathcal W_i$, $i=1,2$, we have
\[
  \beta_{\mathcal L} (\mathcal W_i - \mathcal L) =   \beta_{\mathcal L} (\mathcal W - \mathcal L) \, .
\]
Then  $K(\mathcal W \boxtimes \mathcal F)=0$ if and only if $K(\mathcal W_i \boxtimes \mathcal F)=0$ for $i=1,2$.
\end{cor}

\begin{example}\label{E:PT}\rm
Consider the linear $4$-web $$\mathcal W= \underbrace{[dxdy]}_{\mathcal W_1} \boxtimes \underbrace{[(dx-dy)(dx+dy)]}_{\mathcal W_2} $$ on a two-dimensional complex torus $T$.
Notice that
\[
 \beta_{[dx]}(\mathcal W) = [dy] = \beta_{[dx]}(\mathcal W_1) \quad \text{and} \quad  \beta_{[dy]}(\mathcal W) = [dx] = \beta_{[dy]}(\mathcal W_1) \, .
\]
Similarly
$
 \beta_{[dx\pm dy]}(\mathcal W) = [dx\mp dy] = \beta_{[dx\pm dy ]}(\mathcal W_2)  .
$

In \cite{PT}, germs of exceptional CDQL $5$-webs on $(\mathbb C^2,0)$ of the form $$[dxdy (dx-dy)(dx+dy)]\boxtimes \mathcal F$$ are classified
under the additional assumption that $K([dxdy]\boxtimes \mathcal F)=0$. Mihaileanu's criterion combined with the Corollary \ref{C:w-f} above yields  that the additional assumption
is superfluous if $\mathcal F$ is supposed to be globally defined on {a} complex torus  $T$. Translating the classification of \cite{PT} to our setup,
we obtain that  every flat and  global $5$-web on complex tori of the form $[dxdy (dx-dy)(dx+dy)] \boxtimes \mathcal F$ is isogeneous to
one of the $5$-webs $\mathcal E_{\tau}$ ($\tau \in \mathbb H$) presented in the Introduction. In particular the torus $T$ has to be isogeneous to the square of an elliptic curve.
\end{example}

\subsection{Specialization to CDQL webs on the projective plane} It
would be interesting to extend Theorem \ref{T:curvatura} in order to
deal with more degenerated discriminants. We do not know how to do
it in general. Nevertheless under the assumption that  $\mathcal W$
is a product of linear foliations on the projective plane we have the following weaker result.

\begin{thm}\label{TC:curvatura}
Let $\mathcal F$ be a foliation and  $\mathcal W= \mathcal L_1
\boxtimes \mathcal L_2 \boxtimes \cdots \boxtimes \mathcal L_k$ be a
totally decomposable linear $k$-web, $k \ge 2$, both globally defined on
$ \mathbb P^2$. Suppose that $C$ is an
irreducible component of ${tang}(\mathcal F, \mathcal L_1)$.
{If
 $K(\mathcal F \boxtimes \mathcal W)$ is holomorphic over a generic point of $C$ then
the curve $C$ is invariant by $\mathcal L_1$ or by $\beta_{\mathcal L_1}(
\mathcal W - \mathcal L_1)$}.
\end{thm}
\begin{proof}
If $C$ is not contained in $\Delta(\mathcal W)$ then the result
follows from Theorem \ref{T:curvatura}. Thus, assume that $C \subset
\Delta(\mathcal W)$. The tangency of two linear foliations on $\mathbb P^2$ is
a line invariant by both and, therefore,  $C$ must be a line invariant by at least
two of the defining foliations of $\mathcal W$.

If $C$ is $\mathcal L_1$-invariant
then there is nothing to prove. Thus assume that this is not the case. Because $C \subset \mathrm{tang}(\mathcal F, \mathcal L_1 )$,  we are also
assuming that $C$ is not $\mathcal F$-invariant.

First remark that  Theorem \ref{T:curvatura} implies that  $K(\mathcal F
\boxtimes \mathcal L_i \boxtimes \mathcal L_j)$  is holomorphic over
the generic point of $C$ for every choice of distinct $i, j \in \{
2, \ldots, k\}$. Indeed, on the one hand  if $C \subset \mathrm{tang}(\mathcal F, \mathcal L_i)$ then
$\mathcal L_i$ and $\mathcal L_1$ have to be tangent along $C$. Thus $C$ is $\mathcal L_1$-invariant
contrary to our assumptions. On the other hand if $C \subset \mathrm{tang}(\mathcal L_i, \mathcal L_j)$
then $C$ is invariant by both $\mathcal L_i$ and $\mathcal L_j$ and the triple  $(\mathcal F, \mathcal F_1, \mathcal W)=(\mathcal L_i,\mathcal L_j,\mathcal F  \boxtimes \mathcal L_j)$ satisfies the hypotheses of  Theorem \ref{T:curvatura}. Thus $K(\mathcal F
\boxtimes \mathcal L_i \boxtimes \mathcal L_j)$  is indeed holomorphic over
the generic point of $C$.

Similarly, Theorem \ref{T:curvatura}  implies that   $K(\mathcal F \boxtimes
\mathcal L_1 \boxtimes \mathcal L_i)$ is holomorphic along $C$ whenever $C$ is $\mathcal
L_i$-invariant.

If we write  $\mathcal W = \mathcal L_1 \boxtimes
\mathcal W_0 \boxtimes \mathcal W_1$ with  $\mathcal W_1$ being the
product of foliations in $\mathcal W$ leaving $C$ invariant and $\mathcal W_0$ being the product of foliations in
$\mathcal W - \mathcal L_1$ not leaving $C$ invariant then $K(\mathcal F
\boxtimes  \mathcal L_1 \boxtimes \mathcal W_0 )$  is holomorphic over the generic point of $C$.

Because $C$ is not contained in $\Delta( \mathcal L_1 \boxtimes
\mathcal W_0)$, Theorem \ref{T:curvatura} implies that $C$ is
$\beta_{\mathcal L_1}( \mathcal W_0)$-invariant. From the
definition of the $\mathcal L_1$-barycenter it follows that $C$ is
{also  invariant by $\beta_{\mathcal L_1}( {\mathcal
W_0 \boxtimes \mathcal W_1} )=\beta_{\mathcal L_1}( {\mathcal
W- \mathcal L_1} ) $}. \end{proof}

\medskip

Notice that in Theorem  \ref{TC:curvatura},  unlikely in Theorem \ref{T:curvatura}, the invariance condition imposed on
$C \subset \mathrm{tang}(\mathcal F, \mathcal L_1)$ is no longer a  necessary and sufficient condition for
the regularity of the curvature: it is just necessary.
In fact, the converse to  Theorem \ref{TC:curvatura} does not
hold in general. For instance, if $\mathcal F= [ y dx +
dy]$, $\mathcal L_1 = [ dy]$ and $\mathcal L_2 = [ydx - x dy]$, then
the line $L=\{ y=0 \}$ is invariant by $\mathcal F$, $\mathcal L_1$ and
$\beta_{\mathcal L_1}(\mathcal L_2)= \mathcal L_2$ but $ K(\mathcal F \boxtimes \mathcal L_1 \boxtimes \mathcal L_2)$ is not holomorphic over $L$ since
\[
 K(\mathcal F \boxtimes \mathcal L_1 \boxtimes \mathcal L_2) =
 \frac{dx \wedge dy}{y\, (x+1)^2} \, .
\]


\section{Constraints on flat CDQL webs}
\label{S:constraints}
In this section we start the classification of exceptional CDQL webs on the
projective plane. As already mentioned in the {I}ntroduction the
starting point is Mih\u{a}ileanu criterion:
{\it If $\mathcal W$ is a web  of maximal
rank then $K(\mathcal W)=0$.}

We will  combine this criterion with Theorem
\ref{TC:curvatura} in order to restrict the possibilities for the
pairs $(\mathcal F, \mathcal P)$. For instance,
Theorem \ref{T:bd} below shows that the degree of $\mathcal F$ is bounded by four  when  $\mathcal F \boxtimes \mathcal
W(\mathcal P)$ is flat.

Here, as usual,  the degree of
a holomorphic foliation $\mathcal F$ on $\mathbb P^2$ is the number of tangencies with a
generic line $\ell \subset \mathbb P^2$.  Concretely,  in affine coordinates
$(x,y) \in \mathbb C^2 \subset \mathbb P^2$, a foliation $\mathcal F$ has degree
$d$ if and only if {$\mathcal F$ is defined by}  a polynomial $1$-form $\omega$ with isolated zeros
that can be written
in the form
\[
\omega = a(x,y) dx + b(x,y) dy + h(x,y)(xdy - ydx)
\]
where $h$ is a homogeneous polynomial of degree $d$; $a$ and $b$ are polynomials
of degree at most $d$ and;  when $h$ is the  zero   polynomial   the polynomial $xa + yb$  has degree exactly $d+1$.

We point out that {$h$ vanishes identically} if and only if  the line at infinity is {$\mathcal F$}-invariant. In this case the zeros of the homogenous component of degree $d+1$ of the polynomial $xa + yb$ correspond to the singularities of $\mathcal F$ on the line at infinity. If $h$ is non-zero then the points at infinity determined by $h$
are in one to one correspondence with the tangencies of $\mathcal F$ with the line at infinity.

\subsection{Notations}

The notations below  will be used in the proof of the
classification of {exceptional} CDQL webs on the projective plane.

\smallskip
\begin{tabbing}
aaaaaaaaaaaa \= Description \kill
$\mathcal P$ \>  finite set of points in $\mathbb P^2$;\\
$k$ \> the cardinality of $\mathcal P$; \\
$p_1,\ldots, p_k$ \> the points of $\mathcal P$; \\
$\mathcal P_i $ \> $\mathcal P \setminus \{ p_i \}$;  \\
$\mathcal L_i  $ \> the linear foliation determined by $p_i$; \\
$\mathcal W(\mathcal P)$ \> $\mathcal L_1 \boxtimes \cdots \boxtimes \mathcal  L_k $; \\
$\mathcal W(\mathcal P_i)$ \>  $\mathcal W(\mathcal P) - \mathcal L_i$; \\
$\widehat{\mathcal L_i}$ \> the $\mathcal L_i$-barycenter of $\mathcal W(\mathcal P_i)$ that is,  $\beta_{\mathcal L_i} \big( \mathcal W(\mathcal P_i) \big)$; \\
$\mathcal L_p$ \> the pencil of lines through a point $p \in \mathbb P^2$;\\
$\widehat{\mathcal L_p}$ \> in case $p \in \mathcal P$, the $\mathcal L_p$-barycenter of $\mathcal W(\mathcal P \setminus \{p \} )$; \\
$ \ell$ \>  a line on $\mathbb P^2$; \\
$\mathcal P_{\ell}$ \> $\mathcal P \cap \ell$; \\
$k_{\ell}$ \> the cardinality of $\mathcal P_{\ell}$; \\
$q_1, \ldots, q_{k_{\ell}}$ \> the points of $\mathcal P_{\ell}$; \\
$\widehat{q_{i}}$ \> the $q_i$-barycenter of $\mathcal P_{\ell} \setminus \{ q_i\}$ in $\ell$.
\end{tabbing}
\smallskip

A set $\mathcal P$ is  in
{\it $p_i$-barycentric general position} if the only algebraic leaves of
$\widehat{\mathcal L_i}$ are the lines $\overline{p_i p_j}$ (compare
with Proposition \ref{P:riccati}). We will
{write $b(\mathcal P)$ for the cardinality of}
\[
\mathcal B(\mathcal P) = \big\{ p \in \mathcal P  \, \vert \, \mathcal P
\text{ is in } p\text{-barycentric general position} \big\} \,.
\]

\subsection{Configurations of points in barycentric general position}
As an  immediate consequence of Theorem \ref{T:curvatura} it follows
that a completely decomposable $3$-web  $\mathcal W= \mathcal
F\boxtimes \mathcal L_1 \boxtimes \mathcal L_2$  on $\mathbb P^2$
induced by two pencils of lines and a foliation has  curvature zero
if and only if it is projectively equivalent to a web of the form
\[
 \big[ a(y) dx + b(x) dy\big] \boxtimes [dx] \boxtimes [dy]
\]
where $a$ and $ b $ are rational functions.

In the same vein, the next result  combines Proposition
\ref{P:riccati} with Theorem \ref{TC:curvatura} to show how {\it
generic} configurations of points impose strong restrictions on a
foliation $\mathcal F$ when $\mathcal F \boxtimes \mathcal
W(\mathcal P)$  has curvature zero.

\begin{prop}\label{P:gp}
Let $\mathcal W(\mathcal P)$ be  the $k$-web naturally associated
to  a collection ${\mathcal P}$ of $k$ distinct points in $\mathbb P^2$. If  $\mathcal F$ is a non-linear foliation on $\mathbb P^2$
such that $K\big(\mathcal F \boxtimes \mathcal W(\mathcal P)\big)=0$ then
$b(\mathcal P)$ is at most $4$. Moreover, there exist affine
coordinates $x,y$ such that
\begin{enumerate}
\item[(a)] if $b(\mathcal P) =1$ then $\mathcal F= [ a(y) dx + b(x,y)\,
dy]$ for some $a \in \mathbb C[y]$, $b\in \mathbb C[x,y]$;\vspace{0.15cm}
\item[(b)] if $b(\mathcal P) =2$ then $\mathcal F= \big[ a(y)\, dx + b(x)\, dy\big]$ for some $a,b \in \mathbb C[t]$;\vspace{0.15cm}
\item[(c)] if $b(\mathcal P) =3$ then the points in $\mathcal B(\mathcal P)$ are not aligned and
 $$\mathcal F= \big[ y\,( y^{d-1} - \epsilon_1) \,dx - x\, ( x^{d-1} - \epsilon_2)\, dy\big]$$
for some integer $ d \geq 2$ and $\epsilon_1, \epsilon_2 \in \{
0,1\}$ or
 $$\mathcal F= \big[ y\, dx - \lambda \,x\, dy\big]$$
for some constant $\lambda \in {{\mathbb C} \setminus \{0, 1\}} $;\vspace{0.15cm}
\item[(d)] if $b(\mathcal P) =4$ then the points in $\mathcal B(\mathcal
P)$ are in general position and $\mathcal F$ is the pencil of conics
through them.
\end{enumerate}
\end{prop}
\begin{proof}Suppose that $\mathcal P$ is in $p_1$-barycentric
general position and assume that $p_1=[0:1:0]$. If $K\big(\mathcal F \boxtimes \mathcal W(\mathcal P)\big)=0$
then
Theorem \ref{TC:curvatura} implies that the tangency between $\mathcal F$
and $\mathcal L_1$ is a union of lines through $p_1$. In the affine
coordinates $(x,y)= [x:y:1]$,   $\mathcal L_1= [dy]$ and the lines through $p_1$ correspond to
vertical lines. Therefore $
 \mathcal F = [ a(y) dx + b(x,y) dy ]
$ for some polynomials $a,b$.

If $\mathcal P$ is also in $p_2$-barycentric general position and
$p_2 = [1:0:0]$, the same argument shows that $
 \mathcal F = [ a(y) dx + b(x) dy ]$ for some  polynomials~$a,b$.

Notice that no point  $p \in \mathcal B(\mathcal P) \setminus \{ p_1, p_2 \}$ can be aligned with $p_1$ and $p_2$.
Indeed, suppose the contrary. One can assume that  $p=p_3=[1:1:0]$, or  equivalently $\mathcal L_3 = [dx - dy]$.
Then the tangency  of  $\mathcal F$ and
$\mathcal L_3$ is given by vanishing of
\[
(dx - dy ) \wedge \big( a(y) dx + b(x) dy \big) = \big( b(x) + a(y) \big) dx \wedge dy  \, .
\]
Because   $K\big(\mathcal F \boxtimes \mathcal W(\mathcal P)\big)=0$,  $\{ b(x) + a(y)  = 0\}$  must
be a  union of lines through $p_3$. Explicitly, up to a multiplicative constant,
\[
b(x) + a(y) =  \prod_{j=1}^m ( x-y - c_j )
\]
for suitable constants $c_1,\ldots,c_m$. Such identity is possible if and only if  the homogenous component of higher order of
$a(y) dx + b(x) dy$ is a constant multiple of $x\,dy-y\,dx$. Therefore  $\mathcal F$ has degree zero and consequently is a linear foliation.
This contradicts our assumptions on $\mathcal F$.

Suppose now that $\mathcal P$ is also in $p_3$-barycentric general
position with $p_3 \notin \overline{p_1p_2}$. It is harmless to assume
 that $p_3 = [0:0:1]$.
Since the tangency of $\mathcal F$ and  $\mathcal L_3$ is a union
lines through $p_3=(0,0) \in \mathbb C^2$ then  the polynomial $xa(y) + yb(x)$ must be  homogeneous. Thus
for a certain  $d\in {\mathbb N}^*$ and suitable $c_0, c_1, c_2 \in
 \mathbb C$
\[
 (a(y) , b(x) ) = ( c_1 y^d + c_0 y, c_2x^d - c_0 x)  \, .
\]
It is a simple matter to show that we are in one of the two cases
displayed in part {\it  (c)} of the statement, the first when $d\ge 2$ and
the second when $d=1$.

Finally, suppose that  $b({\mathcal P})\geq 4$.  Since  no three  points in $\mathcal B(\mathcal P)$ are  aligned  we can assume that  $p_1,p_2,p_3$ are as above and $p_4=[1:1:1]$.   Applying again the above
argument to $\mathcal L_4$ and  discarding the solutions corresponding to degree zero foliations, we
 prove that
 \[
{\mathcal F}
  = \big[ a(y) dx + b(x) dy\big] = \big[  y(y-1)\, dx - x(x-1)\, dy
  \big]\, .
 \]
Notice that the rational function $\frac{x(y-1)}{y(x-1)}$ is a first integral of $\mathcal F$, that is, $\mathcal F$ is a pencil
of conics through the four points $p_1, \ldots, p_4$. Notice also  that
$\mathcal F$ leaves invariant exactly six lines: the line at infinity and the five affine lines cut out by the polynomial $xy(x-1)(y-1)(x-y)$.
If  ${\rm tang}({\mathcal L}_p,{\mathcal F})$ is a union of lines through $p$ then $p$ must belong to three of  the six $\mathcal F$-invariant lines. Since there are only four such points ($p_1, p_2, p_3$ and $p_4$)
$b({\mathcal P})$ has at most four elements. This concludes the proof. \end{proof}

\begin{cor}\label{C:onaline}
Assume that  the cardinality of $\mathcal P$ is at least $4$.  If it exists
a non-linear foliation ${\mathcal F}$ such that $K({\mathcal F}\boxtimes{\mathcal W}({\mathcal P}))=0$ then one of the following two situations occurs:
 \begin{enumerate}
\item  there are three  aligned  points in ${\mathcal P}$; \vspace{0.15cm}
\item ${\mathcal P}$ is the union of  of 4 points in general position and ${\mathcal F}$  is the pencil of conics
through them.
 \end{enumerate}
 \end{cor}
\begin{proof} Assume that we are not in case {\it (1)}:
any line contains at most two points of ${\mathcal P}$. Lemma \ref{L:ric} item  {\it (4)} implies that the singularities of $\widehat{\mathcal L_p}$ coincide
with $\mathcal P \setminus \{ p\}$ for any $p \in \mathcal P$. By assumption, the set of points $\mathcal P \setminus \{ p\}$  is not aligned and, according to
Proposition \ref{P:riccati}, $\mathcal P$ is in $p$-barycentric general position. Thus   ${\mathcal P}={\mathcal B}({\mathcal P})$  and
 Proposition  \ref{P:gp} implies the result.
\end{proof}

\subsection{Aligned points versus invariant lines}

Non-generic
configurations of points also impose non-trivial conditions  on  non-linear foliations $\mathcal F$
such that the curvature of
$\mathcal F \boxtimes \mathcal W(\mathcal P)$
vanishes identically.

\begin{prop}\label{P:invariante}
Let $\mathcal P \subset \mathbb P^2$ be a set  of $k$
points and $\mathcal F$ be a non-linear foliation on $\mathbb P^2$  such that $K(\mathcal F \boxtimes \mathcal W(\mathcal P))=0$.
If   $\ell$ is a line that contains  at least three points of
$\mathcal P$ then $\ell$ is $\mathcal F$-invariant.
\end{prop}
\begin{proof}
Remind that $k_{\ell} = {\rm Card}(\mathcal P_{\ell})$ with $\mathcal P_{\ell}=\mathcal P\cap \ell = \{q_1, \ldots, q_{k_{\ell}} \}$. By hypothesis, $k_{\ell}\ge 3$.
If $\ell$ is not invariant by $\mathcal F$ then
\begin{equation}\label{E:l}
|{\rm tang}(\mathcal F, \ell)| \subset |{\rm tang}(\mathcal F,
\mathcal L_i)| \cap \ell
\end{equation}
for every $i=1,\ldots,
k_{\ell}$, since $\ell$ is invariant by $\mathcal L_i=\mathcal L_{q_i}$.

Notice that for every $i$ ranging from 1 to $k_{\ell}$, the Riccati foliation
 $\widehat{ \mathcal L_i} $ leaves $\ell$ invariant and its
singularities on $\ell$ are $q_i$ and $\widehat{q_i}$ according to Lemma \ref{L:ric} items {\it (2)} and {\it (4)}.

If $\ell$ is not $\mathcal F$-invariant, Theorem \ref{TC:curvatura} implies that each irreducible component
of ${\rm tang}(\mathcal F,\mathcal L_i)$ is invariant by
$\mathcal L_i$ or  $\widehat {\mathcal L_{i}}$.
Since the leaves of $\mathcal L_i$ are lines through $q_i$ and because the algebraic curves invariant by
$\widehat{\mathcal L_i}$   must intersect $\ell$  on $\mathrm{sing}(\widehat{ \mathcal L_i } )\cap \ell = \{ q_i, \widehat{q_i} \}$ (according to Lemma  \ref{L:ric}),
it follows from
(\ref{E:l}) that
\[
|{\rm tang}(\mathcal F, \ell)  |\subset \bigcap_{i=1}^{k_{\ell}}
\{q_i, \widehat{q_i} \} \, .
\]
The only possibilities after an eventual  reindexing are
\begin{itemize}
\item[(a)] ${\rm tang}(\mathcal F, \ell) = \emptyset$  or
\item[(b)] ${\rm tang}(\mathcal F, \ell) = \{ q_1 \} $ and $\widehat{q_2} =
\cdots = \widehat{q_{k_{\ell}}} = q_1$.
\end{itemize}

We aim at a contradiction. On the one hand   (a) implies that
$\mathcal F$ is everywhere transversal to $\ell$. Therefore
$\mathcal F$ is of degree zero  what is not the case according to our hypothesis.
On the other
hand  (b) implies that the support of $\beta(q_1, \ldots,
q_{k_{\ell}})$ has a point with multiplicity at least $k_{\ell}-1$,
contradicting Lemma \ref{L:mult}.
\end{proof}

\begin{prop}\label{P:naocontrai}
Let $\mathcal F$ be a non-linear foliation on $\mathbb P^2$. Assume that
$\ell$ is a line that contains at least three points of {a}
set ${\mathcal P}$ of $k$
points in $\mathbb P^2$.  If $\mathcal F \boxtimes  {\mathcal W}({\mathcal P})$ has curvature
zero then the rational map $F: \mathbb P^2 \dashrightarrow \mathbb
P^1$ induced by the linear system $\{ \mathrm{tang}(\mathcal F,
\mathcal L_p) - \ell\}_{p\in \ell}$   does not contract $\ell$.
\end{prop}
\begin{proof}
First of all, a rephrasing of Proposition \ref{P:invariante}, yields that
 $\ell$ is a fixed component of {the pencil}
$\{ \mathrm{tang}(\mathcal F,
\mathcal L_p) \}_{p\in \ell}$. Thus $\{ \mathrm{tang}(\mathcal F,
\mathcal L_p) - \ell\}_{p\in \ell}$ is indeed a linear system.

Concretely, working with affine coordinates $(x,y)$  such that $\ell$
is the line at infinity  and $\mathcal F$ is induced by a polynomial
$1$-form $\omega= a(x,y)\, dx  + b(x,y) \,dy$ with isolated zeros then
$F(x:y:z)= (B(x,y,{z}):-A(x,y,{z}))$, where $A$ and $B$ are homogenizations
of $a$ and $b$ of degree $\max\{ \deg(a), \deg(b) \}=\deg(\mathcal F)$.

Assume that $F$
contracts $\ell$. It means that there  exists a point $p \in \ell$ such that $2\, \ell
\le \mathrm{tang}(\mathcal F,\mathcal L_p)$.  In other words the polynomials $A(x,y,0)$ and $B(x,y,0)$ are linearly dependent over $\mathbb C$.
Therefore
\begin{equation}\label{E:???}
|\mathrm{tang}(\mathcal F, \mathcal L_q) - \ell| \cap \ell = |\mathrm{tang}(\mathcal F, \mathcal L_{q'}) - \ell| \cap \ell 
{\varsubsetneq}
 \ell
\end{equation}
for every $q, q' \in \ell \setminus  \{ p \}$.

For any  $i\in\{ 1,\ldots,k_\ell\}$,
if $C$ denotes an irreducible component of $\mathrm{tang}(\mathcal F,\mathcal L_i)$ distinct from $\ell$, then Theorem \ref{TC:curvatura} implies  that $C$ necessarily is  ${\mathcal L}_i$-invariant or  $\widehat{{\mathcal L}_i}$-invariant.
Therefore, arguing as in the proof of Proposition \ref{P:invariante}, it follows from (\ref{E:???}) that
\[
|\mathrm{tang}(\mathcal F, \mathcal L_q) - \ell| \subset \bigcap_{{i=1}, { q_i \neq p }}^{k_{\ell}}\{ {q_i},{\widehat{q_i}} \} \,
\]
for every $q \in \ell \setminus  \{ p \}$, {in particular for every $q_i\in {\mathcal P}_\ell \setminus  \{ p \}$}.

After an eventual reindexing  we have  four
possibilities:\\
\begin{tabular}{p{2.5in}p{2.5in}}
\begin{enumerate}
\item[(a)] $p \notin \{q_1,\ldots, q_{k_{\ell}} \}$ \vspace{0.1cm}
\begin{enumerate}
\item[(a.1)] $\widehat{q_1}= \widehat{q_2} = \cdots = \widehat{q_{k_{\ell}}}$;\vspace{0.1cm}
\item[(a.2)] ${q_1}= \widehat{q_2} = \cdots = \widehat{q_{k_{\ell}}}$;
\end{enumerate}
\end{enumerate}
&
\begin{enumerate}
\item[(b)] $p = q_1 $\vspace{0.1cm}
\begin{enumerate}
\item[(b.1)] $\widehat{q_2}= \widehat{q_3} = \cdots = \widehat{q_{k_{\ell}}}$;\vspace{0.1cm}
\item[(b.2)] $q_2= \widehat{q_3} = \cdots = \widehat{q_{k_{\ell}}}$.
\end{enumerate}
\end{enumerate}
\end{tabular}

 Lemma \ref{L:mult} excludes the  cases (a.1), (a.2) and (b.1). To deal with the case (b.2) we will choose an identification
$\ell = \mathbb P ( \mathbb C x \oplus \mathbb Cy)$  where $q_1 = [
x ]$, $q_2 = [y ]$ and $q_i = [x + \lambda_i \,y]$ with $\lambda_i\neq 0$ for $i=3,\ldots,
k_{\ell}$.

A straight-forward computation shows that
\begin{equation*}
\begin{array}{lclcl} \widehat{q_1} &=&
\left(\displaystyle{\sum_{i=3}^{k_{\ell}} \frac{1}{\lambda_i}}
\right)x &+& (k_{\ell}-1)\,
y\, , \\
\widehat{q_2} &=& \big(k_{\ell}-1\big) \, x &+&
\left(\displaystyle{\sum_{i=3}^{k_{\ell}}
\lambda_i}\right)\,  y \, ,\\
 \widehat{q_j} &=& \left(-\displaystyle{\frac{1}{\lambda_j}} +
\displaystyle{\sum_{i=3,i\neq j}^{k_{\ell}}
\frac{1}{\lambda_i-\lambda _j}} \right) x &+& \left(1 +
\displaystyle{\sum_{i=3,i\neq j}^{k_{\ell}}
\frac{\lambda_i}{\lambda_i - \lambda_j}} \right) y \;
\end{array}
\end{equation*}
for $j$ ranging from $3$ to $k_\ell$.
 Now  $\widehat{q_j} = q_2$ for any such  $j$ implies that the
coefficient of $x$ in $\widehat{q_j}$ is zero. Summing up these coefficients  for~$j=3, \ldots
,k_{\ell}$,~one~obtains
\[
\sum_{j=3}^{ k_{\ell}} \frac{1}{\lambda_j} = 0 \, .
\]
Therefore $\widehat{q_1}=q_2$. Applying Lemma \ref{L:mult} once
again we conclude that (b.2) is also impossible. This concludes the
proof.
\end{proof}

\subsection{A bound for the degree of $\mathcal F$}

Combining Propositions \ref{P:invariante} and \ref{P:naocontrai}
with Riemman-Hurwitz formula we are able to bound the degree of
$\mathcal F$.

\begin{thm}\label{T:bd}
Let $\mathcal P \subset \mathbb P^2$ be a set  of $k \ge 4$
points and  $\mathcal F$ be a  a non-linear foliation on $\mathbb P^2$. If $K(\mathcal F \boxtimes \mathcal
W(\mathcal P))=0$
 then
$\deg(\mathcal F) \leq  4$. Moreover, if
$\deg(\mathcal F) \geq  2$,
and $\ell$ is a line containing
 $k_{\ell}$ points of $\mathcal P$
then $ k_{\ell}   \le 7 -\deg(\mathcal F) $.
\end{thm}
\begin{proof}
Assume that there is  no line that contains  at least three  points of $\mathcal P$. Then Corollary \ref{C:onaline} implies  that $\mathcal P$ has cardinality four and that $\mathcal F$ is the degree two foliation tangent to the
pencil of conics through $\mathcal P$.

From now on, we assume that there exists a line $\ell$
containing  $k_{\ell}$ points of $\mathcal P$, with $k_{\ell}  \ge
3$. Identifying $\ell$ with ${\mathbb P}^1$, let us note  $f:  \mathbb P^1\to \mathbb P^1 $  the restriction to
$\ell$ of the rational map $F: \mathbb P^2 \dashrightarrow \mathbb
P^1 $ induced by the linear system $\{ \mathrm{tang}(\mathcal F,
\mathcal L_p) - \ell\,|\, {p\in \ell} \}$.  Proposition
\ref{P:naocontrai} ensures that $f$ is a non-constant map.

The map $f$ is
characterized by the   following equalities  between divisors on $\ell$
$$f^{-1}(p) = \Big(\mathrm{tang}(\mathcal F, \mathcal L_p) -
\ell\Big)_{\vert \ell} \, ,$$
with  $p \in \ell$  arbitrary.

Let  $d$ be  the degree of ${\mathcal F}$. Recall from the proof of Proposition \ref{P:naocontrai} that $f$ is defined  by degree $d$ polynomials, that is,
${\rm deg}(f)=d$. Theorem \ref{TC:curvatura}  implies,  for any $i=1,\ldots,k_\ell$,
\begin{equation}\label{E:dois}
   f^{-1}(q_i) =  e_i\,q_i + (d-e_i) \,\widehat{q_i}
\end{equation}
 where  $e_i$ is an integer satisfying  $0\leq e_i\leq d$. Notice that the contribution of each of these fibers in
Riemann-Hurwitz formula is at least $d-2$. Therefore
\begin{equation*}
\label{E:rh}
\chi(\mathbb P^1)
= d\, \chi(\mathbb P^1) -(d-2)\, k_{\ell} - r
\end{equation*}
 for some non-negative
integer $r$. If  $d>2$ then
\[
{k_{\ell}} \le \frac{2d-2}{d-2} \, .
\]

If we keep in mind   that $k_{\ell} \ge 3$ and $d\ge 1$ then we end up  with
the following possibilities
\[
d=4\; \text{ and }\; k_{\ell} =3, \quad  \text{or}  \quad  d=3
\; \text{ and } \; k_{\ell} \le 4,  \quad \text{or}  \quad  1 \le
d\le2 \; \text{ and } \; k_{\ell} \ge 3 \, .
\]
If one realizes that for $d=2$ the map $f$ will have at most three
fixed points and two totally ramified points then one sees that in
this case $k_{\ell} \le 5$.
\end{proof}


The map $f:\mathbb P^1 \to \mathbb P^1$ used in the proof of Theorem
\ref{T:bd} codifies  a lot of information about the foliation $\mathcal F$. From now on we will refer to $f$ as the
{\it $\ell$-polar map of $\mathcal F$}.

\subsection{The polar map: properties and normal forms}
We use here  the same notations than in the  preceding section  and keep   the hypothesis of Theorem \ref{T:bd}.

We first state two  properties of the polar map that will be used in the sequel.

\begin{lemma}
\label{L:singpolarmap}
If the line  $\ell$ is  $\mathcal F$-invariant  then the  singularities of $\mathcal F$ on $\ell$ correspond to
the fixed points of $f$.
\end{lemma}
\begin{proof}
Let $(x,y) \in \mathbb C^2 \subset \mathbb P^2$ be  affine coordinates and assume that  $\ell$ is the line at infinity.
The foliation ${\mathcal F}$ is induced by a polynomial 1-form $\omega=a(x,y)\,dx+b(x,y)\,dy$ where $a(x,y)$ and $b(x,y)$ are relatively prime polynomials of degree $d$. If  $a_d(x,y)$ and $b_d(x,y)$ are the homogeneous components of degree $d$ of $a(x,y)$ and $b(x,y)$ (respectively) then, in the homogeneous coordinates $(x: y:0) \in \ell$, the polar map $f$ is
$$
f ( x:y )  =  \big[ b_d(x,y) : -a_d(x,y)\big] \; .
$$
On the other hand, one has
$$
{\rm sing}({\mathcal F})\cap \ell= \big\{ \, [x:y:0] \in  {\mathbb P}^2\, \big| \;x\,a_d(x,y)+y\,b_d(x,y)=0 \, \big\}\; .
$$
Thus $[x:y:0]\in \ell$ is a fixed point of $f$ if and only if it belongs to ${\rm sing}({\mathcal F})$.
\end{proof}

For $i=1,\ldots,k_\ell$, let  $e_i$ be the non-negative integer  appearing in  (\ref{E:dois}).
\begin{lemma}
\label{L:numberoflines}
There are exactly $e_i+1$  lines invariant by $\mathcal F$ through $q_i$   counted with the multiplicities that appear in
$\mathrm{tang}(\mathcal F, \mathcal L_{q_i})$.
\end{lemma}
\begin{proof}
Let $C$ be  an irreducible component  of ${\rm tang}({\mathcal F},{\mathcal L}_{q_i})$ passing through $q_i$.
According to Theorem \ref{TC:curvatura}, $C$ is ${{\mathcal L}_{q_i}}$-invariant or $\widehat{{\mathcal L}_{q_i}}$-invariant.
Since the only algebraic leaves of $\widehat{{\mathcal L}_{q_i}}$ trough $q_i$ are lines  (see Lemma \ref{L:ric} item {\it (6)})
$C$ must be a line. This fact together with  (\ref{E:dois})  proves   the lemma.
\end{proof}

It turns out that the relations (\ref{E:dois})  determine $f$ and $\mathcal
P_{\ell}$ up to automorphism of $\mathbb P^1$ when $\deg(\mathcal F) \ge 2$. Indeed by routine
elementary computations we arrive at the list presented in TABLE
\ref{T:polarf} below. For the sake of conciseness, we have chosen not to
derive this list here in full generality but just {to} deal with a particular case in the lemma below. All
the other cases follow from similar arguments.

\begin{lemma}
Assume that $k_\ell=3$ and ${\rm deg}({\mathcal F})=4$. Then we are in one of  the two cases {\rm (c.1)} or {\rm (c.2)} of TABLE \ref{T:polarf}.
\end{lemma}
\begin{proof}
In what follows,  $[a:b]$ designates the point $[a:b:0]$ on the line $\ell$ that is supposed to be  at infinity.
Since $k_\ell=3$, one can assume that  $q_1=[1:0]$, $q_2=[0:1]$ and $q_3=[1:-1]$; so $\widehat{q_1}=[ -1: 2 ]$, $\widehat{q_2}=[ 2:-1 ]$  and $\widehat{q_3}=[1 :1 ]$. By hypothesis ${\rm deg}({\mathcal F})=4$ so the polar map is $f(x:y)=(P(x,y):Q(x,y))$ where $P$ and $Q$ are homogenous polynomials in $x,y$, of degree 4. According to (\ref{E:dois}), one have $ (P)_0=e_2\,q_2+(4-e_2)\,\widehat{q_2} $
so $P$ is of the form $\lambda\,x^{e_2}\,(2\,y+x)^{4-e_2}$
for a certain $\lambda\in {\mathbb C}^*$ that can be supposed equal to 1. Similarly,
 $Q=\mu\,y^{e_1}\,(2\,x+y)^{4-e_1}$ with $\mu\in {\mathbb C}^*$. Since $(P+Q)_0=e_3\,q_3+(4-e_3)\,\widehat{q_3}$, there  exists
$\nu\in {\mathbb C}^*$ such that
$$
x^{e_2}\,(x+2\,y)^{4-e_2}+\mu\, y^{e_1}\,(y+2\,x)^{4-e_1}=\nu\, (x+y)^{e_3}(x-y)^{4-e_3}\, .  $$
After  straight-forward computations, it appears that such a relation in the space of homogeneous polynomials in two variables is only possible  when $(e_1,e_2,e_3)$ takes one of the two values:  $(1,1,1)$ or $(3,3,3)$. These correspond respectively to the cases {\rm (c.1)} and  {\rm (c.2)} in TABLE \ref{T:polarf} below.
\end{proof}

\begin{table}[h]
  \centering
  \begin{tabular}{|c|c|l|c|c|}
  \hline
  $\mathbf k_{\ell} $ & ${\mathbf d}$ &\qquad  {\bf action} & {\bf normal form for $f(x:y)$} & {\bf label} \\
  \hline \hline
 &  & $f^{-1}(q_i) = q_i + \widehat{q_i}$  &
   $\displaystyle{\big(x\,(2y+x):-\,y\,(2x+y)\big)}$ & (a.1) \\
  \cline{3-5}
     &  & $f^{-1}(q_1) = 2\, q_1 $  &  &\\

     &  & $f^{-1}(q_2) = 2 \, q_2 $  & $ \big(x^2 : -\,y^2 \big)$ &(a.2)\\

   $3$  & $2$ & $f^{-1}(q_3) = q_3 + \widehat{q_3} $  & &\\
\cline{3-5}
     &  & $f^{-1}(q_1) = 2 \,\widehat{q_1} $  &  & \\

     &  & $f^{-1}(q_2) = 2 \,\widehat{q_2} $  & $\big( (x+ 2y)^2 : - (2x +y)^2\,  \big)  $ & (a.3)\\

     &  & $f^{-1}(q_3) = q_3 + \widehat{q_3} $  & & \\
\cline{2-5}
   & $4$ & $f^{-1}(q_i) = q_i + 3\, \widehat{q_i}$ & $\displaystyle{ \Big(x\,(2y+x)^3 :-\,y\,(2x+y)^3\Big)}$ & (c.1)\\
\cline{3-5}
   &  & $f^{-1}(q_i) = 3\,q_i + \widehat{q_i}$ & $\displaystyle{
 \Big(x^3(2y+x):-\,y^3(2x+y)\Big)} $ & (c.2)\\
\hline \hline $4$ & $3$ & $f^{-1}(q_i) = q_i + 2\,\widehat{q_i}$ &
$\displaystyle{ \left(
3x \left(x + y\,(1-\xi_3^{{2}}) \right)^2
  :- y \big(3x + y\, (1-\xi_3^{{2}})  \big) ^2
 \right)} $    & (b.1)
\\ \hline \hline
 &  & $f^{-1}(q_1) =  2\,\widehat{q_1}$ &  &\\
 &  & $f^{-1}(q_2) =  2\,\widehat{q_2}$ &  & \\
$5$ & $2$  & $f^{-1}(q_3) = q_3+  \widehat{q_3}$ & $ \big(y^2 : -x^2 \big)$ & (a.4)\\
 &  & $f^{-1}(q_4) =  q_4+ \widehat{q_4}$ & &\\
 &  & $f^{-1}(q_5) =  q_5 + \widehat{q_5}$ & & \\ \hline
\end{tabular}
\vspace{0.1cm}
\caption{{\small The $\ell$-polar map of $\mathcal F$ assuming that $K\big(\mathcal F \boxtimes \mathcal W(\mathcal P)\big)=0
$.  The integer $k_{\ell}$ stands for the cardinality of $  \ell\cap \mathcal P$, $d$ designates
${\rm deg}(\mathcal F)$ and the points $q_1,q_2,q_3 \in \ell$  are normalized as
$q_1=[1:0:0]$, $q_2=[0:1:0]$ and $q_3=[1:-1:0]$.}
}\label{T:polarf}
\end{table}

\subsection{Points {of} $\mathcal P$ versus singularities of $\mathcal F$}

We start with a simple observation.

\begin{lemma}\label{L:passareta}
Let $\mathcal P$ be a collection of points of $\mathbb P^2$.  If $\mathcal F$ is a
non-linear foliation on $\mathbb P^2$ such that  $K(\mathcal F \boxtimes
\mathcal W(\mathcal P))=0$ then each point $p \in \mathcal P$ is contained in an $\mathcal F$-invariant line.
\end{lemma}
\begin{proof}
The argument used to settle  Lemma \ref{L:numberoflines} implies that every irreducible component
of $\mathrm{tang}(\mathcal F,\mathcal L_p)$ containing $p$ must be an $\mathcal F$-invariant line.
\end{proof}

TABLE \ref{T:polarf} allows us to restrain the possibilities of
$\mathcal F$  when $K(\mathcal F \boxtimes \mathcal W(\mathcal
P))=0$ and  $\deg(\mathcal F)> 1$. The next result   shows that once
$\mathcal F$ is known there are not many possibilities for $\mathcal
P$.

\begin{prop}\label{P:sing}
Let $\mathcal P$ be a finite set of points of $\mathbb P^2$. Suppose there exists a line $\ell$ containing at least  three points of $\mathcal P$. If $\mathcal F$ is a
non-linear foliation on $\mathbb P^2$ such that  the curvature of $\mathcal F \boxtimes
\mathcal W(\mathcal P)$ vanishes identically then $\mathcal
P \setminus \ell \subset \mathrm{sing}(\mathcal F)$.
\end{prop}
\begin{proof}
Let  $f:\mathbb
P^1 \rightarrow \mathbb P^1$ be the $\ell$-polar map of $\mathcal F$. Recall that  ${\mathcal P}_\ell={\mathcal P}\cap \ell=\{q_1,\ldots,q_{k_\ell}\}$ where the $q_i$'s are pairwise distinct.

For any distinct $i,j \in \{ 1, \ldots,
k_{\ell}\}$,
\begin{equation}
\label{E:singF}
\mathrm{sing}(\mathcal F) \cap (\mathbb P^2 \setminus \ell) =
{\big|}\mathrm{tang}(\mathcal F, \mathcal L_{q_i}){\big|}\cap
{\big|}\mathrm{tang}(\mathcal F, \mathcal L_{q_j}){\big|} \cap (\mathbb P^2
\setminus \ell)\, .
\end{equation}
Let $p$ be a point in $\mathcal P \setminus \ell$.
Assume that $p \notin \mathrm{sing}(\mathcal F)$. After an eventual reordering,  (\ref{E:singF}) implies
 that $p$ does not belong to $\mathrm{tang}(\mathcal F, \mathcal
L_{q_1})$ nor to $\mathrm{tang}(\mathcal F, \mathcal L_{q_2})$.

Since $p \notin \mathrm{tang}(\mathcal F, \mathcal L_{q_1})$, the
line $\overline{p q_1}$ is not $\mathcal F$-invariant. Thus
Proposition \ref{P:invariante} implies that   $\mathcal P \cap
\overline{p q_1} = \{ p, q_1 \}$. Consequently
$p \in \mathrm{sing}
(\widehat{\mathcal L_{q_1}})$
 thanks to  Lemma  \ref{L:ric} item~{\it (4)}.

Let $C$ be an irreducible component of $\mathrm{tang}(\mathcal F, \mathcal L_{q_1})$. If $C$ is not
$\mathcal L_{q_1}$-invariant then it must be $\widehat{\mathcal L_{q_1}}$-invariant by Theorem \ref{TC:curvatura} and
cannot contain $q_1$ by Lemma  \ref{L:ric}  item {\it (6)}. Thus $C$ must intersect
the $\widehat{\mathcal L_{q_1}}$-invariant line $\overline{pq_1}$ at $p$.
 Since  $p \notin \mathrm{tang}(\mathcal F, \mathcal L_{q_1})$ we deduce that
every irreducible component of  $\mathrm{tang}(\mathcal F, \mathcal L_{q_1})$ is $\mathcal
L_{q_1}$-invariant. Lemma \ref{L:numberoflines} implies that   $f^{-1}(q_1)= \deg(\mathcal F)\,  q_1$.
 Similarly, $\mathcal P \cap
\overline{p q_2} = \{ p, q_2 \}$ and $f^{-1}(q_2)= \deg(\mathcal F)\, q_2$.

Every rational self-map of $\mathbb P^1$ has at most two totally ramified
points (or at most two fixed points when the degree is one and the map is not the identity). Consequently  $p$ must belong to $ \mathrm{tang}(\mathcal F,
\mathcal L_{q_i})$ for every $i \in \{ 3, \ldots, k_{\ell}\}$. The
only possibility is that $k_{\ell}=3$ (otherwise $p$ would be in
$\mathrm{sing}(\mathcal F)$ according {to} (\ref{E:singF})).

Lemma \ref{L:passareta} implies that there is a $\mathcal F$-invariant line $\ell_p$ through $p$. Since $\overline{p q_1}$ and
$\overline{p q_2}$ are not $\mathcal F$-invariant, the line $\ell_p$ must be distinct from these.
In particular, $\ell_p \cap \ell$ must be contained in $(\mathrm{sing}(\mathcal F) \cap \ell ) \setminus \{ q_1, q_2 \}$.
Therefore $\mathrm{sing}(\mathcal F) \cap \ell$ has cardinality at least three and
consequently, the degree of $\mathcal F$ is at least two.
 After  analyzing TABLE
\ref{T:polarf}, one concludes that the map $f$ must be as  in case (a.2). Explicitly  $f^{-1}(q_1)=2\,q_1$, $f^{-1}(q_2)=2\, q_2$
and $f^{-1}(q_3)=q_3 + \widehat{q_3}$.
 In particular $\mathcal F$ has degree two;  admits exactly three singularities on $\ell$, namely $q_1,q_2$ and $q_3$; and $\ell_p=\overline{pq_3}$ is the unique $\mathcal F$-invariant line through $p$.

Notice that  $\mathrm{tang}(\mathcal F, \mathcal L_p)$  is an effective divisor of degree three and its support
contains both  $\ell_p$ and the singular points of $\mathcal F$. Since
$q_1, q_2$ and $p$ are not aligned, there exists an irreducible component $C$ of  $\mathrm{tang}(\mathcal F, \mathcal L_p)$
distinct from $\ell_p$ and with degree at most two. According to Theorem \ref{TC:curvatura}  $C$  is invariant by $\widehat{\mathcal L_p}$ or ${\mathcal L_p}$.

If $C$ contains  $p$ then by Lemma \ref{L:ric} item {\it (6)} it must be a line and therefore is equal to $\ell_p$. This is not possible
due to our choice of $C$. Thus $C$ does not contains $p$ and must
be  $\widehat{\mathcal L_{p}}$-invariant.

Recall from above that  $\overline{pq_i}\cap \mathcal P = \{p,q_i\}$ for $i=1,2$. Corollary \ref{C:riccati} implies
that the irreducible curves invariant by $\widehat{\mathcal L_p}$ that are not lines must have degree at least three. Thus
we can assume that  $C$ is a line. Moreover
$$ \qquad {\rm sing}\big( \widehat{\mathcal L_p}\big) \cap \overline{pq_i}=\{p,q_i\} \;\, \text{ for } \, \, i=1,2 $$
thanks to Lemma  \ref{S:riccati} item {\it{(4)}}. Because the intersections of $C$ with   $\overline{pq_1}$ and $\overline{pq_2}$ are
 singularities of $\widehat {\mathcal L_p}$ that are distinct from $p$, we conclude that  $C=\ell$. However  $\ell$ is $\mathcal F$-invariant but
 not $\mathcal L_p$-invariant and consequently cannot be in $\mathrm{tang}(\mathcal F, \mathcal L_p)$.  Thus the
 assumption $p\notin \mathrm{sing}(\mathcal F)$ leads to a contradiction. The  proposition follows. \end{proof}


\section{Exceptional CDQL webs of degree one on $\mathbb P^2$}\label{S:deg1}

The degree of a web $\mathcal W$ on $\mathbb P^2$ is, like in the case of foliations, the number of tangencies
of $\mathcal W$ with a generic line. In particular the degree of a completely decomposable web  is nothing more than the sum
of the degrees of its defining foliations and  the degree of an CDQL web is nothing more than
the degree of its non-linear defining  foliation.

\subsection{Infinitesimal automorphisms}

\begin{prop}
\label{P:cdqldeg1}
Let $\mathcal W=\mathcal F \boxtimes \mathcal W(\mathcal P)$ be a
CDQL  $(k+1)$-webs  of degree one with $k\ge 4$. If $K(\mathcal
W)=0$ then it exists a line $\ell$ containing at least $k-1$ points of $\mathcal P$.
Moreover there is a system of affine coordinates $(x,y) \in \mathbb C^2 \subset \mathbb P^2$ where
$\ell$ is the line at infinity,  $\mathcal F$ is induced by a homogeneous $1$-form $\omega_0$ of degree 1, and
 the radial vector field $R= x \partial_x + y \partial_y$ is an infinitesimal automorphism of
$\mathcal W$.
\end{prop}
\begin{proof}
 If $K(\mathcal
W)=0$ then   Corollary \ref{C:onaline} and
Proposition \ref{P:invariante} imply
that there is a ${\mathcal F}$-invariant line $\ell$ that contains (at least) three points of ${\mathcal P}$.
A classical result by Darboux says that a degree $d(=1)$ foliation on $\mathbb P^2$ has $d^2+d+1(=3)$ singularities
counted with multiplicities.
Since at least two of the three singularities of ${\mathcal F}$ necessarily lie  on $\ell$, it follows that ${\rm sing}({\mathcal F})\setminus \ell$ reduces to a point or is empty.
Proposition \ref{P:sing} yields that
at least  $k-1$ points of
$\mathcal P$ lie on $\ell$. According to Proposition \ref{P:naocontrai}  the $\ell$-polar map of ${\mathcal F}$ does not
contract $\ell$ so one of the singularities
of $\mathcal F$ is not contained in $\ell$.

All that said we can choose affine coordinates where $\ell$ is the
line at infinity and $\mathcal F$ is induced by a homogeneous linear
$1$-form $\omega_0$ that vanishes only at the origin of ${\mathbb  C}^2$.
It is {then} clear  that $R$ is a infinitesimal automorphism of $\mathcal W$.
\end{proof}

It is a simple exercise to show  that after a linear change of coordinates the $1$-form  $\omega_0$ that defines
$\mathcal F$ on the system of affine coordinates given by  Proposition \ref{P:cdqldeg1} can be written as
\begin{equation*}{}^{} \hspace{1cm}
 \omega_0'=y\,dx-(x-y)\,dy \qquad \mbox{ or } \qquad
\omega_{0}^\kappa=y\,dx-\kappa \,x\,dy \quad \mbox{with } \;\kappa\neq 0,1\; .
\end{equation*}
The canonical first integral $u_0$ of ${\mathcal F}$ (see Section \ref{S:reviewMPP}) with respect to the radial vector field $R$ is then
\begin{equation}
\label{E:FInormalform}
u_0'={x}/{y}+\log(y)
\qquad \mbox{ or } \qquad
u_0^\kappa= \frac{1}{1-\kappa} \log\left(\frac{x}{y^\kappa}\right)\; .
\end{equation}

Denote by $\mathcal F_R$  the
foliation induced by the radial vector field.
If $\mathcal P$ is not included in $\ell$ then $\mathcal W(\mathcal P) = \mathcal F_R
\boxtimes \mathcal W( \mathcal P \cap \ell)$. It follows from
Theorem \ref{T:1} that $\mathcal F \boxtimes \mathcal W(\mathcal P)$
has maximal rank if and only if $\mathcal F \boxtimes \mathcal
W(\mathcal P \cap \ell)$ does. Therefore we will restrict our
attention to the case where $\mathcal P \subset \ell$. Until we say
otherwise, $\mathcal W(\mathcal P)$ is a web
induced by    $k\geq 3$ constant $1$-forms $\omega_1, \ldots,
\omega_k$ and $\mathcal W = \mathcal F \boxtimes \mathcal W(\mathcal P)$ is a
CDQL web of degree one admitting $R$ as an infinitesimal automorphism.

\subsection{The action on $\mathcal A(\mathcal W)$ is semisimple}
Recall from Section \ref{S:reviewMPP} that $L_R$ acts on the space of abelian relations of $\mathcal W$.
In our particular setup this action is semisimple as the next proposition shows.

\begin{prop}\label{P:ddd}
The linear map $L_R: \mathcal A(\mathcal W) \to \mathcal A(\mathcal
W)$ is diagonalizable and its eigenvalues are non-negative integers.
Moreover the zero eigenspace of $L_R$, if not trivial,  has
dimension one.
\end{prop}
\begin{proof}
For $i=1, \ldots, k$, the canonical first integrals for the foliations $[\omega_i]$ ({with respect to the radial vector field}) are the  functions
$$
  u_i = \int \frac{\omega_i}{\omega_i(R)} \, .
$$
Clearly  $u_i=\log h_i$ {for  suitable} linear forms $h_i  \in \mathbb
C[x,y]$.

On  a simply connected open set contained  in the complement
of $\Delta(\mathcal W)$, let us consider an abelian relation of the form
\begin{equation}\label{E:m}
\sum_{j=0}^k  \\
P_j (u_j)\, e^{\lambda \,u_j} du_i =0\;
\end{equation}
corresponding to an eigenvalue $\lambda$ of $L_R$  (see Proposition \ref{P:description}).

 Analytic continuation of (\ref{E:m}) {along} closed paths
homotopic to zero in $\mathbb C^2 \setminus \{ h_2 \cdots h_k = 0
\}$ but not homotopic to zero in $\mathbb C^2 \setminus \{ h_1 = 0
\}$ implies that
\[
P_1 (z +  2\pi i  ) \,e^{2\pi i \lambda }= P_1(z)\,
\]
for every $z\in \mathbb C$.
Therefore $P_1$ must be constant and $\lambda$ an integer. In the same way, one proves that $P_2,\ldots,P_k$ are constant polynomials.
Thus $P_0 (u_0)\, e^{\lambda \,u_0} $ must be a rational function. Taking into account  (\ref{E:FInormalform}),  one deduces that  $P_0$  is also constant. According to the last part of  Proposition \ref{P:description}, the linear operator $L_R$ is diagonalizable.

Assume that $\lambda\leq 0$.  Equation (\ref{E:m}) takes the form
\begin{equation}\label{E:m2}
c_0\,e^{\lambda \,u_0} du_0+\sum_{j=1}^k
c_j\,{h_j}^{\lambda-1}dh_j =0
\end{equation}
for certain constants $c_0,c_1,\ldots,c_k\in \mathbb C$. It follows that the line $h_j=0$ is invariant by ${\mathcal F}$ as soon as  $c_j\neq 0$.
But ${\mathcal F}$ admits at most two invariant lines through the origin.
Thus if (\ref{E:m2}) is not trivial, one can assume that it has the form
$$e^{\lambda \,u_0} du_0+c_1\,{x}^{\lambda-1}dx+c_2\,y^{\lambda-1}dy =0\, .$$
Since the curvature of $[dx\,dy]\boxtimes[ydx - (x-y)dy]$ is non-zero, an identity of the
form (\ref{E:m}) holds only  when $\omega_0=\omega_{0}^\kappa$ for a certain $\kappa$ and the eigenvalue $\lambda$ is zero.
\end{proof}

\begin{cor}
If ${rk}(\mathcal F \boxtimes \mathcal W(\mathcal P)) >{rk}(\mathcal W(\mathcal P))+1 $ then $\mathcal F$ admits a
polynomial first integral of the form $x^p y^q$ where $p $ and $q$ are
relatively prime positive integers.
\end{cor}
\begin{proof}
It follows from Proposition \ref{P:ddd} that
there is at least
one  abelian relation of the form (\ref{E:m2}) with $c_0=1$ and  $\lambda\in {\mathbb N}^*$. Integrating these, it follows that
\[
 \int e^{\lambda u_0} du_0  =  \frac{1 }{\lambda}  \sum_{j=1}^k {c_j}\,
 h_j^{\lambda }
\]
is a homogeneous polynomial first integral of $\mathcal F$. Then $\omega_0=\omega_{0}^{p/q}$ where $p $ and $q$ are
relatively prime positive integers.
\end{proof}

\subsection{Characterization of $\mathcal F$}

Let $\delta\in\{0,1,2\}$ be the number of  lines  through
the origin of
$ \mathbb C^2$ that are invariant by $\mathcal F$.
\begin{lemma}
\label{L:rkF*W(P)}
If $\mathcal F$ has a first integral of the form $x^p y^q$ where
$p,q \in \mathbb N$ are relatively prime then
\[
  {rk} \big(\mathcal F \boxtimes \mathcal W(\mathcal P)\big) - {rk} \big(\mathcal W(\mathcal P)\big)
   \le \frac{2k-2}{p+q}\,  .
\]
Moreover if the equality holds then  $\delta \neq 1$.
\end{lemma}
\begin{proof}
The abelian relations involving $\mathcal F$ and corresponding to
strict positive eigenvalues can be written in integrated form as
\[ \qquad \qquad
(x^py^q)^r = \sum_{i=1}^k \mu_i\,  h_i^{(p+q)r}
\qquad (\mu_1,\ldots,\mu_k \in {\mathbb C})\; .
\]
Let $\alpha \in \mathbb C^*$ be sufficiently general and set
$\varphi(x,y)=(x,\alpha \,y)$.   Then
\[
0=(x^py^q)^r - \alpha^{-qr}{\varphi^*\big((x^py^q)^r}\big)= \sum_{i=1}^k
\mu_i \, h_i^{(p+q)r} - {\alpha^{-qr}}\sum_{i=1}^k \mu_i\,
\varphi^*\big(h_i^{(p+q)r}\big)\; .
\]
where the right-hand side involves at most $2k - \delta$ distinct
linear forms.

It is very convenient to  interpret geometrically this equality in terms of the rational normal curve
$\Gamma$ in $\mathbb P^{(p+q)r}$ of degree $(p+q)r$. It says that there are $2k-\delta$
distinct points on $\Gamma$ that are not in general position.  It is classical result that  $m$ distinct points on the rational normal curve of degree $l$ are in general position if $m\leq l+1$. Therefore   $(p+q)r + 2 \le 2k - \delta$. Hence
\[ r \le \frac{2k-\delta - 2}{p+q}
\, .
\]

Recall from the proof  of Proposition \ref{P:ddd} that when $\delta\neq
2$, all the abelian relations are polynomials identities  and when $\delta=2$
there is exactly one extra logarithmic abelian relation. It follows
that
\[
\dim \frac{\mathcal A\big(\mathcal F \boxtimes \mathcal W(\mathcal
P)\big)}{\mathcal A\big(\mathcal W(\mathcal P)\big)} \le \left\{
\begin{array}{lcl}
                             \displaystyle{\frac{2k-\delta - 2}{p+q}} & \text{ when} & \delta \in \{ 0,1\}, \vspace{0.15cm} \\
                              \displaystyle{\frac{2k- 2}{p+q}}  & \text{ when} &
                             \delta = 2.
                           \end{array}
  \right.
\]
\end{proof}

\subsection{The classification}

If a web $\mathcal F \boxtimes \mathcal W(\mathcal P)$ {of degree $1$} is of maximal
rank it follows from  Lemma \ref{L:rkF*W(P)} that $\delta$ must be $0$ or
$2$ and $p=q=1$. It turns out that there exist examples of $(k+1)$-webs of maximal rank and of degree 1 for
any  $\delta \in \{ 0,2\}$ and any $k\ge 3$, namely the webs ${\mathcal A}_{I}^{k}$ and ${\mathcal A}_{III}^{k-2}$ of the Introduction.

To complete the classification of exceptional CDQL webs of degree one on ${\mathbb P}^2$, it suffices to show that these examples are the only ones up to
projective automorphisms.

In what follows, ${\mathcal P}$ and ${\mathcal Q}$ designate two sets of $k$ points on the line at infinity, disjoint of $[1:0:0]$ and $[0:1:0]$, with $k\geq 3$ if $\delta=0$ and $k \geq 1$ if $\delta=2$.
Let   $a_1,\ldots,a_k, b_1,\ldots,b_k\in {\mathbb C}^*$ be such that
$$
{\mathcal W}({\mathcal P})=\big[ d(x+a_1\,y)\cdots \, d(x+a_k\,y)\big] \quad\mbox{ and } \quad
{\mathcal W}({\mathcal Q})=\big[ d(x+b_1\,y)\cdots  d(x+b_k\,y)\big]\; .
$$
In particular,  $a_i\neq a_j$ and $b_i\neq b_j$ for all $i,j$ such that $1\leq i<j\leq k$.

\begin{prop}  If the two   CDQL $(k+1)$-webs
 $[d(xy)]\boxtimes {\mathcal W}({\mathcal P}) $ and $[d(xy)]\boxtimes {\mathcal W}({\mathcal Q}) $ are both of maximal rank, then they are
projectively equivalent.
\end{prop}
\begin{proof}
From the proof of Lemma \ref{L:rkF*W(P)}, it follows that
\begin{eqnarray*}
 \big(xy\big)^{k-1} &=& \sum_{i=1}^k  \lambda_i \,(x+a_i\,y)^{2\,k-2} =  \\
  &=& \sum_{i=1}^k  \mu_i \,(x+b_i\,y)^{2\,k-2}
\end{eqnarray*}
for suitable complex numbers $\lambda_i, \mu_i$.
Since, for any $\nu\in {\mathbb C}^*$, the automorphisms $(x,y)\mapsto (x,\nu\,y)$  preserve the foliation $[d(xy)]$,
one can assume that $a_1=b_1$ with no loss of generality. Subtracting the two summations yields
\begin{equation}
\label{E:toto}
0= (\lambda_1 - \mu_1 )\,(x+a_1\,y)^{2\,k-2} +
\sum_{i=2}^k \lambda_i \,(x+b_i\,y)^{2\,k-2}
- \sum_{i=2}^k \mu_i \,(x+b_i\,y)^{2\,k-2}\; .
\end{equation}
One can interpret    this relation geometrically as in the proof of Lemma \ref{L:rkF*W(P)} by considering the powers $(x+a_i\,y)^{2\,k-2}$ and
$(x+b_i\,y)^{2k-2}$ (for $i=1,\ldots,k)$ as points on the rational normal curve of degree $2k-2$. Notice that the number of these points is at most  $2k-1$. Since  $m$ distinct points on the rational normal curve of degree $2k-2$ are inevitably in general position when $m\leq 2k-1$, {the relation} (\ref{E:toto}) implies that the sets $\{a_1,\ldots,a_k\}$ and $\{b_1,\ldots,b_k\}$ coincide. The proposition follows.
\end{proof}

In the same way, one   proves the
\begin{prop}  If the two   CDQL webs
 $[d(xy)\,dx\,dy]\boxtimes {\mathcal W}({\mathcal P}) $ and $[d(xy)\,dx\,dy]\boxtimes {\mathcal W}({\mathcal Q}) $ are both of maximal rank, then they are
projectively equivalent.
\end{prop}

At this point we have concluded the classification of exceptional CDQL webs of degree one on $\mathbb P^2$. They are all projectively equivalent
to one the webs in the families $\mathcal A^*_k$. We point out that we have made a  heavy use of the structure of the space of abelian relations
of these webs. It would be interesting to find an alternative approach more focused on the curvature. For instance one could try
to classify all the flat CDQL webs of degree 1 on $\mathbb P^2$.


\section{Flat CDQL webs on $\mathbb P^2$ of degree at least two}\label{S:flat}

 Based on the results of Section  \ref{S:constraints} we will derive a complete list of flat CDQL
$(k+1)$-webs on $\mathbb P^2$ of degree at least two  when $k \ge
4$. Up to automorphisms of $\mathbb P^2$, there
are exactly sixteen examples  --- nine   of degree 2, three of degree 3 and
four of degree 4.

\subsection{Flat CDQL webs of degree two}
In the present and   in the next two subsections,
we will treat independently
the three possibilities for the degree of ${\mathcal F}$. We start by considering  flat CDQL webs
of degree two.

\begin{prop}\label{P:deg2}
Let $\mathcal F$ be a  foliation of degree 2 and $\mathcal P
\subset \mathbb P^2$ be a finite set of at least four points. If
$K(\mathcal F \boxtimes \mathcal W(\mathcal P))=0$ then
$\mathcal F$ is projectively equivalent to one of the following foliations:
{\begin{enumerate}
\item[(a.1.h)] $\mathcal F = \big[d \left( xy(x+y) \right)\big]\,; $\vspace{0.15cm}
\item[(a.2.h)] $\mathcal F= \left[d \big( \frac{xy}{x+y} \big) \right];$\vspace{0.15cm}
\item[(a.3.h)] $\mathcal F= \big[d \left( (4y^2 + xy + 4x^2)^3(x+y) \right) \big]; $\vspace{0.15cm}
\item[(a.4.h)] $\mathcal F
= {\big[d \left(x^3+y^3 \right) \big]} ;  $\vspace{0.15cm}
\item[(a.2)] $\mathcal F= \left[d \left( \frac{y^2-1}{x^2-1} \right) \right]
.$\vspace{0.15cm}
\end{enumerate}}
\noindent Moreover  in the cases  (a.1.h), (a.2.h) and (a.3.h),
$\mathcal P$ has cardinality four    and is  equal to
the singular set of $\mathcal F$.  In the case (a.4.h) there are two
possibilities for $\mathcal P$. Either $\mathcal P$ is equal to
 $\mathrm{sing}(\mathcal F) \cup \{ [0:1:0],[1:0:0]
 \}$ or to  $  \left( \mathrm{sing}(\mathcal F) \cup \{ [0:1:0],[1:0:0] \}
 \right) \setminus\{ [0:0:1] \}
 .$ Finally, in case (a.2) the set $\mathcal P$ is any of
the subsets of $\mathrm{sing}(\mathcal F)$ containing the four base
points of the pencil $< x^2 -z^2 , y^2 -z^2>$. Up to the
automorphism group of $\mathcal F$ there are only four
possibilities.
\end{prop}
\begin{proof}

If the points in  $\mathcal P$  are in general position then, according to Corollary \ref{C:onaline}, $\mathcal F$ is the pencil generated by
two reduced  conics intersecting transversally   and $\mathcal P$ is the set of base
points of this pencil. So $\mathcal F \boxtimes \mathcal W(\mathcal P)$ is Bol's web and we are in case {\it (a.2)}.

From now on we will assume that there is a line $\ell \subset
\mathbb P^2$ that contains at least three points of $\mathcal P$. Up to
the end of the proof  we  work with    affine
coordinates $[x:y:1]$ on $\mathbb C^2 \subset \mathbb P^2$, for which  $\ell=\{z=0\}$ is the line at infinity.
We will also assume  that  $\mathcal P\cap \ell$ contains  $q_1 = [1:0:0]$,
$q_2=[0:1:0]$ and $q_3=[1:-1:0]$.

We will deal separately with each one of the four possibilities given by TABLE \ref{T:polarf} for the
$\ell$-polar map {$f$} of $\mathcal F$.\vspace{0.15cm}

\noindent{{\bf Case (a.1). }} In this case $k_{\ell}=3$ and  $f^{-1}(q_i) = q_i + \widehat{q_i}$  for  $i= 1, 2, 3$.
\medskip
Notice that $f^{-1}(q_1) = q_1 +
\widehat{q}_1$ implies that $\mathrm{tang}(\mathcal F,\mathcal
L_{q_1})$ is the union of three lines: the line at infinity $\ell$
together with two other lines, one intersecting $\ell$ at $q_1$ and the
other at  $\widehat{q}_1$. A similar situation occurs for
$\mathrm{tang}(\mathcal F,\mathcal L_{q_2})$ and $\mathrm{tang}(\mathcal F,\mathcal L_{q_3})$.

Therefore $f^{-1}(q_1)=q_1 + \widehat{q_1}$ and $f^{-1}(q_2)=q_2 +
\widehat{q_2}$  imply that  $\mathcal F$ is induced  by a $1$-form like
\[
\omega = (y + c_1) (2x+y+c2)\, dx + (x + c_3) (2y+x + c_4) \,dy \, ,
\]
where $c_1, c_2, c_3$ and $ c_4$ are complex constants. After composing with a
translation we can assume that $c_1=c_3=0$.

It remains to consider the conditions imposed by  $f^{-1}(q_3)=q_3 +
\widehat{q_3}$. Notice that
$\mathrm{tang}(\mathcal F, \mathcal L_{q_3})$ is cut out by
\[
 y  (2x+y+c2) - x  (2y+x + c_4) = y^2 - x^2 +  c_2 y - c_4 x \, .
\]
This latter expression is a product of lines if and only if
$c_2= \pm c_4$. When $c_2=c_4=0$ we arrive at the homogeneous
foliation $$\mathcal F = \left[d \big( xy(x+y) \big)\right]. $$
We are in case {\it (a.1.h)}. Because the cardinality of the singular set of
$\mathcal F$ is four, {there  is} just one possible choice for $\mathcal P$: $\mathcal P= \mathrm{sing}(\mathcal F)$.

 If $c_2 \neq 0$ then after applying a homothety we can assume that
$c_2=1$. We arrive at two possibilities for $\omega$, namely
\[
\omega_{\pm} =  y\,  (2x+y+1)\,dx +  x \, (2y+x \pm 1)\, dy \, .
\]
Let $\mathcal F_{\pm}$ are the corresponding foliations. By hypothesis, $k_\ell=3$ and ${\mathcal P}_\ell=\{q_1,q_2,q_3\}$.
If $\mathcal F_{\pm}\boxtimes {\mathcal W}({\mathcal P})$ is assumed to be flat then  Proposition \ref{P:sing} implies that ${\mathcal P}\setminus \ell$ is included in the support of ${\rm sing}({\mathcal F}_{\pm})\cap \mathbb C^2$. In particular  there are  only a finite number of possibilities for ${\mathcal P}$.
Lengthy, but straight-forward, computations
shows that $K(\mathcal F_{\pm} \boxtimes \mathcal W(\mathcal
Q\cup \{q_1,q_2,q_3\}))\neq 0$ for any non-empty subset ${\mathcal Q}\subset
{\rm sing}({\mathcal F}_{\pm})\cap \mathbb C^2$.
Therefore the  foliations  $\mathcal F = \mathcal F_{\pm}$ are not among the defining foliations of  any flat CDQL webs of order at least five.

\vspace{0.15cm}

\noindent{{\bf Case (a.2).} } In this case $k_{\ell}=3$,   $f^{-1}(q_i) = 2q_i$ for $i=1,2$ and $f^{-1}(q_3) = q_3 + \widehat{q_3}$.
\medskip
Arguing as in  the paragraph above,
we conclude that $\mathcal F$ is induced
by
\begin{equation}
\label{E:caseA2omega}
\omega = y\,(y-1)\, dx + x\,(x-1) \,dy  \qquad \text{ or } \qquad \omega' = y^2
 dx + x^2  dy \, .
\end{equation}

Recall  that ${\mathcal P}\setminus \ell$ is included in  ${{\rm sing}({\mathcal F})}\cap \mathbb C^2$ (according to Proposition \ref{P:sing}).
If ${\mathcal F}$ is induced by $\omega'$, only one  possibility can happen, namely ${\mathcal P}=\{q_1,q_2,q_3,p_4\}$ where $p_4=[0:0:1]$ (since ${\rm sing}({\mathcal F})= \{ q_1,q_2,q_3,p_4\} $). By a direct computation, one verifies that the 5-web defined by ${\mathcal P}$ and $\omega'$ is indeed flat.

Let us now consider the case when ${\mathcal F}$ is the foliation induced by $\omega$. If we set $p_5=[1:1:1:]$, $p_6=[0:1:1]$ and
$p_7=[1:0:1]$ then
$$\mathrm{sing}(\mathcal F)=\{ q_1,q_2,q_3,p_4, p_5,p_6, p_7 \}
\; . $$

A direct computation shows that there are exactly four  subsets $\mathcal P$ of $\mathrm{sing}(\mathcal F)$ that strictly contain${ \mathcal P_\ell}= \{q_1,q_2,q_3\}$ and that verify $K(\mathcal F
\boxtimes \mathcal W(\mathcal P))=0$, namely
\[
\mathcal P = \left\lbrace \begin{array}{lcl}
                             { \mathcal P_\ell} & \cup & \{ p_4,p_5\}\, ,  \\
                            { \mathcal P_\ell} & \cup & \{ p_4,p_5,p_6\}\, , \\
                            { \mathcal P_\ell}& \cup & \{ p_4,p_5,p_7\}\, , \\
                            { \mathcal P_\ell}& \cup & \{ p_4,p_5,p_6,p_7\}\; .
                           \end{array}
\right.
\]
 Notice that $\mathcal F \boxtimes \mathcal
W(\{q_1,q_2,p_4,p_5\})$ is nothing more than Bol's exceptional
$5$-web. The second and the third possibilities for ${\mathcal P}$ are equivalent since they are  interchanged by the $\mathcal
F\boxtimes \mathcal W(  \mathcal P_\ell)$-automorphism
$(x,y)\mapsto(y,x)$. All the cases above lead to exceptional webs. Indeed
they are the webs $\mathcal B_6, \mathcal B_7$ and $\mathcal B_8$ of the Introduction that  have been previously studied in
\cite{Robert,PTese}. \vspace{0.15cm}

\noindent{{\bf Case (a.3).} }  Here $k_\ell =3$, $f^{-1}(q_i) = 2\,\widehat{q_i}$ for $i=1,2$ and  $f^{-1}(q_3)= q_3 +  \,\widehat{q_3}$.
Theorem \ref{TC:curvatura} tell us that every irreducible component $C$ of $\mathrm{tang}(\mathcal F, \mathcal L_{q_1})$ is
 invariant
by $\mathcal L_{q_1}$ or $\widehat{\mathcal L_{q_1} }$. Because $f^{-1}(q_1)=2 \widehat{q_1}$, there exists such $C$ invariant by $\widehat{\mathcal L_{q_1}}$ and
distinct from $\ell$.
The divisor  $\mathrm{tang}(\mathcal F, \mathcal L_{q_1}) -\ell$ is effective and of degree $2$. Consequently the degree of $C$ is at most two.
If it is two then Corollary \ref{C:riccati} implies that for every point $p \in \mathcal P \setminus \{ q_1, q_2, q_3 \}$ the line $\overline{q_1p}$ contains a third
point of $\mathcal P$.  Proposition \ref{P:invariante} implies the $\mathcal F$ invariance of $\overline{q_1p}$ contradicting $f^{-1}(q_1)=2 \widehat{q_1}$. This proves
that for $i=1,2$ every irreducible component of  $\mathrm{tang}(\mathcal F, \mathcal L_{q_i}) -\ell$ must be a $\widehat{\mathcal L_{q_i}}$-invariant line through $\widehat{q_i}$. Because $\mathrm{Card}(\mathcal P)\ge 4$ and $\mathcal P \nsubseteq \ell$, for $i=1,2$, the foliation $\widehat{\mathcal L_{q_i}}$ has only one invariant
line through $\widehat{q_i}$ distinct from $\ell$. Therefore there exists constants $c_1$ and $c_2$  for which  $\mathcal F = [(2x+y+c_1)^2dx+(x+2y+c_2)^2dy]$. Modulo a translation, we can assume that $c_1=c_2=0$. Thus $${\mathcal F}=\big[ (2x+y)^2dx+(x+2y)^2dy\big]=\big[d \left( (x+y)(4y^2 + xy + 4x^2)^3 \right) \big]\,. $$
We are in case {\it (a.3.h)} and {necessarily $\mathcal P= \mathrm{sing}(\mathcal F)$ since ${\rm Card}(\mathrm{sing}(\mathcal F))=4$}.

\vspace{0.15cm}

\noindent{{\bf Case (a.4).}}  We finally arrive at the last case of TABLE \ref{T:polarf} where  $k_\ell =5$, $f^{-1}(q_i) = 2\,\widehat{q_i}$ for $i=1,2$ and  $f^{-1}(q_j)= q_j +  \,\widehat{q_j}$ for $j=3,4,5$.
\medskip

\noindent Arguing exactly as in case (a.3) one can show that   in this case   $\mathcal F$  is also  homogeneous.
Therefore  ${\mathcal F}=[x^2dx+y^2dy]=[d(x^3+y^3)]$ and, as it was
 shown in Section \ref{S:action},  any of the two possibilities for $\mathcal P$, namely
$$\mathcal P = \{ q_1, \ldots, q_5 , [0:0:1] \} =  \mathrm{sing}(\mathcal F) \cup \{ [0:1:0],[1:0:0]
 \} $$
or
$$\mathcal P=   \{ q_1, \ldots, q_5  \} =\left( \mathrm{sing}(\mathcal F) \cup \{ [0:1:0],[1:0:0] \}
 \right) \setminus\{ [0:0:1] \}
 $$
 leads to exceptional, and in particular flat,  webs.
\end{proof}

\subsection{Flat CDQL webs of degree three}
The classification of flat CDQL webs
of degree three is given by the following proposition.

\begin{prop}\label{P:deg3}
Let $\mathcal F$ be a  foliation of degree three  and $\mathcal P
\subset \mathbb P^2$ be a finite set of at least four points. If
$K(\mathcal F \boxtimes \mathcal W(\mathcal P))=0$ then
$\mathcal F$ is projectively equivalent to one of the following foliations:
\begin{enumerate}
\item[(a)] $\mathcal F = \big[
\left( {x}^{3} +{y}^{3}+1 +6\,x{y}^{2}\right) dx -  \left( x^3 +
{y}^{3}+1 +6\,{x}^{2}y \right) dy \big]\,;  $ \vspace{0.15cm}
\item[(b)] $\mathcal F= \big[ d \big(x(x^3 +y^3)\big)\big]$.
\end{enumerate}
Moreover  $\mathcal P = \mathrm{sing}(\mathcal F) \cap \{
x-y=0\}$ in case {\it (a)} and  $\mathcal P = \mathrm{sing}(\mathcal F)$ or
$\mathcal P = \mathrm{sing}(\mathcal F) \setminus \{ [0:0:1] \}$ in case {\it (b)}.
\end{prop}
\begin{proof}
Corollary \ref{C:onaline} implies that there exists a line $\ell$ containing at least three points
of $\mathcal P$. According to TABLE
\ref{T:polarf}, $\ell$ must contain indeed four points of $\mathcal P$, say $q_1, \ldots, q_4$, and the $\ell$-polar  map $f$ of $\mathcal F$ is completely
determined. It  satisfies
\begin{equation}\label{E:cor}
f^{-1}(q_i) = q_i + 2\,  \widehat{q_i} \quad \text{for } \,\, i=1,\ldots,4 .
\end{equation}

Recall from \cite{Br} that  a  foliation of degree $3$  has at most four singularities
on an invariant line. Therefore $\mathrm{sing}(\mathcal F) \cap \ell = \{ q_1, \ldots,
q_4 \}$.  Lemma \ref{L:numberoflines} implies that  through each $q_i$ there is  a $\mathcal
F$-invariant line $\ell_i$ distinct from $\ell$.

From (\ref{E:cor}) one deduces that
\[
\mathrm{tang}(\mathcal F, \mathcal L_{q_i}) = \ell + \ell_i + C_i
\]
where $C_i$ is a conic (not necessarily reduced nor irreducible) intersecting $\ell$ at $\widehat{q_i}$ with
multiplicity two. Theorem \ref{TC:curvatura} implies  moreover that $C_i$ is $\widehat{\mathcal L_{q_i}}$-invariant.

\begin{claim}
None of the conics $C_i$ is reduced and irreducible.
\end{claim}
\begin{proof}
Aiming at a contradiction, suppose that  $C_1$  is reduced and
irreducible.  Then $C_1$ is a $\widehat {\mathcal L_{q_1}}$-invariant curve of degree two. Corollary \ref{C:riccati} implies that
 $\mathcal P$ is contained in the
 union of $\ell$ and $\ell_1$ and that  $\mathcal P \cap \ell_1$ must have the same cardinality of
$\mathcal P \cap \ell$, that is ${\rm Card} (\mathcal P \cap \ell_1)=4$. Recall from above that  $\ell_1$ is $\mathcal F$-invariant
and  $\ell\cap\ell_1 = q_1$. Let $p_5, p_6$ and $p_7$ be the points of $\mathcal P$ in $\ell_1$ distinct from $q_1$, see Figure \ref{F:222}.

A simple computation shows that $\widehat{q_2} \neq q_1$ and  (\ref{E:cor}) implies that  $q_2$ is contained in at most one $\mathcal F$-invariant line
different from $\ell$. Therefore at least two of the three lines $\overline{q_2 p_5}, \overline{q_2 p_6}$ and $\overline{q_2 p_7}$ are not $\mathcal F$-invariant.
Proposition \ref{P:invariante} combined with item {\it (4)} of Lemma \ref{L:ric}  imply that  two of the three points $p_5,p_6$ and $p_7$ are singularities
of $\widehat{\mathcal L_{q_2}}$. Therefore the singularities of $\widehat{\mathcal L_{q_2}}$ are not aligned. Proposition \ref{P:riccati} tell us that  the only
algebraic leaves of $\widehat{\mathcal L_{q_2}}$ are lines through
$q_2$.   Theorem \ref{TC:curvatura} implies  that ${\rm tang}( {\mathcal F},{\mathcal L}_{q_2})$ is constituted of four lines passing trough $q_2$. Consequently $f^{-1}(q_2)=3\,q_2$ by Lemma \ref{L:numberoflines} contradicting (\ref{E:cor}).
 \end{proof}

\begin{center}
\begin{figure}[h]
\resizebox{1.5in}{1.3in}{
\psfrag{l}[][][1.5]{$\ell $}
\psfrag{l1}[][][1.5]{$\ell_1 $}
\psfrag{q1}[][][1.5]{$q_1 $}
\psfrag{q2}[][][1.5]{$q_2 $}
\psfrag{q3}[][][1.5]{$q_3 $}
\psfrag{q4}[][][1.5]{$q_4 $}
\psfrag{qq2}[][][1.5]{$ \widehat{q_2}   $}
\psfrag{p5}[][][1.5]{$p_5 $}
\psfrag{p6}[][][1.5]{$p_6 $}
\psfrag{p7}[][][1.5]{$p_7 $}
\includegraphics{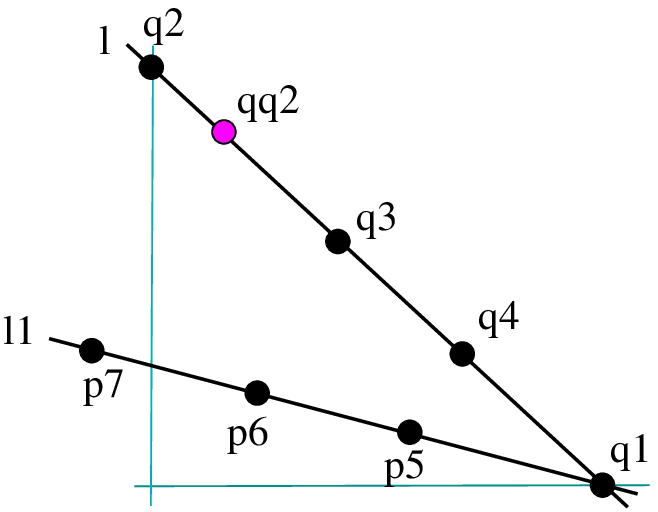}}
\caption{The singularities of $\widehat{\mathcal L_{q_2}}$ are not aligned.}\label{F:222}
\end{figure}
\end{center}

If each $C_i$ is a union of two distinct lines then the linear
system of cubics
\begin{equation}
 \label{E:hesse}
\big\{ \mathrm{tang}(\mathcal F, \mathcal L_p) - \ell \big\}_{p\in \ell}
\end{equation}
contains four totally decomposable fibers. These are triangles
(three lines in general position) with one of the vertices on
$\ell$. This is sufficient (see  \cite[Section 4.4]{Stipins}) to ensure that  (\ref{E:hesse}) is the Hesse
pencil and that $\ell$ is one of its nine harmonic lines. Recall
from \cite{ArtebaniDolgachev} that the Hesse pencil is classically presented as
the one generated by the cubic forms $x^3+y^3+z^3$ and $xyz$. In
these coordinates the harmonic lines are
\[
\begin{array}{lclclc}
\{ x-y = 0 \} && \{ x-\epsilon y = 0 \}  && \{ x-\epsilon^2 y = 0
\} & \\
\{ x-z = 0 \} && \{ x-\epsilon z = 0 \}  && \{ x-\epsilon^2 z = 0
\} &\\
\{ y-z = 0 \} && \{ y-\epsilon z = 0 \}  && \{ y-\epsilon^2 z = 0 \}\, .
&
 \end{array}
\]
The subgroup of $\mathrm{Aut}(\mathbb P^2)$ that preserves the Hesse
Pencil is the Hessian group $G_{216}$  isomorphic to $\left(
\frac{\mathbb Z}{3 \mathbb Z}\right)^2 \rtimes \mathrm{SL}(2,
\mathbb F_3)$. The projective transformations
\[
   a:\, (x:y:z)\mapsto  (y:z:x) \quad \text{ and } \quad b: \,(x:y:z) \mapsto  (x:
   \epsilon y : \epsilon^2 z)
\]
generates a subgroup $\Gamma$ of $G_{216}$ isomorphic to $\left(
\frac{\mathbb Z}{3\, \mathbb Z}\right)^2$ acting transitively on
the set of harmonic lines. Thus we loose no generality by assuming
that $\ell = \{ x-y=0\}$.

Notice that the singular set of $\mathcal F$ contains the base points of the linear system (\ref{E:hesse}).  Thus the singular set of
$\mathcal F$ contains the nine base points of the Hesse
pencil    together with the four fixed points of $f$ on $\ell$.
Since $\mathcal F$ has degree $3$, it has at most $3^2 + 3 + 1= 13$ singular points. Therefore
the singular set of $\mathcal F$ has been completely determined and each of its points has multiplicity
one. In other words the singular scheme of $\mathcal F$ is everywhere reduced.

The main Theorem of \cite{CampilloOlivares} says that a
foliation on $\mathbb P^2$ of degree greater than one is completely
determined by its singular scheme. Therefore  $\mathcal F$ is
determined and it is  equal to the foliation induced by
the $1$-form
\begin{equation}\label{E:omega}
\omega = \left( {x}^{3} +{y}^{3}+1 +6\,x{y}^{2}\right) dx -  \left(
x^3 + {y}^{3}+1 +6\,{x}^{2}y \right) dy .
\end{equation}
Therefore $\mathcal F$ is in case {\it (a)} of the statement.
Concerning the set of points $\mathcal P$ it must be equal to $\{
q_1, \ldots, q_4 \}$. Otherwise Corollary \ref{C:riccati} would imply that  there would exist just one
$\widehat{\mathcal L_{q_i}}$-invariant line through $\widehat{q_i}$ contrary to our assumptions on $C_i$.
A  direct computation shows that $K(\mathcal F \boxtimes \mathcal W( \mathcal P) )=0$.

\smallskip

If at least one of the conics $C_i$ is non-reduced then
\cite[Proposition 3.1]{PY} implies that all the $C_i$'s are double
lines. Therefore $\mathcal F$ is a homogeneous foliation on the
affine chart where $\ell$ is the line at infinity and the
singularity of $\mathcal F$ corresponding to the unique base point
of $ \{ \mathrm{tang}(\mathcal F, \mathcal L_p) - \ell \}_{p\in
\ell}$ is the origin. Thus ${\mathcal F}$ is defined by a homogenous $1$-form with coefficients equal to
 the coefficients of $\ell$-polar map, that is
\[
\mathcal F = \left[  y \big(3x + y\, (1-\xi_3^{{2}})  \big) ^2 dx
+ 3x \left(x + y\,(1-\xi_3^{{2}}) \right)^2 dy
 \right] = \left[ d \big( x y (x+y) (x - \xi_3 y) \big) \right] \, .
\]
A linear change of coordinates envoys $\mathcal F$ to the form presented in  case {\it (b)} of the statement.  Finally, it follows from Proposition \ref{P:sing} that
there are only two possibilities for $\mathcal P$: those  mentioned in the statement of the proposition.
Both cases are exceptional, and in particular flat,  as we have shown in Section \ref{S:action}.
\end{proof}

\subsection{Flat CDQL webs of degree four}\label{S:deg4} Finally we turn our
attention to the flat CQDL webs  $\mathcal F \boxtimes {\mathcal
W}({\mathcal P})$ when ${\rm deg}(\mathcal F)=4$  and the cardinality of
$\mathcal P$ is at least four. Corollary \ref{C:onaline} implies that
$\mathcal P$ cannot be in general position and Theorem
\ref{T:bd} shows that  no four points in $\mathcal P$ are
aligned. Therefore there exists a line $\ell$ such that $\mathcal P
\cap \ell = \{ q_1, q_2, q_3\}$.

According to TABLE \ref{T:polarf} there are only two possibilities
for  the $\ell$-polar map $f$ of $\mathcal F$. In both cases $f$ has
$5$ distinct  fixed points that are cut out by the polynomial
$xy(x+y)(x^2+ xy + y^2)$. In particular,  $\mathcal F$ has exactly
$5$ singular points on $\ell$ according to  Lemma \ref{L:singpolarmap}. Notice that $\mathrm{sing}(\mathcal
F) \cap \ell$ does not intersect $\{ \widehat{q_1},
\widehat{q_2},\widehat{q_3}\}$.

\begin{lemma}\label{L:linear4}
For $i=1,2,3$, the tangency of $\mathcal F$ and $\mathcal L_{q_i}$ is
a union of lines.
\end{lemma}
\begin{proof}
Let's first consider  case (c.2) of TABLE \ref{T:polarf}, that is
$f^{-1}(q_i) = 3 q_i + \widehat{q_i}$ for every $i=1,2,3$. By Theorem \ref{TC:curvatura}, any irreducible component $C$ of
the tangency between $\mathcal F$ and $\mathcal
L_{q_i}$ is invariant by ${\mathcal L}_{q_i}$ or
$\widehat{\mathcal L_{q_i}}$. In the former case  $C$ has to be a line as all the irreducible curves left
invariant by $\mathcal L_{q_i}$. In the latter case,  $C$ is also a line. This follows from   Lemma \ref{L:ric} item {\it (2)}  when
$C$ passes through $q_i$ and from
$f^{-1}(q_i) = 3\, q_i  +
\widehat{q_i}$   when $C$ passes through $\widehat{q_i}$.

\smallskip

We will now deal with  case (c.1) of TABLE \ref{T:polarf}, that is
$f^{-1}(q_i) =  q_i + 3\widehat{q_i}$ for every $i=1,2,3$.
We can assume that
$q_1=p_1=[0:1:0]$, $q_2=p_2=[1:0:0]$,
$q_3=p_3=[1:-1:0]$ and $p_4= [ 0:0:1] \notin \ell$.

We will deal separately two cases: (a)  the cardinality of $\mathcal P$ is four, and (b) the cardinality of
$\mathcal P$ is at least five.

\medskip

\noindent{{\bf Case (a): ${\mathbf k} = \mathbf{Card} (\mathcal P)  =4$.} } In this case we will work in the affine coordinates $(x,y)=[x:y:1]$.
Notice that
\[
\widehat{\mathcal L_{1}} = \left[ d \Big(\frac{(x+2y)^3}{x} \Big) \right] \, , \quad\widehat{\mathcal L_{2}} = \left[d\Big(\frac{(2x+y)^3}{y} \Big)\right]  \quad
 \text{and} \quad \widehat{\mathcal L_{4}} = \Big[ d \big((xy(x+y)\big) \Big].
\]

If we write $\mathcal F = [ a(x,y)dx + b(x,y)dy ]$, where   $a$ and $b$ are relatively prime polynomials, then
$\mathrm{tang}(\mathcal F , \mathcal L_{q_1})$ is defined by the vanishing of $a(x,y)$. Similarly $\mathrm{tang}(\mathcal F , \mathcal L_2)$ is defined
by the vanishing of  $b(x,y)$. Theorem \ref{TC:curvatura}   implies that
\begin{equation}\label{E:???2}
 \mathcal F  = \left[(y-\lambda_1)\left( (2x+y) ^3 - \mu_1 y \right) dx +
(x-\lambda_2)\left( (x+2y) ^3 - \mu_2 x \right) dy \right] \,
\end{equation}
for some $\lambda_1,\lambda_2, \mu_1 ,\mu_2 \in \mathbb C$.

On the one hand,  $\mathrm{tang}(\mathcal F, \mathcal L_4)$ contains the singular
points of $\mathcal F$ on $\ell$. Theorem \ref{TC:curvatura}
implies that its irreducible components must be irreducible cubics
in the pencil $<z^3, xy(x+y)>$ or lines connecting $p_4$ to one of
the $5$ singularities of $\mathcal F$ at $\ell$ (corresponding to
the $5$ fixed points of the $\ell$-polar map of $\mathcal F$).
Thus,
\begin{equation}\label{E:simples}
\mathrm{tang}(\mathcal F, \mathcal L_4) = \big\{ (x^2 + xy + y^2) (
xy(x+y) - \lambda_3 ) = 0 \big\}
\end{equation}
for a certain $\lambda_3 \in \mathbb C$.

On the other hand, the tangency between $\mathcal F$ and $\mathcal L_4$ is
defined by the vanishing of the contraction of the $1$-form in (\ref{E:???2})
with
 $x\partial_x + y \partial_y$. Explicitly
\[
\mathrm{tang}(\mathcal F, \mathcal L_4) = \left\{ x(y-\lambda_1)\left( (2x+y) ^3 - \mu_1 y \right)  +
y(x-\lambda_2)\left( (x+2y) ^3 - \mu_2 x \right)   = 0\right\}.
\]
Comparing the homogeneous components of degree two of the two presentations of $\mathrm{tang}(\mathcal F,\mathcal L_4)$, one concludes that
$\lambda_3=\lambda_1 \mu_1 = \lambda_2 \mu_2 =0$. Plugging $\lambda_3=0$ into (\ref{E:simples}) shows that all
the five lines cut out by  $xy(x+y)(x^2+xy+y^2)$ are  $\mathcal
F$-invariant. The $\mathcal F$-invariance of $\{x=0\}$ and $\{y=0\}$ ensures that  $\lambda_1=\lambda_2=0$. Finally, since
the homogeneous component of degree three of (\ref{E:simples}) is zero, $\mu_1=\mu_2=0$. It is then clear that the expression  of $\mathcal F$ in (\ref{E:???2}) is homogeneous. Consequently $\mathrm{tang}(\mathcal F, \mathcal L_q)$ is a union of lines for every $q \in \ell$.

\medskip

\noindent{{\bf Case (b): $\mathbf{k = {\mathbf Card}(\mathcal P) \ge 5}$.} }
Notice that $\mathcal P$ is not in
barycentric general position with respect to  none of the points
$q_1, q_2, q_3$ because  $f^{-1}(q_i) \neq 4 q_i$ for $i=1,2,3$. Proposition \ref{P:riccati} implies that all the leaves of $\widehat{\mathcal L_i}$ are algebraic.
From the proof of  Corollary \ref{C:riccati}, one deduces  that the leaves of  $\widehat{\mathcal L_1}$ (for instance)
are irreducible components of elements of a pencil of the form
$\mathcal H=< (x+2y + \lambda z)^{\deg(R)},R(x,z)  >$, where
$\lambda \in \mathbb C$ and $R$ is a homogeneous polynomial of
degree $k-1$. The irreducible factors of $R$ correspond to the lines $\overline{q_1p}$ where $p$ ranges in $\mathcal P_1={\mathcal P}\setminus\{q_1\}$ and their
multiplicities correspond to number of points of $\mathcal P_1$
contained in the respective lines.

If $\mathrm{tang}(\mathcal F, \mathcal L_{q_1})$ has an non-linear irreducible component $C$ then  its degree is at most three and  is an irreducible component of
an element of the pencil $\mathcal H$.  But $(x+2y + \lambda z)^{\deg(R)} - \mu R(x,z)$ admits an  irreducible
factor of degree at most three for some $\mu \in \mathbb C^*$ only when  $R$ is a square.  Indeed, on the one hand the square of each linear factor of $R$ must divide $R$ otherwise  Corollary \ref{C:riccati} would imply that $C$ has degree $k-1\ge 4$. On the other hand the third power of any linear factor of $R$ cannot divide $R$, otherwise it would exist four points in $\mathcal P$ on the same line contradicting Theorem \ref{T:bd}.

Since $R$ is a square  it must
exist a third point $p_5 \in \mathcal P$ contained in the the line
$\overline{q_1 p_4}$.
From the fact that  $\mathcal P$ is not in $q_2$-barycentric general position it follows that $\mathrm{sing}(\widehat{\mathcal L_{q_2}})-\{q_2\}$ is contained in a line.
Using that  $\widehat{q_2}\neq q_1$ one deduces that  it must exist
a point $p_6 \in \overline{q_2 p_4}\cap \mathcal P$.  Since $R(x,z)$
is a square, the line $\overline{q_1 p_6}$ must contain another
point of $\mathcal P$ (noted $p_7$ in  Figure \ref{F:incidencia} below).
\begin{center}
\begin{figure}[h]
\resizebox{1.2in}{1.2in}{
\psfrag{l}[][][1.5]{$\ell $}
\psfrag{q1}[][][1.5]{$q_1 $}
\psfrag{q2}[][][1.5]{$q_2 $}
\psfrag{q3}[][][1.5]{$q_3 $}
\psfrag{p4}[][][1.5]{$p_4 $}
\psfrag{p5}[][][1.5]{$p_5 $}
\psfrag{p6}[][][1.5]{$p_6 $}
\psfrag{p7}[][][1.5]{$p_7 $}
\includegraphics{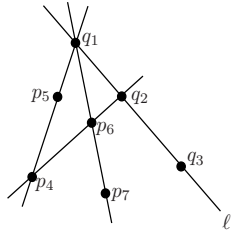}}
\caption{Seven points of $\mathcal P$.}\label{F:incidencia}
\end{figure}
\end{center}\vspace{-0.7cm}
Proposition \ref{P:invariante} tell us that any line containing at least $3$ points of
$\mathcal P$ must be $\mathcal F$-invariant. Thus {there}
are at least three $\mathcal F$-invariant lines through $q_1$. This
contradicts $f^{-1}(q_1) = q_1 + 3\, \widehat{q_1}$ and ends the proof of the lemma.
\end{proof}

We will also need a classical result of Darboux about the degree of foliations induced by
pencil of curves. We state it below as a lemma.

\begin{lemma}\label{L:Darboux}
If  $F, G \in \mathbb C[ x,y,z]$ are relatively prime  homogeneous
polynomials of degree $e$   then
$$FdG - GdF =  \left(\prod_{H} H^{e(H) - 1}\right)  \cdot \omega $$
where $\omega$ is a homogenous polynomial $1$-form with codimension two
singular set;   $H$ runs over the irreducible components of the
polynomials $\left\lbrace s  F + t  G =0 \right\}_{ (s:t) \in  \mathbb P^1 } $; and $e(H)$ denotes de maximum
power of $H$  that divides the member of the pencil that contains $H$. In particular if $\mathcal F = \left[d (F/G) \right]$ then
\[
 \deg(\mathcal F) = 2e - 2 - \sum_{H} {\deg(H)(e(H) - 1)} \, .
\]
\end{lemma}
\begin{proof}
See \cite[Proposition 3.5.1, pages 110--111 ]{Jouanolou}.
\end{proof}

\begin{prop}\label{P:deg4}
Let $\mathcal F$ be a  foliation of degree four   and $\mathcal P
\subset \mathbb P^2$ be a finite set of at least four points. If
$K(\mathcal F \boxtimes \mathcal W(\mathcal P))=0$ then
$\mathcal F$ is projectively equivalent to one of the following foliations:
\begin{enumerate}
\item[(a)] $\displaystyle{\mathcal F= \Big[ d\Big( xy\,(x+y)\,(x^2 + xy +
y^2)^3 \Big)\Big]};$\vspace{0.15cm}
\item[(b)]  $\displaystyle{\mathcal F = \left[  d\Big( \frac{xy\,(x+y)}{x^2+xy+y^2}\Big)  \right]};$
\vspace{0.15cm}
\item[(c)] $\displaystyle{\mathcal F = \left[  d\Big( \frac{x^3+y^3 + 1}{xy}\Big)
\right]}$.
\end{enumerate}
Moreover $\mathcal P = \{ [1:-1:0],[1:0:0],[0:1:0],[0:0:1]\}$  in cases (a) and (b).  In
case (c), $\mathcal P$ is equal to the nine base points of the
pencil  $<xy, x^3+y^3 +1>$ or is equal to the three base
points of this pencil at the line at infinity  union with~$[0:0:1]$.
\end{prop}
\begin{proof}
We keep the notations from the beginning of this section.
According to TABLE \ref{T:polarf} there two possibilities for  the $\ell$-polar map
of $\mathcal F$: (c.1) and (c.2). We will deal with them separately.

\medskip

\noindent{{\bf Case (c.1).}} We are assuming that the $\ell$-polar map of $\mathcal F$ satisfies  $f^{-1}(q_i)
= q_i + 3 \,\widehat{q_i}$ for $i=1, \ldots, 3$. According to Lemma
\ref{L:linear4}, the tangency between $\mathcal F$ and $\mathcal
L_{q_i}$ is a union of lines. Since $\mathcal P$ has cardinality at least four, there exists $p_4\in {\mathcal P} \setminus \ell$.
Moreover $\mathcal P$ is not in $q_i$-barycentric general position.
Proposition \ref{P:riccati} implies that the foliation $\widehat{\mathcal L_{q_i}}$ admits exactly one invariant line $\widehat{\ell_i}$ through $\widehat{q_i}$.
Moreover $f(q_i) = q_i + 3\,\widehat{q_i}$ implies that it exists  exactly one $\mathcal L_{q_i}$-invariant line $\ell_i$ through $q_i$ distinct from $\ell$.  Thus
\[
\qquad
\mathrm{tang}(\mathcal F, \mathcal L_{q_i}) = \ell + \ell_i +
3\, \widehat{\ell_i} \quad \mbox{ for } \; i= 1, \ldots,
3.
\]

If $\mathcal G$ is the foliation induced by the pencil  $\{ (\mathrm{tang}(\mathcal F, \mathcal L_{q})- \ell) \}_{q\in \ell}$
then Lemma \ref{L:Darboux} implies that  $\mathcal G$ has degree at most $2 \cdot 4  - 2 - 3 \cdot (3-1) = 0$.  In  an   affine coordinate system where
$\ell$ is the line at infinity and the origin belongs to $\mathrm{sing}(\mathcal G)$, the foliation $\mathcal F$ is induced by a polynomial $1$-form with homogeneous components. Therefore it is completely determined {by   its} $\ell$-polar map and can be  explicitly presented as
\[
\mathcal F =  \left[ y\,(2x+y)^3 dx + x\,(2y+x)^3 dy \right] .
\]
A simple computation shows that   $xy\,(x+y)(x^2 + xy +
y^2)^3$ is a first integral of $\mathcal F$. Since the  singular set of  $\mathcal F$
has cardinality four it has to be equal to $\mathcal P$.
A direct computation shows that $K(\mathcal F \boxtimes \mathcal W(\mathcal P))=0$.
This example corresponds to case  {\it (a)} of the statement.

\medskip

\noindent{{\bf Case (c.2).}}
Suppose now that the $\ell$-polar map of $\mathcal F$ is in case
(c.2) of TABLE \ref{T:polarf}.  Lemma \ref{L:linear4} implies (for any $i=1,\ldots,3$) that  $\mathrm{tang}(\mathcal
F, \mathcal L_{q_i})$  is the union of five lines:
 $\ell$, one  line through $\widehat{q_i}$ and three lines (counted with
 multiplicities) through $q_i$.
It follows from \cite[Proposition
3.1]{PY} that the  multiplicities appearing in
$\mathrm{tang}(\mathcal F, \mathcal L_{q_i})$ do not depend on the
choice of $i \in \{ 1,2,3\}$. Therefore, if $\mathcal G$ denotes the
foliation associated to the pencil $\{ \mathrm{tang}(\mathcal F,
\mathcal L_p) - \ell\}_{p\in \ell} $ then Lemma \ref{L:Darboux} implies that
the degree of $\mathcal G$ is at most:
\begin{trivlist}
\item[(c.2.1)]  zero when there is one line with multiplicity $3$ in $\mathrm{tang}(\mathcal F, \mathcal L_{q_i})$;\vspace{0.07cm}
\item[(c.2.2)] three when there is one line with multiplicity $2$ in $\mathrm{tang}(\mathcal F, \mathcal L_{q_i})$;\vspace{0.07cm}
\item[(c.2.3)] six when all the lines in $\mathrm{tang}(\mathcal F, \mathcal L_{q_i})$ have multiplicity one.
\end{trivlist}

\medskip

\noindent{{\bf Case (c.2.1).}}
If the degree of $\mathcal G$ is equal to zero then, as in case (c.1) above,  $\mathcal F$ is completely determined by its $\ell$-polar map.
In a suitable affine coordinate system,  the foliation $\mathcal F$ is induced by
\[
  \omega = y^3(2x+y)\, dx + x^3(x+2y)\, dy \, .
\]
One can verify that $\omega$ admits $\frac{xy\,(x+y)}{x^2+xy+y^2}$ as a rational
first integral and, again as in the  case (c.1),  $\mathcal
P=\mathrm{sing}(\mathcal F)$. This example corresponds to case  {\it (b)} of the statement.

\medskip
\noindent{{\bf Case (c.2.2).}}
If the degree of $\mathcal G$ is at most three and distinct from zero then $\mathcal G$  is tangent to a
pencil of quartics with three completely decomposable fibers, each
formed by three distinct lines with one of these lines with multiplicity two. Therefore $\mathcal G$
has at least nine invariant lines. Since a degree $d$ foliation has at most $3d$ invariant lines (see \cite{jvpfourier})
it follows that the degree of $\mathcal G$ is exactly $3$.

 It is not hard to show that, up to  automorphisms of
$\mathbb P^2$, there exists a unique foliation $\mathcal G$  as above. In suitable affine coordinates where  $\ell$ is the
line at infinity  and $q_1=[1:0:0]$, $q_2=[0:1:0]$, $q_3=[1:-1:0]$,
the foliation $\mathcal G$  is defined by the rational function
$$\frac{x^2(x-1)(x+2y-1)}{y^2(y-1)(2x+y-1)}\,.$$
We leave the details to the reader.

It follows that
\[
\mathcal F  =\Big[ y^2 (y-1) (2x + y - 1) dx +  x^2(x-1)(x+2y-1)dy  \Big]  .
\]
By a direct computation, it can checked that the $4$-web $\mathcal F \boxtimes \mathcal W (\{
q_1,q_2,q_3\})$ has curvature zero. Nevertheless a lengthy computation shows that  there is no set $\mathcal P$ verifying  $\{
q_1,q_2,q_3\} \subsetneq \mathcal P \subset \mathrm{sing}(\mathcal
F)$ such that $K(\mathcal F \boxtimes \mathcal W(\mathcal P))=0$.

\medskip

\noindent{{\bf Case (c.2.3).}}
We are now assuming that for each $i=1,\ldots,3$,  $\mathrm{tang}(\mathcal F, \mathcal L_{q_i})$
consists of five distinct lines, four of them being $\mathcal
F$-invariant. It implies that $\mathcal F$ has at least $10$ invariant lines, $\ell$ plus nine
others.

We will further divide this case in two subcases: (c.2.3.a) when $k= {\rm Card }( \mathcal P) = 4$, and (c.2.3.b) when $k= {\rm Card }( \mathcal P )\ge 5$.

\medskip

\noindent{{\bf Case (c.2.3.a).}}
Assume that  $\mathcal P= \{ q_1,q_2,q_3,p\}$ with $p\not \in \ell$. Notice that   $\mathrm{tang}(\mathcal F,
\mathcal L_p)$ intersects $\ell$ at the five singular points of $\mathcal F$ on $\ell$: $q_1, q_2, q_3$ and two other that we will call $s_1$ and $s_2$.
Recall from  Lemma \ref{L:singpolarmap} that these five points coincide with the
fixed points of the $\ell$-polar map of $\mathcal F$. With no loss of generality,
one can assume that the points of ${\mathcal P}$ are normalized such that
$q_1 = [ 1 : -1:0]$, $ q_2=[1: -\xi_3:0]$, $q_3= [1:-\xi^2_3:0]$ and
$p=[0:0:1]$. Then, by  Corollary \ref{C:riccati}, the foliation $\widehat{\mathcal L_p}$ admits ${x^3+y^3}$ as a first integral. Consequently,
any irreducible $\widehat{\mathcal L_p}$-invariant algebraic curve $C$ is of degree less than 3 and satisfies  $C\cap \ell \subset \{ q_1,q_2,q_3\}$.
Observe that  $|\mathrm{tang}(\mathcal F, \mathcal L_p)|$ must contain all the singularities of $\mathcal F$,  in particular $s_1$ and $s_2$.
Because none of the curves $\{ x^3 + y^3 = cst. \}$ contains $s_1$ or $s_2$, Theorem \ref{TC:curvatura} implies that
the lines $\overline{ps_1}$ and $\overline{ps_2}$ are ${\mathcal F}$-invariant irreducible components of $\mathrm{tang}(\mathcal F, \mathcal L_p)$.
Therefore $\mathcal F$ has at least $12$ invariant lines: {$\ell$,  $\overline{ps_1}$,  $\overline{ps_2}$, and the three linear components of  ${\rm tang}(\mathcal F,\mathcal L_i)$ passing trough $q_i$ for each $i\in \{1,2,3\}$.}
It is well-known that a degree $d$ foliation of $\mathbb P^2$ has at most $3d$ invariant lines (see  \cite{jvpfourier} for instance).
Therefore $\mathcal F$ has exactly $12$ invariant lines.
\begin{center}
\begin{figure}[h]
\resizebox{1.3in}{1.3in}{
\psfrag{l}[][][1.5]{$\ell $}
\psfrag{q1}[][][1.5]{$q_1 $}
\psfrag{q2}[][][1.5]{$q_2 $}
\psfrag{q3}[][][1.5]{$q_3 $}
\psfrag{s2}[][][1.5]{$s_2 $}
\psfrag{s1}[][][1.5]{$s_1 $}
\psfrag{p}[][][1.5]{$p $}
\includegraphics{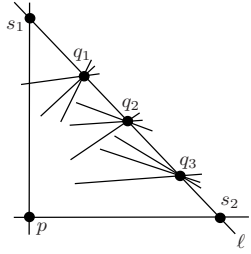}}
\caption{The twelve lines invariant by ${\mathcal F}$.}
\end{figure}
\end{center}
Because $\mathcal F$ has degree $4$, over each $\mathcal F$-invariant lines there are at most $5$ singularities of $\mathcal F$. Notice that over the $\mathcal F$-invariant
line $\overline{ps_1}$ we know already two: $p$ and $s_1$. The three $\mathcal F$-invariant lines through $q_1$ distinct
from $\ell$ must intersect $\overline{ps_1}$ in three distinct singular points of $\mathcal F$, none of them equal to $p$ or  $s_1$ (see Figure \ref{F:lolo} below).
\begin{center}
\begin{figure}[h]
\resizebox{1.3in}{1.3in}{
\psfrag{l}[][][1.5]{$\ell $}
\psfrag{q1}[][][1.5]{$q_1 $}
\psfrag{q2}[][][1.5]{$q_2 $}
\psfrag{q3}[][][1.5]{$q_3 $}
\psfrag{s2}[][][1.5]{$s_2 $}
\psfrag{s1}[][][1.5]{$s_1 $}
\psfrag{p}[][][1.5]{$p $}
\includegraphics{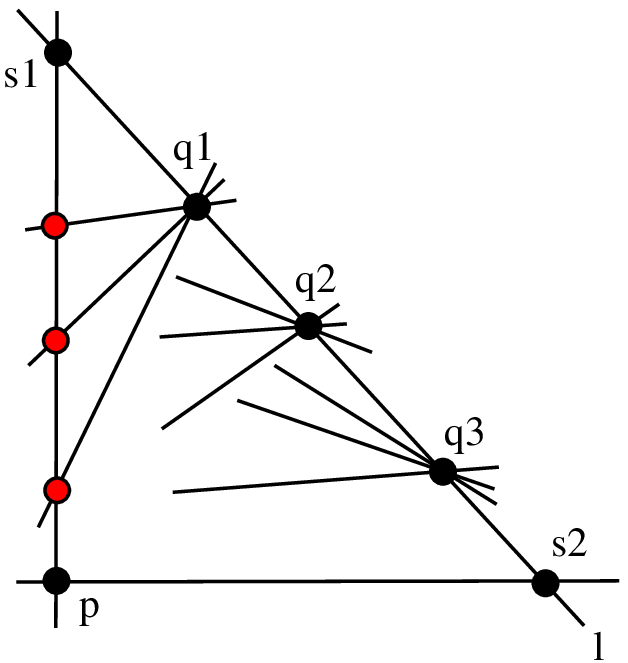}}
\caption{}
\label{F:lolo}
\end{figure}
\end{center}
 The same being true
for the $\mathcal F$-invariant lines through $q_2$ and $q_3$ it follows that on $\overline{ps_1}$ there are three singularities of $\mathcal F$ distinct from $p$ and
$s_1$ such that through each passes four $\mathcal F$-invariant lines. Of course the line $\overline{ps_2}$ has the same property. Thus we have a set $\mathcal Q \subset \mathbb P^2$ of cardinality $9$ such that  each of the points of $\mathcal Q$ {is contained in  four of the the twelve $\mathcal F$-invariant lines}. It is a then a simple
combinatorial exercise to show that these twelve lines support a $(4,3)$-net in the sense of Section \ref{S:net}. Therefore (see \cite[Section 4.4]{Stipins}) the arrangement of twelve $\mathcal F$-invariant lines is projectively equivalent to the Hesse arrangement.  Because the foliation determined by the  Hesse pencil also has degree   four and the tangency of two distinct foliations of degree four has degree nine it follows that $\mathcal F$ is the Hesse Pencil.

With the normalizations made above on the points $q_1,q_2,q_3$ and $p$, we obtain
\[
\mathcal F = \left[  d\Big( \frac{x^3+y^3 + 1}{xy}\Big) \right] .
\]
This $5$-web appeared in the introduction under the label $\mathcal H_5$. In Section \ref{S:net} it is shown that it is an exceptional web and in
particular has curvature zero.

\medskip

\noindent{{\bf Case (c.2.3.b).} }
Suppose now that  $\mathcal P$ has cardinality greater than
four. As in case (c.2.3.a) we will denote by $s_1$ and $s_2$ the two other singularities of $\mathcal F$ on $\ell$  distinct from $q_1,q_2$ and $q_3$.

\begin{claim}
\label{Cl}
There exists a pair of points $p,s \in \mathcal P \setminus \{ q_1, q_2, q_3 \}$  such that the
line $\overline{ps}$ intersects $\ell$ in one of the points $q_1,q_2,q_3$.
\end{claim}
\begin{proof}
Suppose that the claim is not true and let $p_4, p_5$ be any two points in $\mathcal P \setminus \{ q_1, q_2, q_3 \}$.
Proposition \ref{P:invariante} combined with Theorem \ref{T:bd} implies that   the line  $\overline{p_4p_5}$ intersects $\mathcal P$ in
at most three points. Thus there are only two possibilities for
$\mathcal P$: ({\it i}) $\overline{p_4p_5}\cap \mathcal P = \{ p_4,p_5\}$ or
({\it ii}) $\overline{p_4p_5}\cap \mathcal P = \{ p_4,p_5,p_6\}$ for some point $p_6 \in \mathcal P$ distinct of $p_4$ and $p_5$.

If we are in case ({\it i}) then  $\mathcal P$ is in $p_4$ and $p_5$-barycentric general position because, by assumption, the lines $\overline{p_4q_i}$ and
$\overline{p_5q_i}$ (for $i=1,2,3$) have only two elements of $\mathcal P$ each and the points $q_1,q_2,q_3$ are not aligned with $p_4$ nor with $p_5$.
Theorem \ref{TC:curvatura}  ensures that  $|\mathrm{tang}(\mathcal F, \mathcal L_{4})|$ is a union of five $\mathcal F$-invariant lines.
Since $|\mathrm{tang}(\mathcal F, \mathcal L_{4})|$ contains $p_4$ and the singularities of $\mathcal F$,
these lines have to be  $\overline{p_4 s_1},\overline{p_4 s_2},\overline{p_4 q_1},\overline{p_4 q_2}$ and $\overline{p_4 q_3}$.
Similarly the irreducible components of $\mathrm{tang}(\mathcal F, \mathcal L_{p})$ are   the $\mathcal F$-invariant lines
$\overline{p_5 q_i}$ for $i=1,\ldots, 3$ and $\overline{p_5 s_i}$ for $i=1,2$.

Through at least one of the points $s_1,s_2$, say $s_1$,
passes three $\mathcal F$-invariant lines:  $\overline{p_4 s_1}$, $\overline{p_5 s_1}$
and $\ell$.
This  contradicts the behavior of the $\ell$-polar map {because $s_1$ appears in $f^{-1}(s_1)$ with multiplicity one as a simple computation shows}.

Suppose now that we are in case ({\it ii}). Because the barycenter transform of three distinct points in $\mathbb P^1$ is still three distinct points,
$\mathcal P$ is in
 barycentric general position with respect to  at least
two points in $\{ p_4,p_5,p_6\}$. Exactly as before we arrive at a contradiction. The claim follows. \end{proof}

By  Claim \ref{Cl} we can suppose that $p_4,p_5$ are two points in $\mathcal P \setminus \ell$ such that the line $\ell' = \overline{p_4p_5}$ intersect $\ell$ at $q_1$.
Notice that $\ell'$  is $\mathcal F$-invariant (by Proposition \ref{P:invariante}) and  that the $\ell'$-polar map of $\mathcal F$ must  also be in  case (c.2) of TABLE \ref{T:polarf}. Therefore Lemma \ref{L:numberoflines} implies that through each of the points  $p_4$ and $p_5$ passes four $\mathcal F$-invariant lines. Since these intersect $\ell$ at $\mathrm{sing}(\mathcal F)$,  there will be one $\mathcal F$-invariant line through $s_1$ (say $\overline{p_4 s_1}$) and one through $s_2$, say $\overline{p_5 s_2}$. In the total $\mathcal F$ has the maximal number of invariant lines for a degree $4$ foliation: twelve.
\begin{center}
\begin{figure}[h]
\resizebox{1.5in}{1.5in}{
\psfrag{l}[][][1.5]{$\ell $}
\psfrag{lp}[][][1.5]{$\ell' $}
\psfrag{q1}[][][1.5]{$q_1 $}
\psfrag{q2}[][][1.5]{$q_2 $}
\psfrag{q3}[][][1.5]{$q_3 $}
\psfrag{s2}[][][1.5]{$s_2 $}
\psfrag{s1}[][][1.5]{$s_1 $}
\psfrag{p4}[][][1.5]{$p_4 $}
\psfrag{p5}[][][1.5]{$p_5 $}
\includegraphics{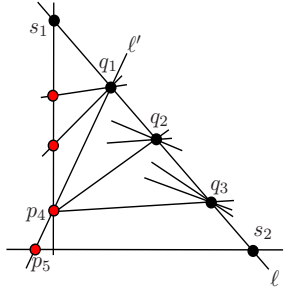}}
\caption{The twelve lines invariant by ${\mathcal F}$ in case {\bf c.2.3.b}.}
\end{figure}
\end{center}\vspace{-0.5cm}
Consider the {effective divisor} $\mathrm{tang}(\mathcal F, \mathcal L_4)$. It has degree $5$,  contains four lines through $p_4$ (namely  $\overline{p_4 q_2}$, $\overline{p_4 q_3}$, $\overline{p_4 s_1}$ and $\ell'=\overline{p_4 q_1}=\overline{p_4 p_5}$) and  the point $s_2$. Since the four lines through $p_4$ do not contain $s_2$ there
is a line $\ell'' \subset|\mathrm{tang}(\mathcal F, \mathcal L_4)|$ through $s_2$. By Theorem \ref{TC:curvatura}, $\ell''$ must be $\widehat{\mathcal L_4}$ invariant
and
Lemma \ref{L:ric} item {\it (4)} implies that  $\ell''$ contains $\widehat{p_4}$:  the $p_4$-barycenter of $\{p_5, q_1\}$  in $\ell'$. In particular, $\ell''=\overline{s_2 \widehat{p_4}}$. Clearly $q_2 \notin \ell''$. Consequently Lemma \ref{L:ric} item {\it (4)} ensures the existence of
 an  extra point in $\mathcal P$, say $p_6$, such that
$p_6 \in \overline{p_4 q_2}$ and the $p_4$-barycenter of $\{q_2,p_6\}$ in  $\overline{p_4 q_2}$ lies in $\ell''$. Similarly, there exists another extra point  $p_7 \in \mathcal P$ contained in $\overline{p_4 q_3}$
such that the $p_4$-barycenter of $\{q_3,p_7\}$ in  $\overline{p_4 q_3}$ also lies in $\ell''$.

Notice that the line  $\overline{p_4q_2}$ contains three points of $\mathcal P$: $q_2,p_4$ and $p_6$. Therefore the $\overline{p_4q_2}$-polar map of
$\mathcal F$  must be  also  in the case (c.2) of TABLE \ref{T:polarf}. Consequently through $p_6$ {pass} four $\mathcal F$-invariant lines. Remark that $p_4$,$p_6$ and $s_1$
are  not aligned and that through $s_1$ pass just two $\mathcal F$-invariant lines ($\ell$ and $\overline{p_4s_1}$). Thus one of the four $\mathcal F$-invariant lines through $p_6$ must be  the line $\overline{p_6s_2}$.  Similarly, through $p_7$ {pass} four $\mathcal F$-invariant lines and
the line $\overline{p_7s_2}$ is among these four lines. Since through
$s_2$ passes just one $\mathcal F$-invariant distinct from $\ell$ it follows that $\overline{p_6s_2}=\overline{p_7s_2}=\overline{p_5s_2}$.

Changing the role of $p_4$ and $p_5$ in the preceding argument it follows that there exist $p_8, p_9 \in \mathcal P \setminus \{ q_1, q_2,q_3,p_4, \ldots, p_7 \}$
in the lines $\overline{p_5 q_2}$ and $\overline{p_5q_3}$ respectively.  As before,  through each of these points passes four $\mathcal F$-invariant lines.

Putting all together we have just proved that $\mathcal F$ leaves invariant an arrangement  of twelve lines and $\mathcal P$ contains a subset of
at least nine points such that  each of these points is contained in  four  distinct lines of the arrangement. At this point it is clear that the arrangement is
the Hesse arrangement (see \cite{Stipins}), that  $\mathcal F$ is projectively equivalent to the Hesse pencil (it is the unique degree $4$ foliation leaving the Hesse arrangement invariant because the tangency of two degree four foliations has degree nine) and that $\mathcal P$ contains the nine base points of it. It remains
to show that $\mathcal P$ cannot be larger than the base points of the Hesse pencil. Indeed if there exists a point $p_{10} \in \mathcal P$ distinct from
the nine base points it would exist a line in the arrangement containing four points of $\mathcal P$ contradicting Theorem \ref{T:bd}. Therefore there exists only one
flat CDQL $(k+1)$-web of degree four with $k \ge 5$: the
$10$-web  $\mathcal H_{10}$ from the Introduction.
\end{proof}

\subsection{Proof of Theorem \ref{T:2}}
According to Section \ref{S:deg1}, the exceptional  CDQL webs 
 of degree one are projectively equivalent to one of the webs $\mathcal A^k_I,\mathcal A^k_{II},\mathcal A^k_{III}, \mathcal A^k_{IV}$.

Propositions \ref{P:deg2},\ref{P:deg3},\ref{P:deg4} putted together give  a complete classification of  flat  CDQL $(k+1)$-webs
of degree bigger than two, on the projective plane, when $k \geq 4$. There are only sixteen such webs (up to projective transformations).
Thirteen of these have been presented in the Introduction and their exceptionality has been put in evidence in Sections \ref{S:action} and  \ref{S:net}.

It can be verified the  that $5$-web  described in Proposition \ref{P:deg2} case {\it (a.3.h)}, the $5$-web described in Proposition \ref{P:deg3} case {\it (a)} and the $5$-web described in Proposition \ref{P:deg4} case {\it (a)} are not exceptional. For this sake one can use, as we did,  the criterion \cite[Proposition 4.3]{PT} or H\'{e}naut's curvature as indicated by Ripoll in \cite[Theorem 5.1]{Ri} or even Pantazi's criterion. Aiming at conciseness we decided not to reproduce the  lengthy computations here.
\qed

\begin{remark}\rm As already mentioned in the Introduction  the non-linear defining foliation of all the exceptional CDQL webs on $\mathbb P^2$ admits a rational
first integral. {The  non-linear} defining foliations of the flat CDQL $5$-webs of degrees two and four also admit rational first integrals. The situation is different
for the flat CDQL $5$-web of degree three.

\begin{prop}
The foliation $\mathcal F=\left[ y\,(2x+y)^3 dx + x\,(2y+x)^3 dy \right]$  does not admit a rational first integral.
\end{prop}
\begin{proof}
Let $\omega = y\,(2x+y)^3 dx + x\,(2y+x)^3 dy$. If we set $X= x+y$ and $Y= x-y$ then  $\omega = 2(X^3+1) dY + ({3}/{2})( Y^3 - X^2 Y) dX$.
It can be promptly verified that the algebraic function $ Y^3\sqrt{X^3 +1}$ is an integrating factor of $\omega$, that is,
$\omega / (Y^3\sqrt{X^3 +1} )$ is a closed $1$-form. If we set $Z = \sqrt{X^3 + 1}$ then
\[
\frac{\omega}{Y^3Z } = -d \left( \frac{Z}{Y^2} \right) +  \frac{dZ}{(Z^2-1)^{\frac{2}{3}}}\, .
\]

Let $\wp(t)$ be a Weierstrass function satisfying
$ \wp'(t)^2 = 4\wp(t)^3 + 1 \, .$ Notice that $\wp''(t) = 6 \wp(t)^2$.
If we set $Z=\wp'(t)$ then, in the coordinates $(t,Y)$,
\[
\frac{\omega}{Y^3 \wp'(t) } =  -d \left( \frac{\wp'(t)}{Y^2} \right) +  \frac{3}{2^{\frac{1}{3}}}dt \, .
\]
Hence the  generic leaf of $\mathcal F$ is defined by the equation
\[
 - \frac{\wp'(t)}{Y^2} + \frac{3}{2^{\frac{1}{3}}} t = \lambda \quad \implies \quad Y = \sqrt{\frac{\wp'(t)}{\frac{3}{2^{\frac{1}{3}}}t+\lambda}}
\]
where $\lambda$ is a constant.  Since $\wp'(t)$ is doubly-periodic while $(3/{2^{\frac{1}{3}}})t+\lambda$ is not,  it follows
that back in the coordinates $(X,Y)$ the generic leaf of $\mathcal F$ cuts the lines with constant $X$ coordinate
in infinitely many distinct points. In particular they are non-algebraic and therefore $\mathcal F$ does not admit a rational first integral.
\end{proof}
\end{remark}


\section{From global to local\ldots}\label{S:globloc}

\subsection{Degenerations}
Let $\mathcal W_t$ be a holomorphic family of webs in the sense that
it is defined by an element
\[
W(x,y,t) = \sum_{i+j=k} a_{ij}(x,y,t) dx^i dy ^j
\]
in $\mathrm{Sym}^k \Omega^1(\mathbb C^2)$ with coefficients  in $\mathcal O=\mathbb C\{x,y,t\}$ (convergent power series) and such that
   $ W(\cdot,\cdot ,t)$ defines a (possibly singular) $k$-web on $(\mathbb C^2,0)$ for every $t \in (\mathbb C,0)$.

We do not claim originality on the next result. Indeed the first author, modulo memory betrayals, first heard about it
in a talk delivered by H\'{e}naut at CIRM in 2003. Anyway it follows  almost immediately from the main result of \cite{Henaut}. Since it would
take us too far afield to recall  the notations and  the results of \cite{Henaut}, we include a sketchy proof below freely using them.
We refer to this work for more precisions.

\begin{thm}\label{T:MR}
The set
$ \{t \in (\mathbb C,0) \,|\, \mathcal W_t \text{ has maximal rank}\, \}$  is closed.
\end{thm}
\begin{proof}
 The differential system $M_t(d)$ can be defined
over ${\mathcal O}$ (with $t$ considered as a constant of derivations) and the
restriction of $M_t(d)$ to a parameter $t_0$ coincides with the
definition of $M_{t_0}(d)$.

 The prolongations $p_k$ of the associated morphism are  morphisms
of $\mathcal O$-modules and the  kernels $R_k$ of the morphisms   $p_k$ are $\mathcal O$-modules locally free outside
the discriminant. Notice that the discriminant is a hypersurface in $(\mathbb C^2,0) \times (\mathbb C,0)$ that does not contain
any  fiber of the projection $(x,y,t) \mapsto t$  by our definition of family of webs.

 If $r_k = \dim R_k$ then Cartan's Theorem B  implies the existence of $r_k$
sections of $R_k$ over a polydisk $D \subset  (\mathbb C^2,0)\times (\mathbb C,0)$  that generates  $R_k$  on
a Zariski open subset of $(\mathbb C^2,0)\times (\mathbb C,0)$. Moreover this subset can be supposed to contain
any given point on the complement of the discriminant. Therefore we can find a meromorphic  inverse of the morphism
$\overline{\pi}_{k-4}$  holomorphic at any given point
in the complement of the discriminant.

Following \cite{Henaut}, we can construct a holomorphic family of  meromorphic
connections $\Delta_t$ such that $\mathcal W_t$ has maximal rank if and only  if $\Delta_t^2 = 0$. The theorem follows.
\end{proof}

\subsection{Singularities of certain  exceptional webs}

Theorem \ref{T:MR} combined with the classification of CDQL exceptional webs in $\mathbb P^2$ yields
the following result.

\begin{cor}[Corollary \ref{CC:local} of the {I}ntroduction]\label{C:local}
Let $\mathcal W$ be a smooth $k$-web, $k\ge 4$, and  $\mathcal F$ be a
singular holomorphic foliation, both  on $(\mathbb C^2,0)$,
such that the $(k+1)$-web $\mathcal W \boxtimes \mathcal F$ has maximal rank.
Then one of the following holds:
\begin{enumerate}
\item the foliation $\mathcal F$ is of the form $\big[ {H}(x,y) ( \alpha d x + \beta dy) + {h.o.t.} \big]$
where $H$ is a {non-zero} homogeneous polynomial and $(\alpha, \beta) \in \mathbb C^2 \setminus \{ 0 \}$;\vspace{0.1cm}
\item the foliation $\mathcal F$ is of the form $\big[  {H}(x,y) ( y d x - x  dy) + {h.o.t.} \big]$
where $H$ is a {non-zero} homogeneous polynomial;\vspace{0.1cm}
\item $\mathcal W \boxtimes \mathcal F$ is exceptional and its first non-zero jet is one of the following  webs
\[
\mathcal A_I^k, \,  \mathcal A_{III}^{k-2}, \,  \mathcal A_5 ^d\;  (\text{only when } k=4 ) \, \mbox{ and } \,  \mathcal A_6^b \; (\text{only when } k=5) \, .
\]
\end{enumerate}
\end{cor}
\begin{proof}
Suppose that $\mathcal W= [\Omega]$ where $\Omega$ is a germ at the origin of a holomorphic $k$-symmetric $1$-form.
Consider the expansion of $\Omega$ in its homogeneous components:
\[
 \Omega = \sum_{i=0}^{\infty} \Omega_i \,
\]
where $\Omega_i$ is a  $k$-symmetric $1$-form with homogenous coefficients of  degree $i$.
According to our assumptions $\Omega_0 \neq 0$  and
$
 \Omega_{0} = \prod_{i=1}^{k} dL_i \, ,
$
where the  $L_i$'s are linear forms defining the tangent spaces of the leaves of  $\mathcal W$ at {the origin}.

Similarly, suppose that $\mathcal F = [ \omega ]$ where $\omega$ is a germ of holomorphic $1$-form with codimension two
zero set. Let
\[
\omega = \sum_{i=i_0} ^{ \infty} \omega_i \, ,   \quad \omega_{i_0} \neq 0
\]
be the expansion of $\omega$ in its homogeneous components, with $i_0>0$ according to the hypothesis made on ${\mathcal F}$. If $\alpha_t(x,y)=(tx,ty)$ then
\[
W(x,y,t) = \frac{\alpha_t^* ( \Omega \cdot \omega )}{t^{k+i_0+2}}  = \left( \sum_{i=0}^{\infty} t^{  i} \Omega_i \right) \left( \sum_{i=i_0}^{\infty} t^{ i -i_0} \omega_i \right)
=  \Omega_0 \cdot \omega_{i_0} + t( \cdots )
\]
is an element of $\mathrm{Sym}^k \Omega^1(\mathbb C^2)$ with coefficients  in $\mathcal O=\mathbb C\{x,y,t\}$.
For every $t\neq 0$, the
web $W_t=[W(\cdot,\cdot, t)]$ is isomorphic to $\mathcal W \boxtimes \mathcal F$.

If $\omega_{i_0}$ is a multiple  of a constant $1$-form (equivalently if $\mathcal F_0 = [ \omega_{i_0} ]$ is a smooth foliation)  then  $\mathcal F$ must be like in item (1) of the statement. Notice that when $W(x,y,0)$ does not
define a $(k+1)$-web we are in this situation. Otherwise the foliation $\mathcal F_0= [\omega_{i_0}]$ has a singularity at {the origin} and $\mathcal W(x,y,0)$ is a $(k+1)$-web.
Since for every $t \neq 0$ the web $W_t $ is {of maximal rank},  $W_0$ also has  maximal rank thanks to Theorem \ref{T:MR}. If $\mathcal F_0$ is linear then we are in case (2) of the statement. Otherwise $W_0=  [ \Omega_0] \boxtimes \mathcal F_0$ is the product of a parallel $k$-web with a non-linear foliation. Since $k\ge 4$, Proposition \ref{P:exc} implies that $W_0$ is exceptional.
  Therefore it must be one of the thirteen sporadic exceptional  CDQL webs  or belong to one of the four infinite families of exceptional CDQL webs. The only
  ones that are the product of a parallel web with a non-linear foliation are listed in {\it  (3)}.
\end{proof}


\section{\ldots and back: quasi-linear webs on  complex tori}
\label{S:oncomplextori}

\subsection{First integrals of linear foliations on tori}
Let $T$ be a two-dimensional complex torus.  The set of linear foliations
on $T$  is naturally identified with the 1-dimensional
projective space $\mathbb{P}{H}^0(T,\Omega^1_T)$. We are interested in the set
$\mathcal I(T) \subset\mathbb{P} {H}^0(T,\Omega^1_T)$ corresponding to
linear foliations which admit a holomorphic  first integral.

\begin{prop}\label{P:abel}
The cardinality $i(T)$ of $\mathcal I(T)$ is $0$, $1$, $2$ or $\infty$. Moreover
\begin{enumerate}
\item If $i(T)=0$ then $T$ is a simple complex torus;
\item If $i(T)=1$ then $T$ is a non-algebraic complex torus;
\item If $i(T)=2$ then $T$ is isogenous to the product of two
non-isogenous elliptic curves;
\item If $i(T)=\infty$ then $T$ is isogeneous to the square of
an elliptic curve $E$. Moreover if $\omega_1, \omega_2$ is a pair of
linearly independent  $1$-forms on $T$ admitting rational first integrals then
\[
\big\{ \lambda \in \mathbb C \, \,  \big\vert \;  \omega_1 + \lambda \,\omega_2 \,
\, \, \text{has a holomorphic   first integral}\, \big\} = \mathrm{End}(E)
\otimes \mathbb Q. \]
\end{enumerate}
\end{prop}
\begin{proof}
Let $\mathcal F$ be  a linear foliation on $T$. It is induced by a $1$-form with constant coefficients
$\omega = s dx + t dy$ on ${\mathbb C}^2$ viewed here as  the universal covering of $T$.

Notice that $\omega$ is invariant by the action of $T$ on itself. Therefore, since this action is transitive,
$\mathcal F$ admits a compact leaf if and only if it has a compact leaf through $0$. Notice also that
a compact leaf is  nothing more than a subtorus of $T$. Reciprocally if $T$ contains a subtorus $T'$
then translations of $T'$ by elements in $T$  form a linear foliation on $T$
admiting  a holomorphic  first integral given by the  quotient map $T\rightarrow T/T'$.

Therefore if $i(T)$ is equal to zero, $T$ has no closed subgroups of dimension one that is, $T$ is a simple complex torus.
If $i(T)$ is equal to one then $T$ admits exactly one closed subgroup of dimension one. It implies that $T$ is non-algebraic otherwise $T$ would be isogeneous to a product of two {elliptic} curves (according to Poincar\'{e}'s reducibility
Theorem) then  would be such that  $i(T)>1$.
If $i(T)=2$ then $T$ admits two closed subgroup $T'$ and $T''$ of dimension one. The natural map
\begin{eqnarray*}
(x,y)\in T' \times T'' \longmapsto x + y \in T
\end{eqnarray*}
has finite kernel equal to $T' \cap T''$ therefore is an isogeny between $T' \times T''$ and $T$. Notice that $T'$ can't be isogenous
to $T''$ otherwise $\mathcal I(T) = \mathcal I(T' \times T'') =  \mathcal I( T' \times T')$ and the latter set has infinite cardinality
since it {is invariant under the induced   action of  $\mathrm{Aut}(T'\times T') \supseteq \mathrm{PSL}(2,\mathrm{End}(T'))\supseteq \mathrm{PSL}(2,\mathbb Z)$ on $\mathbb{P} {H}^0(T,\Omega^1_T)\simeq  \mathbb P^1$}.

If $\mathcal I(T)$ has cardinality at least three then there exist three pairwise distinct subtorus $T'$,$T''$ and $T'''$ passing through the origin of $T$. As before one get that $T$   is isogenous to $T' \times T''$. The existence of  the natural projections $T''' \to T / T'$ and $T''' \to T / T''$ implies that
all the three curves are isogenous. Moreover, up to an  isogeny, $T$ can be assumed to be $T' \times T'$ with $T'$, $T''$ and $T'''$ identified with
the horizontal, vertical and diagonal subtori respectively. It follows  that $\mathcal I(T)$ is an orbit of the natural action of $\mathrm{PGL}(2,\mathrm{End}(T'))  $, hence $i(T)=\infty$.
\end{proof}

\begin{remark} \rm
Item {\it (4)} of Proposition \ref{P:abel} can be traced back to Abel, see \cite[\S X]{Abel2}. According
to Markushevich \cite[p. 158]{LRM}, it is the first appearance of the so-called {\it complex multiplication}
in the theory of elliptic functions.
\end{remark}

\begin{lemma}\label{L:cross}
Let $T$ be a complex torus isogeneous to the square of an elliptic curve $E$. If $[\omega_1], \ldots, [\omega_4] \in \mathbb P H^0(T,\Omega^1_T)$
are linear foliations {on T} with holomorphic first integral then the cross-ratio $([\omega_1],[\omega_2]:[\omega_3],[\omega_4])$ belongs
to $\mathrm{End}(E)\otimes \mathbb Q$.
\end{lemma}
\begin{proof}
According to the proof of Proposition \ref{P:abel} we can assume that $T = E \times E, \omega_2=dx-dy, \omega_3=dy$ and $\omega_4=dx$. Since
the leaves of $\omega_1$ are algebraic they must be translates of  $E_{\alpha,\beta}$ (defined by  (\ref{E:Ealphabeta}) in Section \ref{S:ell57}) for suitable $\alpha, \beta \in \mathrm{End}(E)$. Thus
$\omega_1 = [ \beta dx - \alpha dy ]$. Therefore
\[
([\omega_1],[\omega_2]:[\omega_3],[\omega_4]) = \frac{\beta}{\alpha} .
\]
The lemma follows.
\end{proof}

\subsection{Flat CDQL  webs on complex tori}
Let $\mathcal W$ be a linear $k$-web on $T$. Clearly it is  a completely decomposable
web. Thus we can write
$\mathcal W=\mathcal L_1 \boxtimes \cdots \boxtimes \mathcal L_k$ where the $\mathcal L_i$'s are
linear foliations. For $i=1,\ldots,k$, set
$
\widehat{\mathcal L_i} = \beta_{\mathcal L_i} ( \mathcal W - \mathcal L_i)$.
and  define the {\it polar map} of
a foliation  $\mathcal F$  on $T$ as the rational map
$
P_{\mathcal F} : T  \dashrightarrow \mathbb P  H^0(T, \Omega^1_T)
$
characterized by the property
\[
 P_{\mathcal F}^{-1} ( \mathcal L) = \mathrm{tang}(\mathcal F, \mathcal L)
\]
for every $\mathcal L \in \mathbb P  H^0(T, \Omega^1_T)$.

Recall from the Introduction {that} a  fiber of a rational map from a two-dimensional complex torus {onto} a curve is linear if it is
set-theoretically equal to a union of subtori.

\begin{lemma}\label{L:abelian}
Let $\mathcal W= \mathcal L_1 \boxtimes \cdots \boxtimes \mathcal L_k$ be a linear $k$-web  on  $T$, with $k\geq 2$. If $\mathcal F$ is a non-linear
foliation  on $T$ such that  $K(\mathcal W \boxtimes \mathcal
F)=0$ then the rational map $P_{\mathcal F}$ has at least $k$ linear fibers, one for each  $\mathcal L_i$.
Moreover, if $k \ge 3$ then  each of the fibers $P_{\mathcal F}^{-1} ( \mathcal L_i)$  {contains}
at least one elliptic curve invariant by $\mathcal L_i$ and at least one invariant by $\widehat{\mathcal L_i}$.
\end{lemma}
\begin{proof}
By Theorem \ref{T:curvaturatoro}, any irreducible component of $\mathrm{tang}(\mathcal F, \mathcal L_i)$ is $\mathcal L_i$ or $\widehat{\mathcal L_i}$-invariant. Since $\mathcal L_i$ and $\widehat{\mathcal L_i}$ are linear foliations, it follows   that   the
fibers $P_{\mathcal F}^{-1} ( \mathcal L_i)$  are  linear for  $i=1,\ldots,k$. This  proves the first part of the lemma.

Suppose now that $k \ge 3$. Aiming at a contradiction, assume that all the irreducible components of
$\mathrm{tang}(\mathcal F, \mathcal L_1)$ are $\widehat{\mathcal L_1}$-invariant.   Proposition \ref{P:abel}
implies that $\widehat{\mathcal L_1}$ is tangent to an elliptic fibration.

Since both $K_T$ and $N \mathcal L_i$ are trivial,
$
\mathcal O_{T}(\mathrm{tang}(\mathcal F,\mathcal L_i) )=K_T\otimes  N
\mathcal F \otimes N \mathcal L_{{i}}= N \mathcal F
$
for every $i=  1, \ldots, k$.
{Taking $i=1$, we get} that  $N \mathcal F$ is linearly equivalent to a divisor supported on some fibers
of the fibration $\widehat{\mathcal L_1}$.
{Taking  $i=2,\ldots, k$}, we see that the divisors  $\mathrm{tang}(\mathcal F, \mathcal L_i)$  are  linearly
equivalent to  $N\mathcal F$ and consequently to $\mathrm{tang}(\mathcal F, \mathcal L_1)$.
Therefore, being all of them effective, they  also have to be   supported
on elliptic curves invariant by  $\widehat{\mathcal L_1}$.

Since two distinct linear foliations on $T$ are everywhere transversal, Theorem \ref{T:curvaturatoro} implies that for every $i=  2, \ldots, k$,
$\mathcal L_i$ or $\widehat{\mathcal L_i}$ is equal to $\widehat{\mathcal L_1}$. By hypothesis the linear foliations $\mathcal L_1, \ldots,
\mathcal L_k$ are pairwise distinct. Therefore  at least $k-1$ of the foliations $\widehat{\mathcal L_i}$ ($i=1,\ldots, k$) coincide.
This contradicts Lemma \ref{L:mult}.

If one assumes that all  the irreducible components of
$\mathrm{tang}(\mathcal F, \mathcal L_1)$ are invariant by ${\mathcal L_1}$ then the same argument with minor modifications
also leads to a contradiction. The lemma follows.
\end{proof}

\begin{prop}\label{P:iso}
Let $\mathcal W= \mathcal L_1 \boxtimes \cdots \boxtimes \mathcal L_k$ be a linear $k$-web on  $T$, with $k \ge 3$.
If $\mathcal F$ is a non-linear
foliation  on $T$ such that  $K(\mathcal W \boxtimes \mathcal
F)=0$ then
\begin{enumerate}
\item $T$ is isogenous to the square of an elliptic curve. In particular $T$ is an abelian surface;
\item the foliations $\mathcal L_1, \ldots, \mathcal L_k$  are tangent to elliptic
fibrations;
\item the foliations $\widehat{\mathcal L_1}, \ldots, \widehat{\mathcal L_k}$  are tangent to elliptic
fibrations;
\item $\mathcal P_{\mathcal F}$ has $k$ linear fibers.
\end{enumerate}
\end{prop}
\begin{proof}
 The points {\it (2)}, {\it (3)} and {\it (4)} follow from Lemma \ref{L:abelian}
since a linear foliation on $T$ is tangent to an elliptic fibration if and only if  it leaves an elliptic curve invariant.
Since $k \ge 3$, Proposition \ref{P:abel} implies {\it (1)}.
\end{proof}

\subsection{On the number of linear fibers of a pencil on a complex torus}
Let $F: T \dashrightarrow \mathbb P^1$ be a meromorphic  map on a
two-dimensional complex torus $T$. We are interested in the number
$k$ of linear fibers of $F$.

\begin{thm}[Theorem \ref{T:PYontori} of the Introduction]\label{T:bound}
If $k$ is finite then $k\le 6$. Moreover, if $k=6$ then every fiber of $F$ is reduced.
\end{thm}
\begin{proof}
If $x,y$ are  homogeneous coordinates on ${\mathbb P}^1$   then
$xdy -ydx\in H^0({\mathbb P}^1,\Omega^1_{{\mathbb P}^1}\otimes {\mathcal O}_{\mathbb P^1}(2))$.
Therefore  $\omega =
F^*(xdy -ydx) \in  H^0({T},\Omega^1_{{T}} \otimes \mathcal
N^{\otimes 2})$ with $\mathcal N = F^* \mathcal O_{\mathbb P^1}(1)$.

Let also  $X \in  H^0({T},T_{T} \otimes ({\mathcal
N^*)^{\otimes 2}})$ be dual to $\omega$, that is  $ \omega = i_X \Omega$ where $\Omega$ is
a non-zero global holomorphic $2$-form on $T$. The twisted vector field $X$ can be represented by  a covering of
$\mathcal U=\{ U_i \}$ of $T$ and holomorphic vector fields $X_i \in
T_T(U_i)$ subjected to the conditions
\[
X_i = g_{ij} X_j
\]
on any  non-empty $U_i\cap U_j$,
where $\{ g_{ij} \}$ is a cocycle in $ H^1(\mathcal U,
\mathcal O_T^*)$ representing ${\mathcal N}^{\otimes 2} $.

If $\frac{\partial}{\partial z}$ {and} $ \frac{\partial}{\partial w}$ form a basis of $
H^0( T, T_T)$  then  $X_i = A_i
\frac{\partial}{\partial z} + B_i \frac{\partial}{\partial w}$ for
suitable holomorphic functions $A_i,B_i \in \mathcal O(U_i)$.
Consider the divisor $\Delta$  locally cut out by
\[
\det  \left(
              \begin{array}{cc}
                A_i & B_i \\
                X_i(A_i) & X_i(B_i) \\
              \end{array}
            \right)  .
\]
Clearly these local expressions patch together to form an element  of $
H^0(T, {\mathcal N}^{\otimes 6})$.

{Any  divisor   corresponding to a  fiber of $F$  is defined
by the vanishing of a non-zero element of} $
H^0(T, F^* \mathcal O_{\mathbb P^1}(1))=H^0(T, {\mathcal N})$. By the very definition of $X$, (the closures of) its 1-dimensional orbits are irreducible components of fibers of $F$.
Outside the zero locus of ${X_i}$, the divisor $\Delta_{|U_i}$ corresponds to
the inflection points of the orbits of  ${X_i}$. Indeed, if $\gamma: (\mathbb C,0) \to U_i$ is
an orbit of $X_i$, {that is if $X_i(\gamma(t))=\gamma'(t)$ for $t\in (\mathbb C,0)
$}, then (with an obvious abuse of notation)
\[
  \det  \left(
              \begin{array}{cc}
                A_i & B_i \\
                X_i(A_i) & X_i(B_i) \\
              \end{array}
            \right) ( \gamma )  \equiv {\gamma'} \wedge {\gamma''} .
\]
{Let $L$ be a  linear irreducible components of a fiber of $F$. Its  generic point belongs to $\Delta$ since it is an inflection point of $L$ relatively to $X$ (see \cite[\S6]{jvpfourier}). It follows that $L\leq \Delta$. }

{From the preceding discussion, it follows  that  to prove the theorem } it suffices
to show that {any effective divisor $D_1,\ldots,D_k$ corresponding to a linear fiber of $F$} is
smaller than $\Delta$. Indeed,  the support of distinct fibers of $F$ do not share irreducible
components in common and consequently
\[
  \sum_{i=1}^k D_i \le \Delta  \, .
\]
Since $\sum_{i=1}^k D_i$ is defined by the vanishing of an element in $ H^0 (T,\mathcal N ^{\otimes k})$ while
$\Delta$ is defined by an element in $ H^0 (T,\mathcal N ^{\otimes 6})$, it {would follow  that $k \le 6$ as wanted.} It remains to show that  $D_i \le \Delta$ {for any $i=1,\ldots,k$}.

The divisorial
components of the zero locus of $X_i$ correspond to multiple
components of the fibers of $F$ just like in Darboux's Lemma \ref{L:Darboux}.
If there is a fiber of $F$ containing an irreducible component with multiplicity
$a\ge 2$ and locally cut out over $U_i$ by a reduced holomorphic function $f$ then  we can write $X_i = f^{a-1}
\widetilde{X_i}$ with $\widetilde{X_i}=\widetilde{A_i}
\frac{\partial}{\partial z} + \widetilde{B_i} \frac{\partial}{\partial w}$ holomorphic. Therefore $\Delta$ is locally defined by
\[
 \det  \left(
              \begin{array}{cc}
                f^{a-1} \widetilde{A_i} & f^{a-1} \widetilde{B_i} \\
                f^{a-1}\widetilde{X_i}(f^{a-1}\widetilde{A_i}) & f^{a-1}\widetilde{X_i}(f^{a-1}\widetilde{B_i}) \\
              \end{array}
            \right) = f^{3a - 3 }\det  \left(
              \begin{array}{cc}
                 \widetilde{A_i} & \widetilde{B_i} \\
                \widetilde{X_i}(\widetilde{A_i}) & \widetilde{X_i}(\widetilde{B_i}) \\
              \end{array} \, \right) \, .
\]
Since $3a-3 > a$ when $a\ge 2$ it follows that every  linear fiber of $F$ is smaller than  $\Delta$ as
wanted. Moreover  if $k=6$ then
$F$ cannot have non-reduced fibers.
\end{proof}

Theorem \ref{T:bound} combined with Proposition \ref{P:iso} yields the following corollary.

\begin{cor}\label{C:resumo}
Let $\mathcal W= \mathcal L_1 \boxtimes \cdots \boxtimes \mathcal L_k$ be a linear $k$-web on  $T$.
If $\mathcal F$ is a non-linear foliation  on $T$ such that  $K(\mathcal W \boxtimes \mathcal
F)=0$ then $T$ is isogenous to the square of an elliptic curve and  $k\le 6$.
\end{cor}

\subsection{Constraints on the {l}inear {w}eb}
Let $\mathcal W \boxtimes \mathcal F$ be a {flat} CDQL $(k+1)$-web on a complex torus $T$. If $\mathcal P_{\mathcal F}$ {denotes the polar map of ${\mathcal F}$ and if $\mathcal W= \mathcal L_1 \boxtimes \cdots \boxtimes \mathcal L_k$ then the fibers
$\mathcal P_{\mathcal F}^{-1}(\mathcal L_i)$ are all linear and supported on a union of elliptic curves invariant by  $\mathcal L_i$ or by $\widehat{\mathcal L_i}$ according to Proposition \ref{P:iso}}.
From the very definition of $\mathcal P_{\mathcal F}$ it is clear that the singular set of $\mathcal F$ coincides with the
indeterminacy set~of~$\mathcal P_{\mathcal F}$.

In order to determine the linear web $\mathcal W$ under the assumption that $\mathcal W \boxtimes \mathcal F$ has maximal rank we will
take a closer look at the singularities of $\mathcal F$. It will be convenient  to consider the natural affine coordinates $(x,y)$ {on} the  universal
covering $\mathbb C^2 \to T$.

\begin{lemma}\label{L:alternative}
Let $\mathcal W= \mathcal L_1 \boxtimes \cdots \boxtimes \mathcal L_k$ be a linear $k$-web, with $k \ge 3$, and  $\mathcal F$ be  a non-linear
foliation, {both defined}  on $T$. Suppose that  $K(\mathcal W \boxtimes \mathcal F)=0$. If  $p \in \mathrm{sing}(\mathcal F)$  {is  the origin in the affine coordinate system $(x,y)$}   then one of the following two
alternatives holds:
\begin{enumerate}
\item the foliation $\mathcal F$ is locally given  by $[xdy - ydx + h.o.t.]$. In this case, for each  $i= 1, \ldots, k$, the divisor $\mathrm{tang}(\mathcal F, \mathcal L_i)$ has multiplicity one at $p$ and there exists an elliptic curve through ${p}$ invariant by $\mathcal L_i$ and {by} $\mathcal F$.
\item the foliation $\mathcal F$ is locally given by $[\omega_d + h.o.t.]$ where $\omega_d$ {is a non-zero} homogeneous $1$-form of degree $d \ge 1$ in the coordinates $x,y$  with singular set reduced to $(0,0)$ and not proportional to $xdy-ydx$. In particular the foliation $[\omega_d]$ is non-linear.
\end{enumerate}
\end{lemma}
\begin{proof}
{According to the proof of Proposition \ref{P:abel}, one can assume that}
 $\mathcal L_1= [dx]$,  $\mathcal L_2=[dy]$ and  $\mathcal L_3=[dx - dy]$. If $\mathcal F$ is locally given
by $[a(x,y)dx + b(x,y) dy]$ {where $a$ and $b$ designate holomorphic functions  without common factors, then $\mathrm{tang}(\mathcal F,\mathcal L_1)= \{ b=0\}$,  $\mathrm{tang}(\mathcal F,\mathcal L_2)= \{ a=0\}$ and  $\mathrm{tang}(\mathcal F,\mathcal L_3)= \{ a + b=0\}$. Notice that the assumption $p \in \mathrm{sing(\mathcal F)}$ implies that $a(0,0)=b(0,0)=0$.}

Recall from Proposition \ref{P:iso} that  $\mathrm{tang}(\mathcal F, \mathcal L_1)$ is supported on a  union of elliptic {curves}. Therefore the first non-zero jet of $b$ will be a constant multiple of $x^k \cdot h(x,y)^l$, $k, l \in \mathbb N$, where $h$ is a linear form vanishing defining the tangent space of $\widehat{s \mathcal L_1}$ at zero.  Similarly for ${a}$ and $a+b$.

The first non-zero jet of $a,b$ and $a+b$  have the same degree  and are pairwise without common factor. Otherwise  the supports of $\mathrm{tang}(\mathcal F, \mathcal L_i)$ and $\mathrm{tang}(\mathcal F, \mathcal L_j)$
would  share an  irreducible component in common for some pair $(i,j)$ satisfying $1 \le i < j \le 3$. But this is impossible since $\mathrm{tang}(\mathcal L_i, \mathcal L_j)$ is empty as soon as $i\neq j$.

Thus we can write $\omega =  \omega_d + h.o.t.$  where $\omega_d$ is homogeneous and with singular set equal to the origin.
We are in the first case of the statement when $\omega_d$ is proportional to $xdy-ydx$ and in the second case otherwise.
\end{proof}

We are now in position to use Corollary \ref{C:local} to restrict the possibilities of the maximal linear subweb of an exceptional CDQL web on complex tori.

\begin{prop}\label{P:lista}
Let $\mathcal W= \mathcal L_1 \boxtimes \cdots \boxtimes \mathcal L_k$ be a linear $k$-web, with $k \ge 4$, and  $\mathcal F$ be  a non-linear
foliation  on $T$. If $\mathcal W \boxtimes \mathcal F$ has maximal rank then, up to isogenies, one of the following alternatives holds:
\begin{enumerate}
\item The torus $T$ is the square of an elliptic curve, $k=4$ and $ \mathcal W = [ dx dy (dx^2-dy^2) ];$
\item The torus $T$ is  $E_{i}^2$, $k=6$ and $ \mathcal W = [ dx dy (dx^2 - dy^2) (dx^2 + dy^2)];$
\item The torus $T$ is  $E_{\xi_3}^2$, $k=5$ and $ \mathcal W = [ dx dy (dx^3  + dy^3) ];$
\item The torus $T$ is  $E_{\xi_3}^2$, $k=4$ and $ \mathcal W = [ dx dy (dx + dy) (dx - \xi_3 \, dy)].$
\end{enumerate}
\end{prop}
\begin{proof}
Corollary \ref{C:resumo} tell us that  $T$ is isogeneous to the square of
an elliptic curve $E$ and that $k \le 6$.   Lemma \ref{L:mult} implies that we can assume, after an eventual reordering,  that $\widehat{\mathcal L_1} \neq \widehat{\mathcal L_2}$.
For $i=1,2$, let $\widehat{E_i}$  be an elliptic curve contained in $\mathrm{tang}(\mathcal F, \mathcal L_i)$ that
is $\widehat{\mathcal L_i}$-invariant. Notice that the existence of these curves is ensured by Lemma \ref{L:abelian}.

Since $\widehat{\mathcal L_1} \neq \widehat{\mathcal L_2}$ there exists $p \in \widehat{E_1} \cap \widehat{E_2}$. Notice that $p$ belongs
to $\mathrm{sing}(\mathcal F)$.   Moreover our choice  of $p$  implies  that  it fits in  the second alternative  of Lemma \ref{L:alternative}.
Therefore we can apply Corollary \ref{C:local} to conclude that the first non-zero jet of $\mathcal W\boxtimes \mathcal L$ at $p$ is
equivalent, under a linear change of the affine  coordinates $(x,y)$, to one of the following webs:
\begin{equation}\label{E:list}
\mathcal A_I^4,\mathcal A_I^5, \mathcal A_I^6 , \,  \mathcal A_{III}^{2}, \mathcal A_{III}^{3}, \mathcal A_{III}^{4}, \,  \mathcal A_5^d ,    \mathcal A_6^b.
\end{equation}

To prove the proposition we will analyze the constraints imposed on the torus $T$ by the above {\it local models}.

\smallskip

Notice that the $5$-web  $\mathcal A_I^4 = [ (dx^4 - dy^4)] \boxtimes [d(xy) ] $ is isomorphic (via a linear map) to $[ dx dy (dx ^2 - dy^2)] \boxtimes [ d (x^2 + y^2) ]$. All the defining foliations of $[ dx dy (dx ^2 - dy^2)]$ are tangent to elliptic fibrations on the square of an arbitrary elliptic curve $E$. Therefore these local models
do not impose restrictions on the curve $E$.  Similarly the $5$-web $\mathcal A_{III}^2=[dxdy(dx^2-dy^2)]\boxtimes[d(xy)]$ also does not impose restrictions on $E$. Indeed these two local models coexist in distinct singular points of the exceptional CDQL $5$-webs $\mathcal E_{\tau}$.

\smallskip

The $6$-webs $\mathcal A_{III}^{3}=[dx dy (dx^3- dy^3)]\boxtimes [d(xy)]$ and $ \mathcal A_6^b =  [dx dy (dx^3 + dy^3)]\boxtimes [d(x^3 + y^3)]$ share
the same  linear $5$-web (after the change of coordinates $(x,y) \mapsto (x,-y)$ on $\mathcal A_6^b$). On the one hand Proposition \ref{P:iso} implies that all
the defining foliations of the linear $5$-web $[dx dy (dx^3- dy^3)]$ must be tangent to elliptic fibrations. On the other hand Lemma \ref{L:cross} implies
that $\xi_3 \in \mathrm{End}(E) \otimes \mathbb Q$.
Therefore   $T$ must be isogenous to $E_{\xi_3}^2$. Notice that both  local models  coexist in distinct singular points of the exceptional CDQL $6$-web $\mathcal E_6$.

\smallskip

The same argument shows that  the $5$-web $\mathcal A_5^d = [ dxdy (dx + dy) (dx - \xi_3 dy)] \boxtimes [d(xy(x+y)(x- \xi_3 y))]$ can only be a local model for an exceptional CDQL web when  $T$ is isogenous to $E_{\xi_3}^2$.  Similarly
the $7$-web $\mathcal A_{III}^4 = [ dx dy (dx^4 - dy^4)] \boxtimes [d(xy)]$ can only be a local model for an
exceptional CDQL web when $T$ is isogeneous to $E_i^2$.

\medskip

To conclude the proof of the Proposition it suffices to show that the two remaining possibilities  in the list (\ref{E:list}) (namely $\mathcal A_I ^5$ and $\mathcal A_I ^6$) cannot appear as  local models for  exceptional CDQL webs on a torus.

\smallskip

We will first deal with  the $6$-web $\mathcal A_I^5=[ (dx^5 - dy^5)] \boxtimes [d(xy) ]$.
If $\xi_5$ is {a} primitive $5$th root of the unity then the  cross-ratio $(1,\xi_5:\xi_5^2,\xi_5^3)$ is a  root of the polynomial $p(x)=x^2-x-1$. Notice that  the roots of $p(x)$ are the golden-ratio and its  conjugate: $1/2 \pm \sqrt{5}/2$. In particular they are  irrational real numbers and, as such, cannot induce an endomorphism on any elliptic curve $E$. Lemma \ref{L:cross} implies that does not exist a two-dimensional complex torus $T$  where  all the defining foliations
of  $[ (dx^5 - dy^5)]$ are tangent to elliptic fibrations. Proposition \ref{P:iso} implies that $\mathcal A_I^5$ cannot appear as  a local model {of}  an exceptional CDQL web on a torus.

\smallskip

We also claim that the  $7$-web $\mathcal A_I^6=[ (dx^6 - dy^6)] \boxtimes [d(xy) ] $  cannot appear as a local model of an exceptional CDQL web on a torus $T$. Using Lemma \ref{L:cross} it is a simple matter to show that $T$ is isogenous to $E_{\xi_3}^2$. Assume now that $\mathcal L_1$ and $\mathcal L_2$ are such that $\mathcal L_1 \neq \widehat{\mathcal L_2}$. Lemma \ref{L:abelian} ensures that there are: an elliptic curve $E_1$ in $\mathrm{tang}(\mathcal F, \mathcal L_1)$ invariant by $\mathcal L_1$ and an elliptic curve $\widehat{E_2}$ in $\mathrm{tang}(\mathcal F, \mathcal L_2)$  invariant by $\widehat{\mathcal L_2}$. Since $\mathcal L_1 \neq \widehat{\mathcal L_2}$ there exits $ p \in  E_1 \cap \widehat{E_2}$. Since $p \in |\mathrm{tang}(\mathcal F, \mathcal L_1)| \cap |\mathrm{tang}(\mathcal F, \mathcal L_{2})|$,  it  is a singular point of  $\mathcal F$.

Notice that our choice of $p$  implies  that the first non-zero jet of  $\mathcal F$ at $p$ is non-linear, see  Lemma \ref{L:alternative}.
Since $E_1$ is also $\mathcal F$-invariant, the linear polynomial defining it on the affine coordinates $(x,y)$ will be also invariant by the first
jet of $\mathcal F$. But for the $7$-web $\mathcal A_I^6$ none of the invariant lines through $0$ of the non-linear  foliation is invariant by
 any of the linear  foliations. Therefore  the local model at $p$ must be the only other $7$-web appearing in the list (\ref{E:list}): $\mathcal A_{III}^4 = [ dx dy (dx^4 - dy^4)]\boxtimes[d(xy)]$. But this implies that $T$ is isogenous to $E_i^2$. Since $E_i^2$ is not isogenous to $E_{\xi_3}^2$ the claim follows and so does the proposition.
\end{proof}

\subsection{The classification of exceptional   CDQL webs on tori}

To obtain the classification of exceptional CDQL webs on tori we will analyze in Sections \ref{S:continuous}, \ref{S:77}, \ref{S:66} and \ref{S:55}
the respective alternatives {\it(1)},{\it(2)},{\it(3)} and {\it(4)}  provided by Proposition \ref{P:lista}.
\subsubsection{{\bf The {c}ontinuous {f}amily of exceptional CDQL $5$-webs}}\label{S:continuous}

In case {\it(1)} of Proposition \ref{P:lista}, the torus  $T$ is isogenous to the square of an elliptic curve, $k=4$ and the linear web is  $\mathcal W=[ dx dy (dx^2-dy^2) ]$.
As we have proved  in  Example \ref{E:PT}  every flat (in particular exceptional) CDQL $5$-web of the form $\mathcal W \boxtimes \mathcal F$ must be isogenous to one of the  $5$-webs $\mathcal E_{\tau}$ {(with $\tau \in \mathbb H$)} presented in the Introduction.

\subsubsection{{\bf The  {e}xceptional CDQL $7$-web on $E_i^2$}}\label{S:77}

In the second alternative of Proposition \ref{P:lista}, the torus $T$ is isogenous to  $E_{i}^2$, $k=6$ and the linear web
is $\mathcal W=\mathcal W_1 \boxtimes \mathcal W_2$ where $\mathcal W_1= [ dx dy (dx^2 - dy^2)]$ and $\mathcal W_2 =  [dx^2 + dy^2].$
This decomposition of $\mathcal W$ satisfies the hypothesis of Corollary \ref{C:w-f}. Therefore a non-linear foliation $\mathcal F$  satisfies   $K(\mathcal F \boxtimes \mathcal W)=0$
if and only if $K(\mathcal F \boxtimes \mathcal W_1)=K(\mathcal F \boxtimes \mathcal W_2)=0$. Thus the subweb $\mathcal F \boxtimes \mathcal W_1$ is isogenous to a web of the continuous family $\mathcal E_{\tau}$. We loose no generality by assuming that $\mathcal F \boxtimes \mathcal W_1 = \mathcal E_{\tau}$ for some $\tau \in \mathbb H$. It remains to determine $\tau$. Since $T$ is isogenous  to $E_i^2$ we know that
$\tau = \alpha + \beta i$ for suitable rational numbers $\alpha, \beta$.  Set $\Gamma = \mathbb Z \oplus (\alpha + \beta i) \mathbb Z$.

Recall from Section \ref{S:ell57} that the non-linear foliation ${\mathcal F}$ is equal to $[dF_{\tau}]$ where
\[
  F_{\tau}(x,y) = \left(\frac{\vartheta_1(x,\tau)  \vartheta_1(y,\tau) }{\vartheta_4(x,\tau)  \vartheta_4(y,\tau) }\right)^2   .
\]
Recall also  that $\mathrm{Indet}(F_{\tau}) = \big\{ (\tau/2,0) , (0,\tau/2) \big\}$ and that these indeterminacy points correspond to
radial singularities of $\mathcal F$.

The tangency of $\mathcal F$ with the linear foliation $[dx + i dy]$ at $(0,\tau/2)$ has first non-zero jet equal to $(x + i \,y)$
since ${(x\,dy - y\,dx) \wedge (dx + i\, dy) = - (x + i \,y) \,dx \wedge dy}$. Therefore, Theorem \ref{T:curvaturatoro} implies that there exists
an elliptic curve $C$ through $(0,\tau/2)$  invariant by $\mathcal F$ and by $[dx + i dy]$. Notice that $C$ is the image of the entire map
\begin{eqnarray*}
\varphi : \mathbb C &\longrightarrow& E_{\tau}^2 = \left({\mathbb C}/{\Gamma}\right) ^2 \\
z &\longmapsto& \!\!(-i\, z ,  z + \tau/2 )   \, .
\end{eqnarray*}
Thus $C \cap E_{0,1} = \varphi(i \Gamma )$. {The curve $E_{0,1}$ is  also $\mathcal F$ invariant (but do not coincide with $C$) so the set   $C\cap E_{0,1}$ is  contained in $\mathrm{sing}(\mathcal F)$.  But the singularities of $\mathcal F$ over $E_{0,1}$ are  $(0,0)$ and $(0,\tau/2)$. Moreover the singularity at $(0,0)$  has only two separatrices, namely $E_{1,0}$ and $E_{0,1}$. It follows that $C \cap E_{0,1} = \varphi(i \Gamma )$ is equal to the radial singularity $(0,\tau/2)$ of $\mathcal F$ on $E_{0,1}$.
Therefore $i \, \Gamma + \tau/2 \subset \Gamma + \tau/2$. Consequently   $i\Gamma \subset \Gamma$
and  $-\Gamma \subset i\,\Gamma$. Thus $i\,\Gamma=\Gamma$, that is $i \in \mathrm{Aut}(E_{\tau})$.} This is sufficient to show that the elliptic curve
$E_{\tau}$ is isomorphic to $E_i$.

Recall that
\[
\Gamma_0(2) = \left\{ \left(
                                     \begin{array}{cc}
                                       a & b \\
                                       c & b \\
                                     \end{array}
                                   \right) \in \mathrm{SL}(2,\mathbb Z) \, \Big\vert \, \,  b \equiv 0 \mod 2 \right\}.
\]
Thus, modulo the action of $\Gamma_0(2)$ we can assume that $\tau \in \{ i, 1+i, (1+i)/2 \}$. Moreover  the $\mathbb Z_2$-extension of $\Gamma(2)$ by the transformation $z \mapsto -2/z$ identifies $1+i$ with $(1+i)/2$ because $-2((1+i)/2)^{-1}=-2+2i$. Therefore we can assume that $\mathcal F \boxtimes \mathcal W_1$ is equal to
$\mathcal E_{1+i}$ or to $\mathcal E_{i}$. If $\tau =i$ then
$(i/2,0)$ is a radial singularity of $\mathcal F$ and, as above, the
 curve $L_{(i/2,0)} E_{1,i}$ invariant by $[dx+idy]$ is also $\mathcal F$-invariant. But this curve  intersects the $\mathcal F$-invariant curve
   $E_{0,1}$ at $(0,1/2)$ which is not
a singularity of $\mathcal F$. This contradiction implies that, up to isogenies, $\mathcal E_7 = [dx^2 + dy^2]\boxtimes \mathcal E_{1+i}$ is the
unique exceptional CDQL $7$-web on complex tori.

\subsubsection{{\bf The  exceptional CDQL $6$-web on $E_{\xi_3}^2$}  }\label{S:66}
In the third alternative of Proposition \ref{P:lista} the torus $T$ is  isogenous to $E_{\xi_3}^2$, $k=5$ and  the linear web is $\mathcal W=\mathcal W_1 \boxtimes \mathcal W_2$ with $\mathcal W_1 = [ dx dy]$ and $\mathcal W_2 = [(dx^3  + dy^3) ]$.
As in the previous case this decomposition satisfies the hypothesis of Corollary \ref{C:w-f}. Therefore
 $\mathcal F$ is a non-linear foliation on $T$ satisfying  $K(\mathcal W \boxtimes \mathcal F)=0$ if and only if $K(\mathcal W_1 \boxtimes \mathcal F)=K(\mathcal W_2 \boxtimes \mathcal F)=0$.

If $K(\mathcal F \boxtimes [dx dy]) = 0$ then Theorem \ref{T:curvaturatoro} (see also \cite{PT})  implies that
$\mathcal F = [ a(x) dx + b(y) dy ]$  for suitable rational functions $a, b \in \mathbb C({E_{\xi_3}})$.
Moreover, {according to item {\it (3)} of Corollary \ref{C:local},}  we can assume
that the singularity of $\mathcal W \boxtimes \mathcal F$ at $(0,0)$ has  first non-zero jet  equivalent to
$\mathcal A_6^b=[dx dy (dx^3 + dy^3)]\boxtimes [x^2dx + y^2dy]$. In particular,  {interpreting} $x,y$ as coordinates
on the universal covering of $T$, we can assume that the meromorphic functions $a$,$b$ satisfy
$a(x) = x^2 +O(x^3)$ and $b(y) = y^2 + O(y^3)$. In particular, $a(0)=a'(0)=b(0)=b'(0)=0$.

A tedious (but  trivial) computation shows that $K(\mathcal F \boxtimes [dx^3+ dy^3])$ is equal to
\[
6 \frac{ba^3\left(aa'' - 2a'^2\right) -ab^3\left(bb'' - 2b'^2\right)
+a^4\left( b'^2 + b b ''\right) - b^4\left( a'^2 + a a'' \right)
  }{(a^3-b^3)^2}\, dx \wedge dy.
\]
After deriving twice the numerator with respect to $y$, one obtains
\[
a^3b''\left( 2(a')^2 - a\,a'' -3\,b''a \right) - 4\,a^4{b'}b''' + b R
\]
where $R$ is a polynomial in $a(x), b(y)$ and theirs derivatives up to order four. Evaluation of this expression at $y=0$ yields the following
second order differential equation {identically satisfied by $a$}:
\begin{equation}\label{E:2order}
a^3\big( 2\,(a')^2 - aa'' -6 \,a \big){=0}  .
\end{equation}

\begin{lemma}
If   $a :(\mathbb C, 0) \to \mathbb C$ is  a germ of  solution of (\ref{E:2order}) satisfying the boundary conditions
$a(0)=a'(0)=0$ and $a''(0)=2$ then
\[ \qquad \qquad \qquad
a(x) = x^2   \qquad  \text{ or }  \qquad a(x) = \frac{\lambda^2}{\wp(\lambda^{-1} x,\xi_3)}
\]
for a suitable $\lambda \in \mathbb C^*$.
\end{lemma}
\begin{proof}
Notice that the $6$-web $[dxdy(dx^3 + dy^3)]\boxtimes[a(x) dx + a(y) dy]$ with
$a(x)=x^2$  is the $6$-web  $\mathcal A_6^b$ from the introduction. Similarly when
$a(x) = {{\lambda^2}/\wp(\lambda^{-1} x,\xi_3)}$ then the $6$-web $[dxdy(dx^3 + dy^3)]\boxtimes[a(x) dx + a(y) dy]$
can be obtained from $\mathcal E_6$ by the change of coordinates $(x,y)\mapsto (\lambda x,\lambda y)$. Since both
$\mathcal E_6$ and $\mathcal A_6^b$ are exceptional, the corresponding $a$'s  are  solutions of (\ref{E:2order}).
Clearly they all satisfy the boundary conditions. To prove the lemma it suffices to verify that they are the only solutions.

If $a(x)$ is a solution of (\ref{E:2order}) satisfying the boundary conditions then it  is indeed a solution of
$2(a')^2 - aa'' -6 a =0$. Therefore $\gamma(t) = (a(t), a'(t))$ is an orbit of the following vector field
\[
Z(x,y) =  y \frac{\partial}{\partial x} +  \frac{2y^2 - 6 x}{x} \frac{\partial}{\partial y} \, ,
\]
that is $Z(\gamma(t))=\gamma'(t)$.

Notice that $Z$ admits as a rational first integral the function $\frac{y^2 - 4x }{x^4}$. Therefore  every solution $a(x)$ of (\ref{E:2order}) satisfying
$a(0)=a'(0)=0$ and $a''(0)=2$ must parameterize (through the  map $t \mapsto (a(t), a'(t))$) {a branch of}
one of the curves $y^2 - 4x + \mu x^4$ for some $\mu \in \mathbb C$. When $\mu=0$, the corresponding curve is parameterized by $a(x)=x^2$. For
$\mu \neq 0$ it is parameterized by $a(x) = {\lambda^2}/{\wp(\lambda^{-1} x,\xi_3)}$ with $\lambda$ satisfying $\mu\lambda^6=1$.  Notice that the different choices
for $\lambda$ leads to the same function $a$. Indeed, the symmetry $-\xi_3 ( \mathbb Z \oplus \xi_3 \mathbb Z ) = \mathbb Z \oplus \xi_3 \mathbb Z$ combined with (\ref{E:homogeneo}) implies that
\begin{equation}\label{E:kkk}
  \frac{(-\xi_3)^2}{\wp\big((-\xi_3)^{-1} x, \xi_3\big)} = \frac{ 1 } {\wp( x , \xi_3)}.
\end{equation}
Since   each of the curves $\{ y^2 - 4x - \mu x^4=0\}$ admits only one parametrization of the  form
$t\mapsto (a(t), a'(t))$ with $a''(0)=2$, the lemma
follows.
\end{proof}

Keeping in mind that the coefficients of the defining $1$-form of $\mathcal F$ must be
doubly-periodic functions and the symmetry of our setup, so far we have proved that $K(\mathcal F \boxtimes \mathcal W)=0$
implies that, up to homotethies,
\[
  \mathcal F= \left[ \frac{dx}{\wp( x,\xi_3)}   + \frac{\lambda^2 dy }{\wp(\lambda^{-1} y,\xi_3)} \right]
\]
for a suitable $\lambda \in \mathbb C^*$. Computing again $K(\mathcal F \boxtimes [dx^3+ dy^3])$ yields
\[
 \frac{ 9 \lambda^2(\lambda^6 -1)\wp(x,\xi_3)^2 \wp(y/\lambda,\xi_3)^2  }{ \lambda^{12} \wp(x, \xi_3)^6-2\lambda^6\wp(y/\lambda, \xi_3)^3\wp(x, \xi_3)^3+\wp(y/\lambda, \xi_3)^6}{\, dx \wedge dy.}
\]

The vanishing of the curvature, taking into account (\ref{E:kkk}), implies that
\[
  \mathcal F= \left[ \frac{dx}{\wp( x,\xi_3)}   + \frac{ dy }{\wp(y,\xi_3)} \right]  .
\]
It follows that the $6$-web $\mathcal F \boxtimes \mathcal W$  is isogenous to the $6$-web $\mathcal E_6$
from the Introduction.

\subsubsection{{\bf The exceptional CDQL $5$-web on $E_{\xi^3}^2$. Combinatorial patchwork}}\label{S:55}

In the last case of Proposition \ref{P:lista} (transformed via the change of coordinates $(x,y)\mapsto(y,-x)$), the complex torus $T$ is isogenous to   $E_{\xi_3}^2$, $k=4$ and the linear web $\mathcal W$  is ${ [ dx dy (dx -dy) (\xi_3\,dx + dy)]}$. Unlikely in the previous case the web $\mathcal W$ does not admit a decomposition satisfying the hypothesis of Corollary \ref{C:w-f}. We have not
succeeded in dealing with this case using  analytic methods as  in the previous section and  in \cite{PT}. We were lead to adopt a geometrical/combinatorial approach.

If  $\mathcal L_1 = [dx]$ ,  $\mathcal L_2 = [dy]$, $\mathcal L_3 = [dx-dy]$ and  $\mathcal L_4 = [\xi_3\,dx +  dy ]$ then
  straight-forward computations using formula (\ref{E:adef}) show  that
\begin{align}
\label{E:Lhat}
\widehat{\mathcal L_1}   = &  \big[ dx + ( \xi_3^2-1)\, dy \big]  & &   \hspace{-1cm} \widehat{\mathcal L_3}   =   {\big[dx -\xi_3\,  dy \big] }       \\
 \widehat{\mathcal L_2}   =  & \big[(\xi_3 -1)\, dx + dy    \big]  & & \hspace{-1cm} \widehat{\mathcal L_4}   =  { \big[ dx + \xi_3\, dy  \big].} \nonumber
\end{align}
For $i=1,2,3,4$, the leaves of $\widehat{\mathcal L_i}$ are translates of the elliptic curve $\widehat{E_i}$ where
\[
\widehat{E_1}   =  E_{1-\xi_3^2,1} , \quad
\widehat{E_2}   =  E_{1,1-\xi_{3}} , \quad
\widehat{E_3}   =  E_{{\xi_3,1}}  \quad  \text{and} \quad
\widehat{E_4}   =  E_{{\xi_3,-1} } \, .
\]

Suppose that $\mathcal F$ is a non-linear foliation on $T$ such that $\mathcal W \boxtimes \mathcal F$ has maximal rank.
According to Corollary \ref{C:local} {and taking into account the change of coordinates $(x,y)\mapsto(y,-x)$, there are only two possibilities for a singularity $p$ of $\mathcal F$: either
$p$ is a radial singularity or the first non-zero jet of  $\mathcal F \boxtimes \mathcal W$ at $p$  is equivalent to  }
\[
\big[  dx\, dy\, (dx-dy)\,(\xi_3\,dx +  dy) \big] \boxtimes
\big[d\big(xy(x-y)(\xi_3\,x+y) \big)  \big]  .
\]
{ We will say that the former singularities are of type $A$  whereas the latter are of type $B$. We will
write $\mathrm{sing}^A(\mathcal F)$ (resp.   $\mathrm{sing}^B(\mathcal F)$) for the set of singularities of type $A$ (resp. of type
$B$).}

{By the very definition, the first non-zero jet of $\mathcal F$ at a point $p \in \mathrm{sing}^B(\mathcal F)$ is}
$$ \mathcal F_0 = \Big[d\big(xy(x-y)(\xi_3\,x +y) \big)  \Big].$$
Simple computations~show~that
\begin{equation}\label{E:ttt}
\mathrm{tang}(\mathcal F_0, \mathcal L_i) = \left\lbrace
                                              \begin{array}{lcc}
                                                \big\{\, x \left( x+ (\xi_3^2-1)y\right)^2 = 0\,\big\} & \text{when}  & i = 1 \vspace{0.08cm}\\
                                                \big\{\, y \left((\xi_3-1)x + y\right)^2 = 0\, \big\} & \text{when} & i=2 \vspace{0.08cm}\\
                                               { \big\{\, (x-y) \left(x -\xi_3\, y\right) ^2  = 0\, \big\}} & \text{when} & i=3\vspace{0.08cm} \\
                                               { \big\{\, (\xi_3\,x+  y) \left(  x + \xi_3\,y \right)^2 = 0\,\big\}} & \text{when} & \; \; i=4\; . \\
                                              \end{array}
                                            \right.
\end{equation}

Being aware of the first nonzero jets of the singularities of $\mathcal F$, we are able to
describe the first non-zero jets of $\mathrm{tang}(\mathcal F, \mathcal L_i)$. This is  the
content of the two following lemmata.

{\begin{lemma}\label{L:csing}
Let  $p \in \mathrm{sing}^A (\mathcal F)$.
For every $i \in \{ 1, \ldots, 4\}$, there is an unique irreducible component of the divisor $\mathrm{tang}(\mathcal F, \mathcal L_i)$ passing through $p$: it is an irreducible curve $C$ invariant by  $\mathcal L_i$. In particular, there is no  $\widehat{\mathcal L_i}$-invariant curve in  $\mathrm{tang}(\mathcal F, \mathcal L_i)$ passing through  $p$.
\end{lemma}}
\begin{proof}
Since the first non-zero jet  of $\mathrm{tang}(\mathcal F_, \mathcal L_i)$  at $p$  coincides with $\mathrm{tang}([xdy-ydx], \mathcal L_i)$
the lemma follows from  Theorem \ref{TC:curvatura}.
\end{proof}

\begin{lemma}\label{L:csingB}
Let  $p \in \mathrm{sing}^B (\mathcal F)$.
For every $i \in \{ 1, \ldots, 4\}$, the divisor $\mathrm{tang}(\mathcal F, \mathcal L_i)$
contains in its support two distinct irreducible curves {$C_i$ and $\widehat{C}_i$  both containing $p$. Moreover,  $C_i$ (resp. $\widehat{C}_i$) is invariant by $\mathcal L_i$ (resp. by
 $\widehat{\mathcal L_i}$).}
\end{lemma}
\begin{proof}
Since the first non-zero jet  of $\mathrm{tang}(\mathcal F_, \mathcal L_i)$  at $p$ coincides with $\mathrm{tang}(\mathcal F_0, \mathcal L_i)$ the lemma follows
from {Theorem \ref{TC:curvatura} combined with (\ref{E:Lhat}) and (\ref{E:ttt}).}
\end{proof}

The core  of our  argument to characterize $\mathcal E_5$ is contained in the next lemma.

\begin{lemma}\label{L:primol}
Let $\mathcal F$ be a non-linear foliation on the torus $T=E_{\xi_3}^2$.
Suppose that  the 5-web $\mathcal F \boxtimes [ dx dy (dx - dy) {(\xi_3\,dx + d y)}]$ has
maximal rank. If  $0\in \mathrm{sing}^B(\mathcal F)$  then
\begin{enumerate}
\item[(a)] $(0,y) \in \mathrm{sing}(\mathcal F)$ if and only if $(y,0) \in \mathrm{sing}(\mathcal F)$;\vspace{0.08cm}
\item[(b)] If $(y,0) \in \mathrm{sing}(\mathcal F)$ then $(2y,0) \in \mathrm{sing}^B(\mathcal F)$;\vspace{0.08cm}
\item[(c)] If $(y,0) \in \mathrm{sing}(\mathcal F)$ then $(-\xi_3^{{2}}\, y,0) \in \mathrm{sing}(\mathcal F)$;\vspace{0.08cm}
\item[(d)] Both $\mathrm{sing}(\mathcal F) \cap E_{1,0}$ and  $\mathrm{sing}^B (\mathcal F) \cap E_{1,0}$
are subgroups of $E_{1,0}$.
Similarly $\mathrm{sing}(\mathcal F) \cap E_{0,1}$ and  $\mathrm{sing}^B (\mathcal F) \cap E_{0,1}$ are subgroups of $E_{0,1}$.
\end{enumerate}
\end{lemma}
\begin{proof}To prove {\it (a)}, {we will use  that the  curves  $E_{0,1}$, $E_{1,0}$ and $E_{1,1}$ passing through $(0,0)$ are $\mathcal F$-invariant (what is ensured by
Lemma \ref{L:primol}). If $(0,y) \in \mathrm{sing}(\mathcal F)$ then Lemma \ref{L:csing} implies that $L_{(0,y)} E_{1,0}$ is $\mathcal F$-invariant.  Therefore $L_{(0,y)} E_{1,0} \cap
E_{1,1} =(y,y)$  is the intersection of two distinct leaves of ${\mathcal F}$. It follows that
$(y,y)\in \mathrm{sing}(\mathcal F)$.
Consequently $L_{(y,y)} E_{0,1}$ is also $\mathcal F$-invariant. Since $(y,0) = L_{(y,y)} E_{0,1} \cap E_{1,0}$, item {\it (a)} follows. }

{To prove {\it (b)}, start by  noticing that $(0,y) \in \mathrm{sing}(\mathcal F)$  by
 {\it (a)}. Therefore $L_{(0,y)} E_{1,0}$ is $\mathcal F$-invariant according to  Lemma \ref{L:csing}.
By hypothesis $(0,0) \in \mathrm{sing}^B(\mathcal F)$ thus Lemma \ref{L:csingB} ensures that the elliptic curve $\widehat{E_1}=E_{1-\xi_3^2,1}$ belongs to $\mathrm{tang}(\mathcal F, \mathcal L_1)$. The curve  $L_{(0,y)} E_{1,0}$ being  invariant by $\mathcal L_2$ and ${\mathcal F}$ (since $(0,y) \in \mathrm{sing}(\mathcal F)$), it is necessarily an irreducible component of ${\rm tang}({\mathcal F},\mathcal L_2)$. As a consequence,  the intersection $ \widehat{E_1}\cap L_{(0,y)} E_{1,0}$ is included in $\mathrm{sing}^B(\mathcal F)$. In particular,   the point $p=((1-\xi_3^2)y,y)$ belongs to $\mathrm{sing}^B(\mathcal F)$. Considering now the $\widehat{ \mathcal L_3}$-invariant curve through $p$, that is
$L_p E_{\xi_3,1}$, we see that it intersects $E_{0,1}$ at $(2y,0)$.
Thus $(2y,0)\in \mathrm{sing}^B(\mathcal F)$   proving item {\it (b)}.}

To prove item {\it (c)},  recall from the previous paragraph that  $L_{(0,y)} E_{1,0}$ is $\mathcal F$-invariant. The curve  $E_{1,-\xi_3}$ intersects $L_{(0,y)} E_{1,0}$ at $p=(-\xi_3^2 \,y , y )$. Since $E_{1,-\xi_3}$ is $\mathcal F$-invariant (by Lemma \ref{L:csing}) it follows $p \in \mathrm{sing}(\mathcal F)$. Consequently $L_p E_{0,1}$ is
$\mathcal F$-invariant (again by Lemma \ref{L:csing}) and $(-\xi_3^2\,y,0) = L_p E_{0,1} \cap E_{1,0} \in \mathrm{sing}(\mathcal F).$ Item {\it (c)} follows.\vspace{0.1cm}

It remains to prove item {\it (d)}. We will first prove that $S= \mathrm{sing}(\mathcal F) \cap E_{1,0}$ is a subgroup of $E_{1,0}$. From item {\it (c)} it follows   $(y,0) \in S$ if and only if $(-y,0) \in S$. Thus it suffices to show that, given  two elements  $(y_1,0)$ and $ (y_2,0)$ of $S$, their sum $(y_1+y_2,0)$ is also in $S$. Item {\it (a)} implies that $(0,y_2) \in \mathrm{sing}(\mathcal F)$ and consequently  the curve $L_{(0,y_2)} E_{1,0}$ is $\mathcal F$-invariant (by Lemma \ref{L:csing}). For the same reason the curve $L_{(y_1,0)} E_{1,1}$ is also $\mathcal F$-invariant thus the point $p=(y_2+y_1,y_2) \in L_{(0,y_2)} E_{1,0} \cap L_{(y_1,0)} E_{1,1}$ belongs to $\mathrm{sing}(\mathcal F)$. Since $L_p E_{0,1}$ intersects $E_{0,1}$ at $(y_1 + y_2,0)$ and because these two curves are $\mathcal F$-invariant, it follows that  $(y_1 + y_2,0) \in S$. Therefore $\mathrm{sing}(\mathcal F) \cap E_{1,0}$ is a subgroup of $E_{1,0}$.

Consider now the group homomorphism
\begin{eqnarray*}
  S &\longrightarrow& S \\
  x &\longmapsto& x + x \, .
\end{eqnarray*}
Item {\it (b)} implies that its image is $\mathrm{sing}^B (\mathcal F)\cap E_{1,0}$. Therefore $\mathrm{sing}^B (\mathcal F)\cap E_{1,0}$ is also a subgroup of $E_{1,0}$.

 Mutatis mutandis we obtain the same statements for $\mathrm{sing} (\mathcal F)\cap E_{0,1}$ and $\mathrm{sing}^B (\mathcal F)\cap E_{0,1}$: both are subgroups of $E_{0,1}$.
\end{proof}

\begin{thm}\label{T:E5}
Let $\mathcal F$ be a non-linear foliation on $T=E_{\xi_3}^2$.
If the 5-web $   [ dx dy (dx - dy) (\xi_3\,dx + dy)] \boxtimes\mathcal F  $ has
maximal rank then it is isogenous to  $\mathcal E_5$.
\end{thm}
\begin{proof}
Let us denotes by $\equiv$ the numerical equivalence of divisors on $T$.
Since $\mathcal O_T(\mathrm{tang}(\mathcal F, \mathcal L_i))=N\mathcal F$ for $i=1,\ldots, 4$, all the divisors $\mathrm{tang}(\mathcal F, \mathcal L_i)$  are
pairwise  linearly equivalent. Moreover, Theorem \ref{T:curvaturatoro} implies that
\[
\mathrm{tang}(\mathcal F, \mathcal L_i) \equiv a_i E_i + b_i \widehat{E_i}
\]
for $i=1,\ldots, 4$, where $E_i$ and $\widehat{E_i}$ are elliptic curves in $T$ invariant by $\mathcal L_i$ and  $\widehat{\mathcal L_i}$ respectively and $a_i,b_i$ are non-negative integers.
Indeed Lemma \ref{L:abelian} implies that $a_i,b_i$ are positive integers.
In particular we obtain that
\[
  a_1 E_{0,1} + b_1 E_{1-\xi_3^2,1} \equiv  a_2 E_{1,0} + b_2 E_{1,1-\xi_{3}} \, .
\]
Intersecting both members with $E_{0,1}$, $E_{1,0}$ and $E_{1,1}$  we obtain respectively
\[
 3 b_1 = a_2 + b_2 \, , \quad a_1 + b_1 = 3 b_2 \, \quad \text{and} \quad a_1+ b_1 = a_2 + b_2 \, .\qquad
\]
Thus $a_1/b_1 = a_2/b_2 = 2$.

Assume, without loss of generality, that $0 \in T$ is point in $\mathrm{sing}^B (\mathcal F)$. Notice
that $E_{1,0}$ is $\mathcal F$-invariant and $\mathrm{sing}(\mathcal F) \cap E_{1,0}$ is equal to the
set of intersection points of $\mathrm{tang}(\mathcal F, \mathcal L_1)$ with $E_{1,0}$. Moreover $\mathrm{sing}^B (\mathcal F)$
corresponds to the intersection  with $E_{1,0}$   of the irreducible components of $\mathrm{tang}(\mathcal F, \mathcal L_1)$ that are invariant by $\widehat{\mathcal L_1}$. Equation (\ref{E:ttt}) implies that  each of the $\widehat{\mathcal L_1}$-invariant curves  in $\mathrm{tang}(\mathcal F, \mathcal L_1)$
appears with multiplicity two. From  $a_1/b_2=2$ it follows that  the cardinality of $\mathrm{sing}(\mathcal F) \cap E_{1,0}$ is four times the cardinality of
$\mathrm{sing}^B(\mathcal F) \cap E_{1,0}$. Recall from Lemma \ref{L:primol} item\;{\it (d)}  that $S=\mathrm{sing}(\mathcal F) \cap E_{1,0}$ and $S^B=\mathrm{sing}^B (\mathcal F) \cap E_{1,0}$ are subgroups of $E_{1,0}$. It is now clear that  the kernel of the map $S \to S^B$ given by multiplication by two is the subgroup of two-torsion points of $E_{1,0}$.

Notice that we can reconstruct the divisors $\mathrm{tang}(\mathcal F, \mathcal L_i)$, for $i=2,3,4$, from the subgroups $S$ and $S^B$. Indeed
\[
\mathrm{tang}(\mathcal F,\mathcal  L_i) = \sum_{p \in S} L_p E_i + 2  \left( \sum_{p \in S^B}  L_p \widehat{E_i}  \right) \, .
\]
It follows that the foliation $\mathcal F$ is invariant by the natural action of $S^B\subset E_{1,0}$  in $T$, that is,
\begin{eqnarray*}
S^B \times T &\longrightarrow&  T  \\
(g,0), (x,y) &\mapsto& (x+g,y)  \, .
\end{eqnarray*}
Indeed, due  to the symmetry of our setup, $\mathcal F$ is left invariant by the following action of $(S^B )^2$ ,
\begin{eqnarray*}
(S^B)^2 \times T &\longrightarrow&  T  \\
\big((g,0),(h,0), (x,y)\big) &\mapsto& (x+g,y+h)  \, .
\end{eqnarray*}
The quotient of $\mathcal F\boxtimes \mathcal W$ by this action is
a CDQL $5$-web on $E_{\xi_3}^2$ of the form $\mathcal G \boxtimes \mathcal W$. If $E_{0,1}(2)$ denotes the two-torsion points on $E_{0,1}$ then, by construction,
\[
\mathrm{tang}(\mathcal G, \mathcal L_i) =  2 \, \widehat{E_i} + \sum_{p \in E_{0,1}(2)   } L_p E_i  \, .
\]
for $i=2,3,4$. This is sufficient to show that $\mathcal G\boxtimes \mathcal W $ is the $5$-web $\mathcal E_5$ of the Introduction.
\end{proof}

With Theorem \ref{T:E5} we complete the classification of exceptional CDQL webs on complex tori and, consequently, on compact complex
surfaces.

\end{document}